%% file: main.tex
\documentclass[11pt,twoside]{book}
\usepackage{project}
\graphicspath{{graphics/}}
\addto\captionsenglish{}

\def\projectauthor{Carola-Bibiane Sch\"onlieb and Zakhar Shumaylov}
\def\projecttitle{{\small CIME 2023 \\ Machine Learning: From Data to Mathematical Understanding}\\ Data-driven approaches to inverse problems}
\def\projectmonth{\today}
\def\projectdepartment{University of Cambridge}
\begin{document}
    \prefrontmatter
    \frontmatter
    \pagenumbering{roman} 
    \input{abstract}

        \markboth{Abstract}{}
    \pdfbookmark{\contentsname}{toc}
    \renewcommand{\sectionmark}[1]{\markright{#1}}
    \addtolength{\parskip}{-\baselineskip}  
    \tableofcontents
    \begingroup
        \let\clearpage\relax
        \vspace*{60pt}
        \listoffigures
        \vspace*{20pt}
    \endgroup
    \addtolength{\parskip}{\baselineskip}
    \renewcommand{\sectionmark}[1]{\markright{\thesection\ #1}}
    \mainmatter
    \input{chapter1}

    \input{chapter2}

    \input{chapter3}

    \input{chapter4}

    \cleardoublepage
    \renewcommand{\sectionmark}[1]{\markright{#1}}
    \addcontentsline{toc}{chapter}{Bibliography}
    \bibliographystyle{abbrvnat}
    \bibliography{refs}

\end{document}

%% file: abstract.tex
\chapter*{Abstract}

\noindent Inverse problems are concerned with the reconstruction of unknown physical quantities using indirect measurements and are fundamental across diverse fields such as medical imaging (MRI, CT), remote sensing (Radar), and material sciences (electron microscopy). These problems serve as critical tools for visualizing internal structures beyond what is visible to the naked eye, enabling quantification, diagnosis, prediction, and discovery. However, most inverse problems are ill-posed, necessitating robust mathematical treatment to yield meaningful solutions. While classical approaches provide mathematically rigorous and computationally stable solutions, they are constrained by the ability to accurately model solution properties and implement them efficiently. 

\noindent A more recent paradigm considers deriving solutions to inverse problems in a data-driven manner. Instead of relying on classical mathematical modeling, this approach utilizes highly over-parameterized models, typically deep neural networks, which are adapted to specific inverse problems using carefully selected training data. Current approaches that follow this new paradigm distinguish themselves through solution accuracy paired with computational efficiency that was previously inconceivable. 

\noindent These notes offer an introduction to this data-driven paradigm for inverse problems, covering methods such as data-driven variational models, plug-and-play approaches, learned iterative schemes (also known as learned unrolling), and learned post-processing. The first part of these notes will provide an introduction to inverse problems, discuss classical solution strategies, and present some applications. The second part will delve into modern data-driven approaches, with a particular focus on adversarial regularization and provably convergent linear plug-and-play denoisers. Throughout the presentation of these methodologies, their theoretical properties will be discussed, and numerical examples will be provided for image denoising, deconvolution, and computed tomography reconstruction. The lecture series will conclude with a discussion of open problems and future perspectives in the field.

%% file: chapter1.tex
\chapter{Introduction to Inverse Problems}
Inverse problems arise in a wide variety of scientific fields, from medical imaging and geophysics to finance and astronomy, where one is faced with the task of inferring information about an unknown object of interest from observed indirect measurements. Common features of inverse problems include the need to understand indirect measurements, and to overcome extreme sensitivity to noise and inaccuracies arising due to imperfect modeling. Knowledge-driven approaches traditionally dominated the field, relying on first-principles to derive physical models, leading to mathematically grounded reconstruction methods. However, the past decade has witnessed a paradigm shift towards data-driven methods, particularly with the advent of deep learning. While these data-driven approaches have achieved remarkable empirical success in image reconstruction, they often lack rigorous theoretical guarantees. In these notes we will examine the data-driven paradigm through a mathematical lens, presenting some state of the art methods, covering their provable properties and theoretical limitations.
\begin{figure}[htb] 
    \newlength{\commonimageheight}
    \setlength{\commonimageheight}{0.225\linewidth}
    \newlength{\commonimagewidth}
    \setlength{\commonimagewidth}{0.225\linewidth}
    \centering 

    \begin{subfigure}[b]{0.48\textwidth}
        \centering
        \includegraphics[height=\commonimageheight,width=\commonimagewidth]{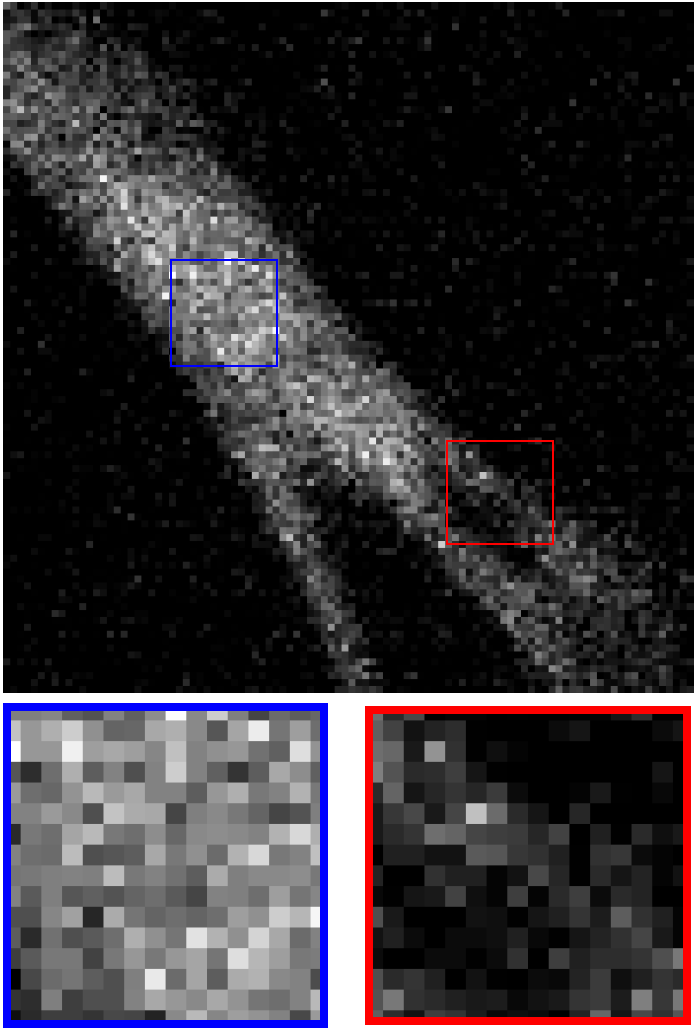}
        \hspace{0.025\linewidth} 
        \includegraphics[height=\commonimageheight,width=\commonimagewidth]{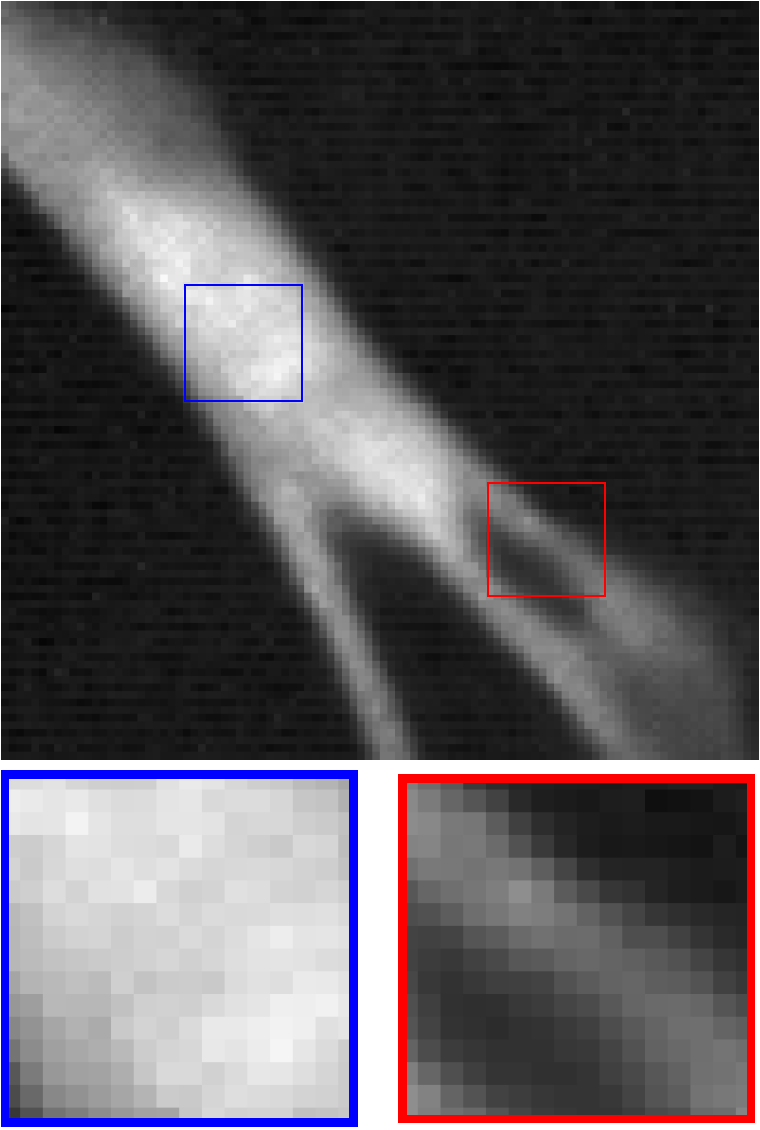}
        \caption{\textcolor{blue}{Image Denoising}: Given noisy image $y=u+n$, the task is to compute a denoised version $u^*$. \citep{ke2021unsupervised}}
        \label{fig:denoising}
    \end{subfigure}
    \hfill 
    \begin{subfigure}[b]{0.48\textwidth}
        \centering
        \includegraphics[height=\commonimageheight]{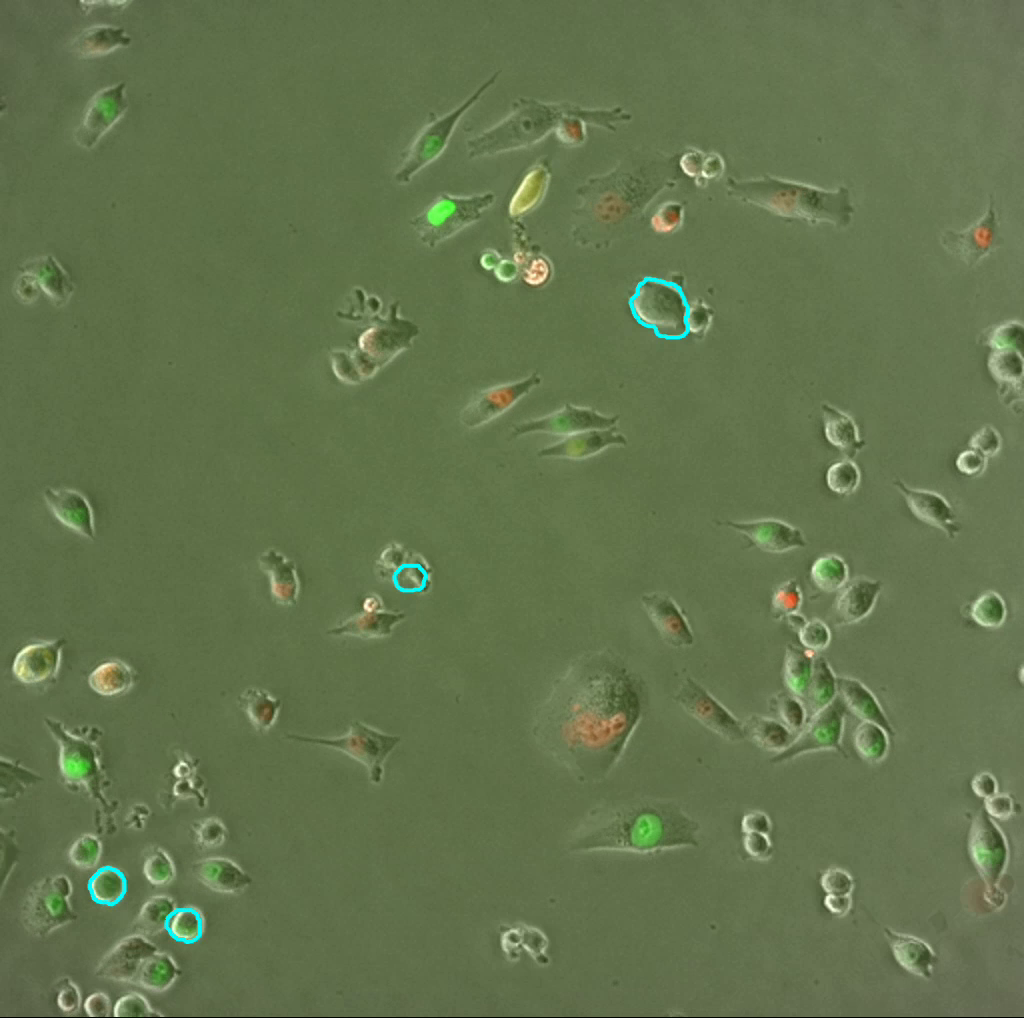} 
        \caption{\textcolor{blue}{Image Segmentation}: Given image $u$ on domain $\Omega$, compute the characteristic function $\chi_S$ of a region of interest $S\subset\Omega$ \citep{grah2017mathematical}\footnotemark.}
        \label{fig:segmentation}
    \end{subfigure}

    \vspace{1cm} 

    \begin{subfigure}[b]{0.48\textwidth}
        \centering
        \includegraphics[trim={0.1cm 0cm 0.1cm 0.1cm}, clip, height=0.45\linewidth,width=0.45\linewidth]{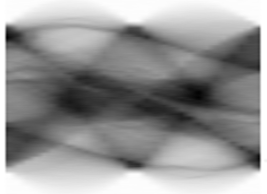}
        \hspace{0.025\linewidth} 
        \includegraphics[trim={2.5cm 6.25cm 15cm 6.5cm}, clip, height=\commonimageheight,width=\commonimagewidth]{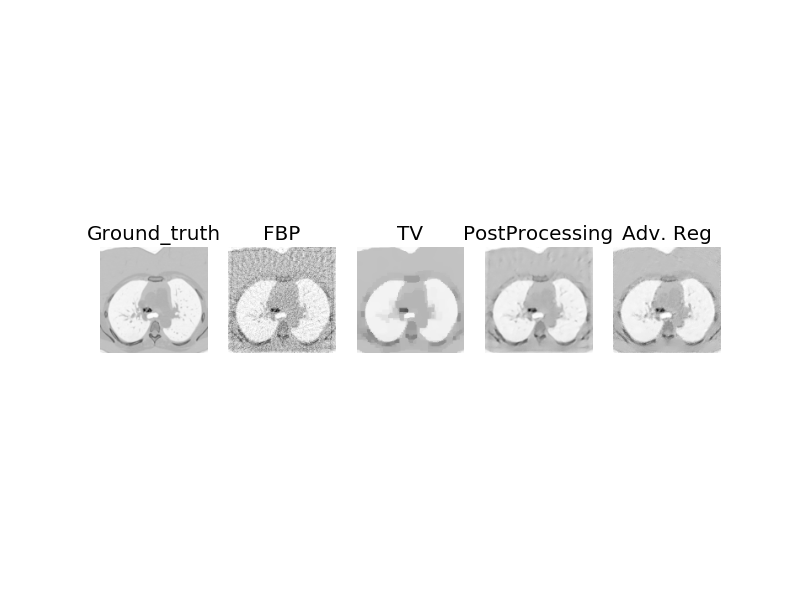}
        \caption{\textcolor{blue}{Image Reconstruction}: Compute $u$ from indirect, noisy measurements $y=A(u) + n$, with $A$ known operator  \citep{benning2018modern,arridge2019solving}.}
        \label{fig:reconstruction}
    \end{subfigure}
    \hfill 
    \begin{subfigure}[b]{0.48\textwidth}
        \centering
        \includegraphics[height=\commonimageheight,width=\commonimagewidth]{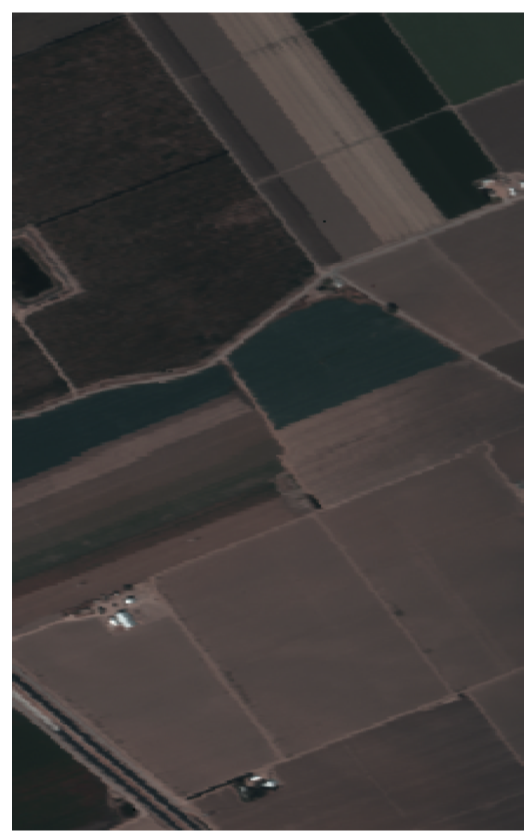} 
        \hspace{0.025\linewidth}
        \includegraphics[height=\commonimageheight,width=\commonimagewidth]{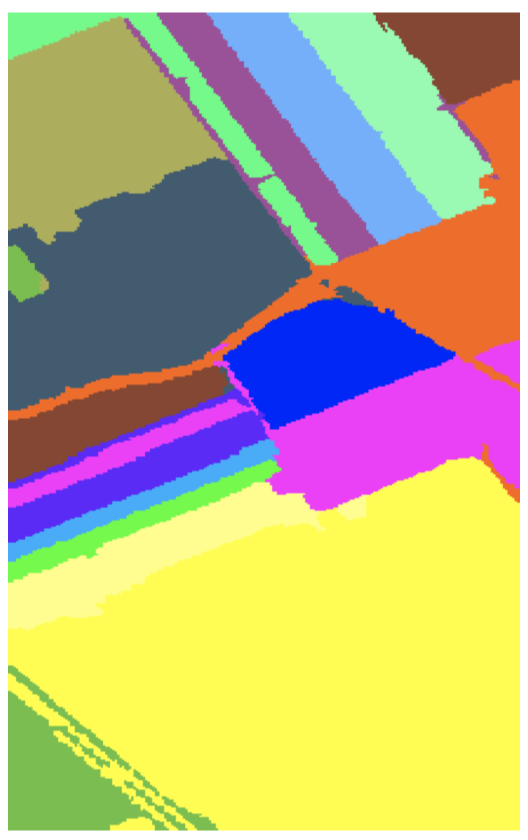} 
        \caption{\textcolor{blue}{Image Classification}: { Given a set of images $u_i$, the task is to assign appropriate labels $y_i$ to each image \citep{aviles2019graphx}.}} 
        \label{fig:classification}
    \end{subfigure}

    \caption{An overview of various fundamental image processing tasks.}
    \label{fig:image_tasks_overview}
    \footnotetext{\url{github.com/JoanaGrah/MitosisAnalyser}}
    \vspace{-1em}
\end{figure}

\noindent In this chapter, we begin by exploring the concept of well-posedness and its significance in the context of inverse problems. This will lead us to the notion of ill-posedness of inverse problems. These are often characterized by high sensitivity to noise, meaning even small errors or perturbations in the input data can lead to large variations in the solution. We will then investigate a range of knowledge-driven regularization techniques designed to mitigate the effects of ill-posedness and stabilize the solution process. 

\subsection{Examples of Inverse problems}

Before delving into the mathematical details, it's important to stress that the concepts discussed here are deeply rooted in real-world applications. Inverse problems form the fundamental language of numerous scientific and engineering disciplines, with some examples illustrated in \Cref{fig:biomedical_tasks_overview,fig:combined_apps_overview,fig:physical_sciences_overview}. To illustrate their complexity and wide-ranging applicability, let's consider several notable examples, mostly focusing on imaging problems, which often focus on four fundamental tasks illustrated in \Cref{fig:image_tasks_overview}.

\begin{itemize}
    \item \textit{Bio-Medical Imaging (Computed Tomography, MRI)}: Techniques such as Computed Tomography (CT) \citep{natterer2001mathematical} and Magnetic Resonance Imaging (MRI) \citep{lustig2008compressed,fessler2008image} serve as fundamental examples of inverse problems in action. These methods reconstruct detailed internal images of the human body from externally collected data, such as X-ray attenuation profiles or magnetic resonance signals, see e.g. \Cref{fig:reconstruction}. The challenge lies in accurately determining internal structures from these measurements, which are often incomplete or corrupted by noise due to various factors like patient movement or limitations on radiation exposure.
    Beyond these foundational examples, the scope of inverse problems in biomedical and biological imaging is vast, with various modalities and objectives, as illustrated in \Cref{fig:biomedical_tasks_overview}. This diversity includes other primary imaging modalities like Positron Emission Tomography (PET) (\Cref{fig:pet_imaging}), as well as more complex imaging scenarios such as spatio-temporal MRI, which incorporates additional dynamical information to capture changes over time (\Cref{fig:spatio_temporal_mri}). 
    Furthermore, solutions of inverse problems are often used for subsequent quantitative analysis and clinical decision-making. Examples of such downstream applications include the detailed analysis of cellular processes, such as mitosis analysis (\Cref{fig:mitosis_analysis}) and the estimation of cell dynamics (\Cref{fig:estimating_dynamics}), as well as clinical support tools like automated tumor segmentation (\Cref{fig:tumour_seg}) and systems for aiding diagnosis or prognosis from X-ray data (\Cref{fig:diag_prog}). The applications of inverse problem extend even further into broader biological research, including studies in zoology (\Cref{fig:zoology}) and molecular biology (\Cref{fig:molecular_biology}), underscoring their wide-ranging impact.

    \begin{figure}[p] 
    \centering 
    \newlength{\commonimageheighttwo}
    \setlength{\commonimageheighttwo}{0.23\linewidth}
    \newlength{\commonimagewidthtwo}
    \setlength{\commonimagewidthtwo}{0.225\linewidth}
    \begin{subfigure}[b]{0.48\textwidth} 
        \centering
        \includegraphics[height=\commonimageheighttwo]{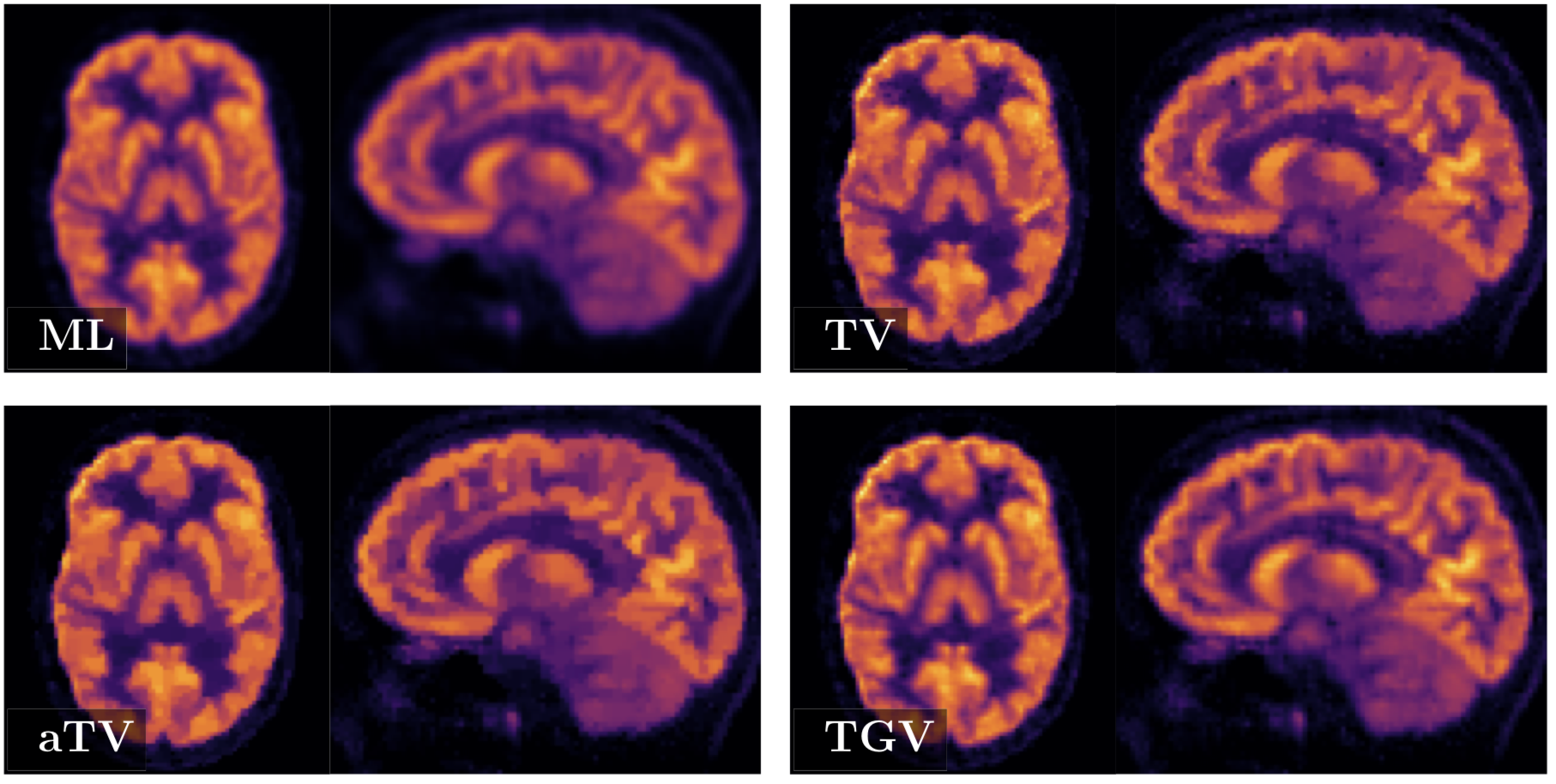} 
        \caption{PET imaging: \\\citep{ehrhardt2019faster,chambolle2018stochastic}} 
        \label{fig:pet_imaging}
    \end{subfigure}
    \hfill
    \begin{subfigure}[b]{0.48\textwidth}
        \centering
        \includegraphics[height=\commonimageheighttwo]{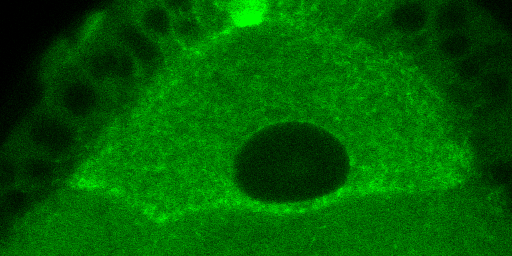} 
        \caption{Estimating cell dynamics:\\   \citep{drechsler2020optical}}
        \label{fig:estimating_dynamics}
    \end{subfigure}

    \vspace{1cm} 

    \begin{subfigure}[b]{0.24\textwidth}
        \centering
        \includegraphics[height=\commonimageheighttwo]{graphics/Fucci/ffmpeg_35} 
        \caption{Mitosis analysis:\\ \citep{grah2017mathematical}}
        \label{fig:mitosis_analysis}
    \end{subfigure}
    \hfill 
    \vspace{1cm}
    \begin{subfigure}[b]{0.48\textwidth}
        \centering
        \includegraphics[height=\commonimageheighttwo]{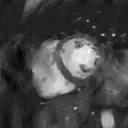} 
        \includegraphics[height=\commonimageheighttwo]{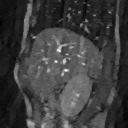} 
        \caption{Spatio-temporal MRI:  \citep{aviles2018csm,aviles2021compressed}}
        \label{fig:spatio_temporal_mri}
    \end{subfigure}
    \hfill
    \begin{subfigure}[b]{0.24\textwidth}
        \centering
        \includegraphics[height=\commonimageheighttwo]{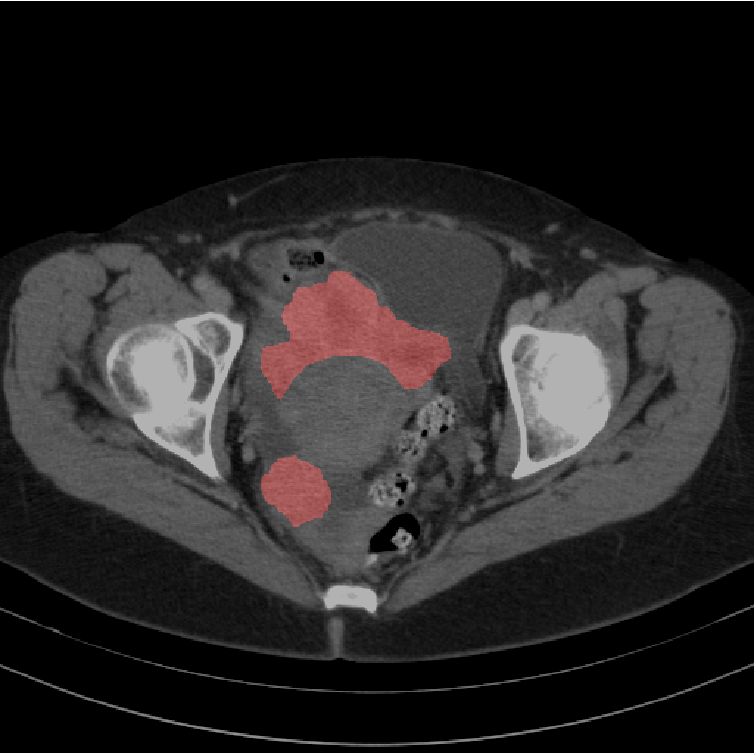} 
        \caption{Tumour segmentation: \\\citep{buddenkotte2023calibrating}} 
        \label{fig:tumour_seg}
    \end{subfigure}
    \vspace{1cm} 
    \begin{subfigure}[b]{.96\textwidth}
        \centering
        \includegraphics[height=\commonimageheighttwo]{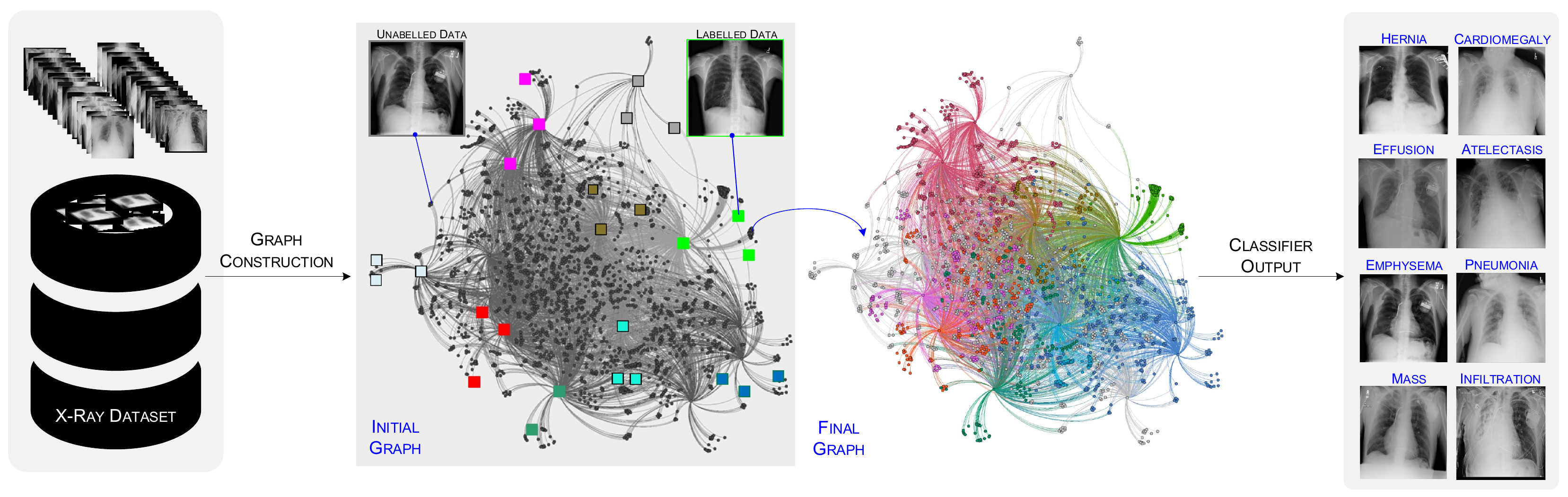} 
        \caption{X-Ray diagnosis and prognosis: \citep{aviles2019graphx}}
        \label{fig:diag_prog}
    \end{subfigure}
    \vspace{1cm}
    \begin{subfigure}[b]{0.48\textwidth} 
        \centering
        \includegraphics[height=\commonimageheighttwo]{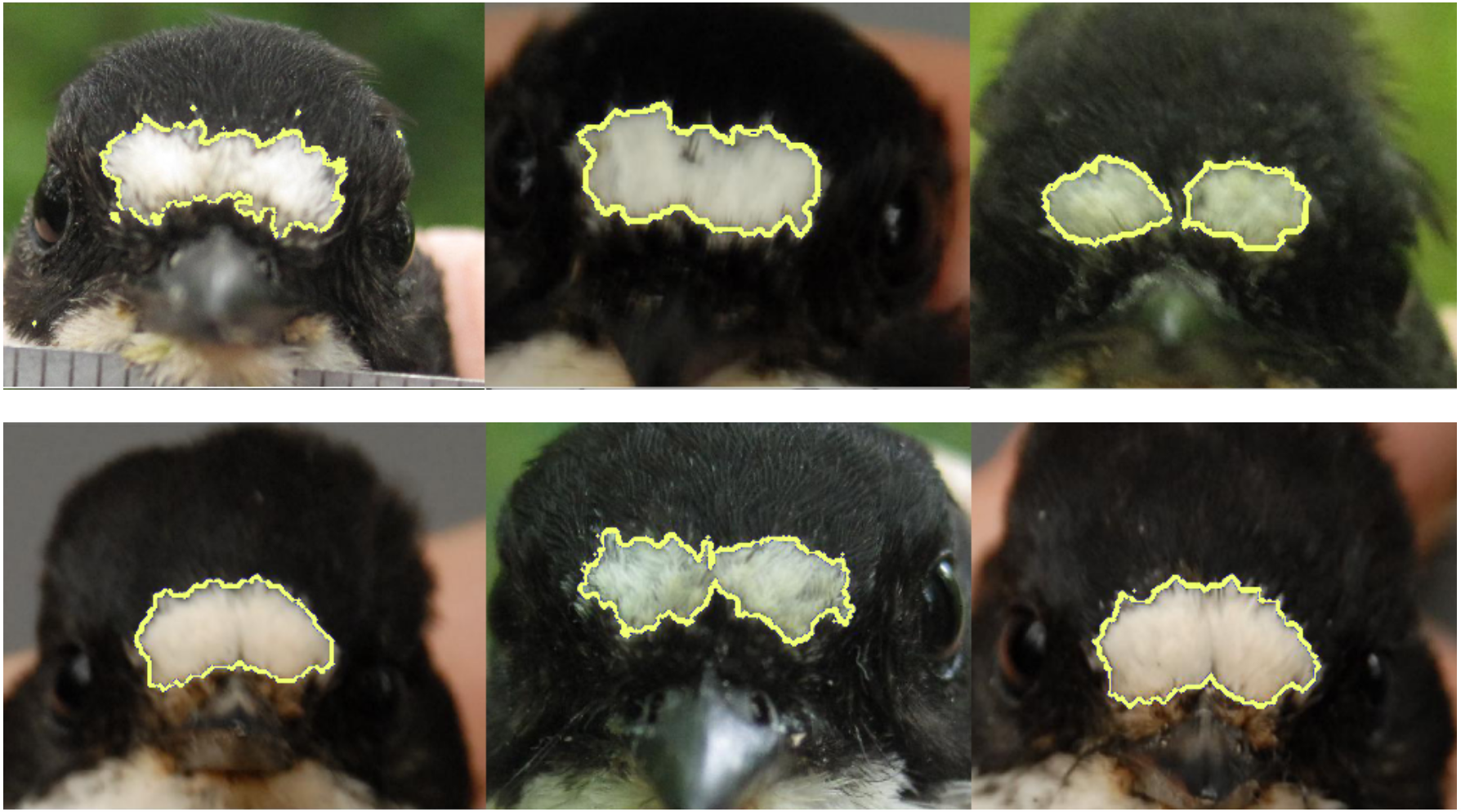} 
        \caption{Zoology: \\ \citep{calatroni2017graph}}
        \label{fig:zoology}
    \end{subfigure}
    \hfill
    \begin{subfigure}[b]{0.48\textwidth} 
        \centering
        \includegraphics[height=\commonimageheighttwo]{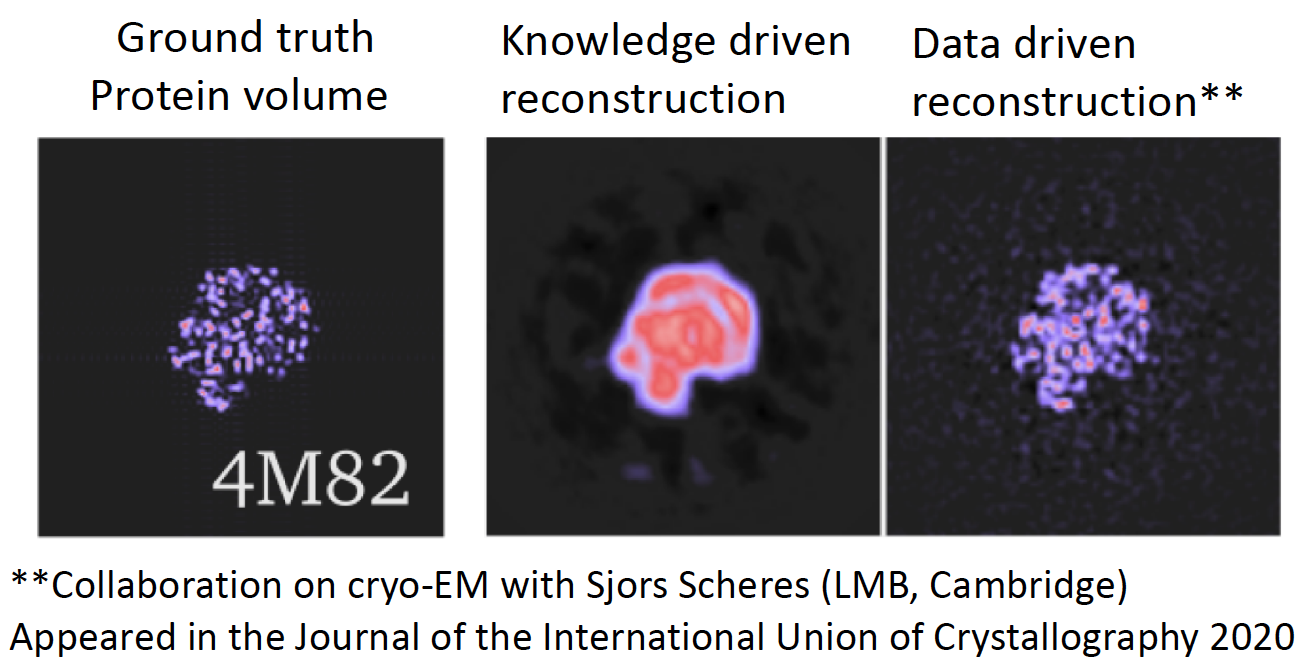} 
        \caption{Molecular biology: \\\citep{diepeveen2023regularizing,esteve2023spectral}}
        \label{fig:molecular_biology}
    \end{subfigure}
    \caption{Overview of various biological, biomedical and clinical research applications using image analysis.}
    \label{fig:biomedical_tasks_overview}
    \end{figure}
    

    \item \textit{General Image Processing}:
    Image processing more generally involves solving multiple intertwined inverse problems, with many applications spanning environmental conservation, remote sensing, an  digital humanities, as illustrated in \Cref{fig:combined_apps_overview}.
    For instance, in conservation and environmental science, LiDAR data can be used for detailed tree monitoring and forest assessment (\Cref{fig:tree_monitoring}), or multispectral and hyperspectral imagery for landcover analysis (\Cref{fig:landcover_analysis}). Remote sensing data, coupled with image processing, also aids in understanding dynamical systems, such as analyzing traffic flow for urban planning and infrastructure management (\Cref{fig:traffic_analysis}).
    Multi-modal image fusion (\Cref{fig:multimodal_fusion}), with data from different sensors, can be used for improving data representations, applicable in many fields ranging from remote sensing to medical imaging. In the realm of digital humanities, computational image processing plays a vital role for virtual art restoration and interpretation (\Cref{fig:virtual_art_restoration}), where imaging can help unveil hidden details, analyze materials, or digitally restore damaged cultural heritage artifacts.
    These varied applications all rely on extracting meaningful information from image data, often necessitating a range of image processing steps such as image reconstruction, enhancement, segmentation, feature extraction, deblurring, denoising and registration, many of which can be formulated as inverse problems.
    \begin{figure}[p]
    \centering
    \newlength{\commonimageheightthr}
    \setlength{\commonimageheightthr}{0.225\linewidth}
    \newlength{\commonimagewidththr}
    \setlength{\commonimagewidththr}{0.225\linewidth}
    \begin{subfigure}[b]{0.48\textwidth}
        \centering
        \includegraphics[height=\commonimageheightthr,width=2\commonimagewidththr]{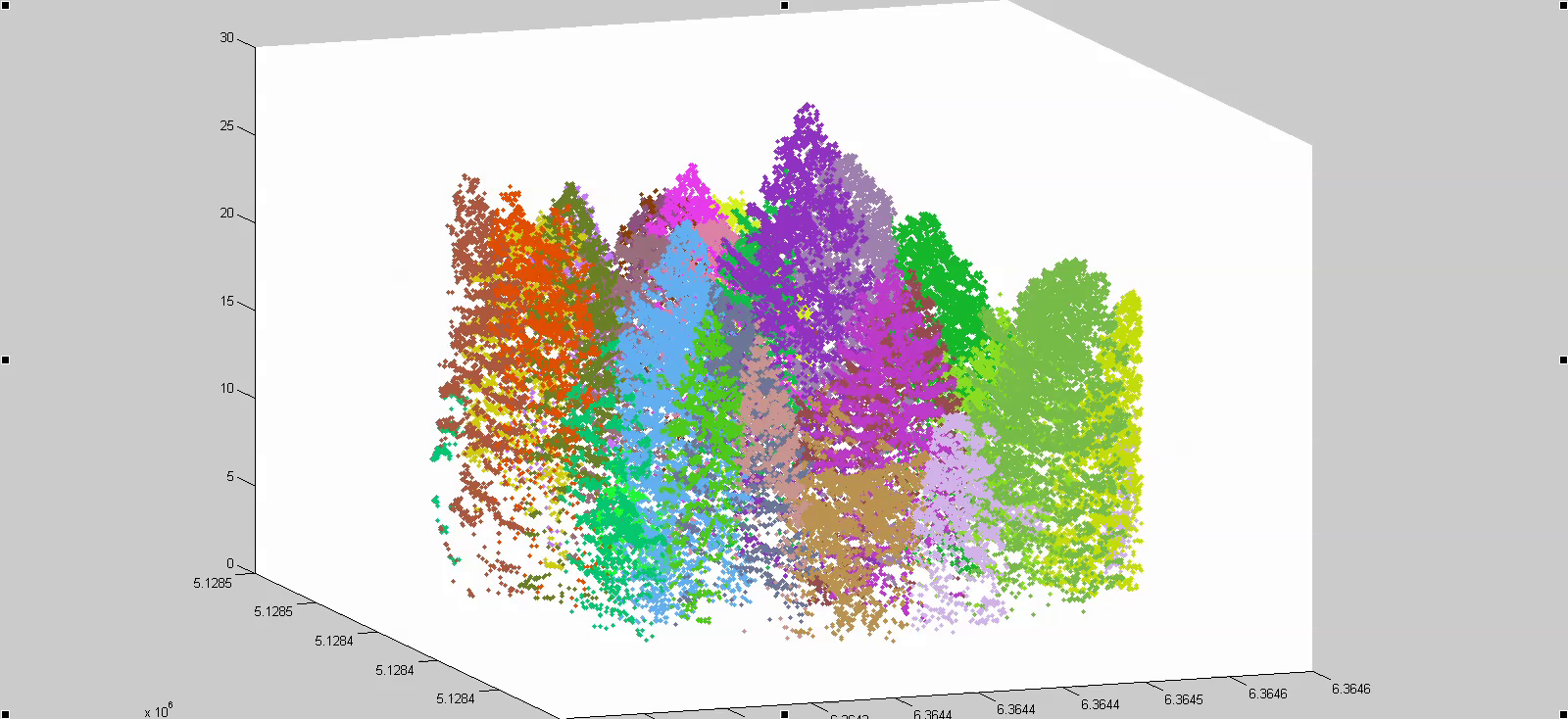} 
        \caption{Tree monitoring w/ LiDAR: \citep{lee2015nonparametric,lee2016individual}}
        \label{fig:tree_monitoring}
    \end{subfigure}
    \hfill
    \begin{subfigure}[b]{0.48\textwidth}
        \centering
        \includegraphics[height=\commonimageheightthr,width=\commonimagewidththr]{sellarshypsuppicture} 
        \hspace{0.005\linewidth}
        \includegraphics[height=\commonimageheightthr,width=\commonimagewidththr]{sellarshypsupclassification} 
        \caption{Landcover analysis: \citep{sellars2019semi}}
        \label{fig:landcover_analysis}
    \end{subfigure}

    \vspace{1cm} 

    \begin{subfigure}[b]{0.48\textwidth}
        \centering
        \includegraphics[height=\commonimageheightthr,width=2\commonimagewidththr]{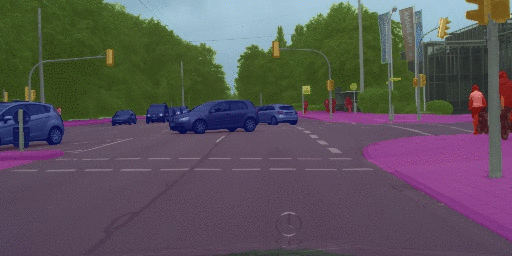} 
        \caption{Analysing traffic: \href{https://sites.google.com/view/remote-sens-research-for-india/home?authuser=0}{EPSRC project}}
        \label{fig:traffic_analysis}
    \end{subfigure}
    \hfill
    \begin{subfigure}[b]{0.48\textwidth}
        \centering
        \includegraphics[height=\commonimageheightthr,width=1.5\commonimagewidththr]{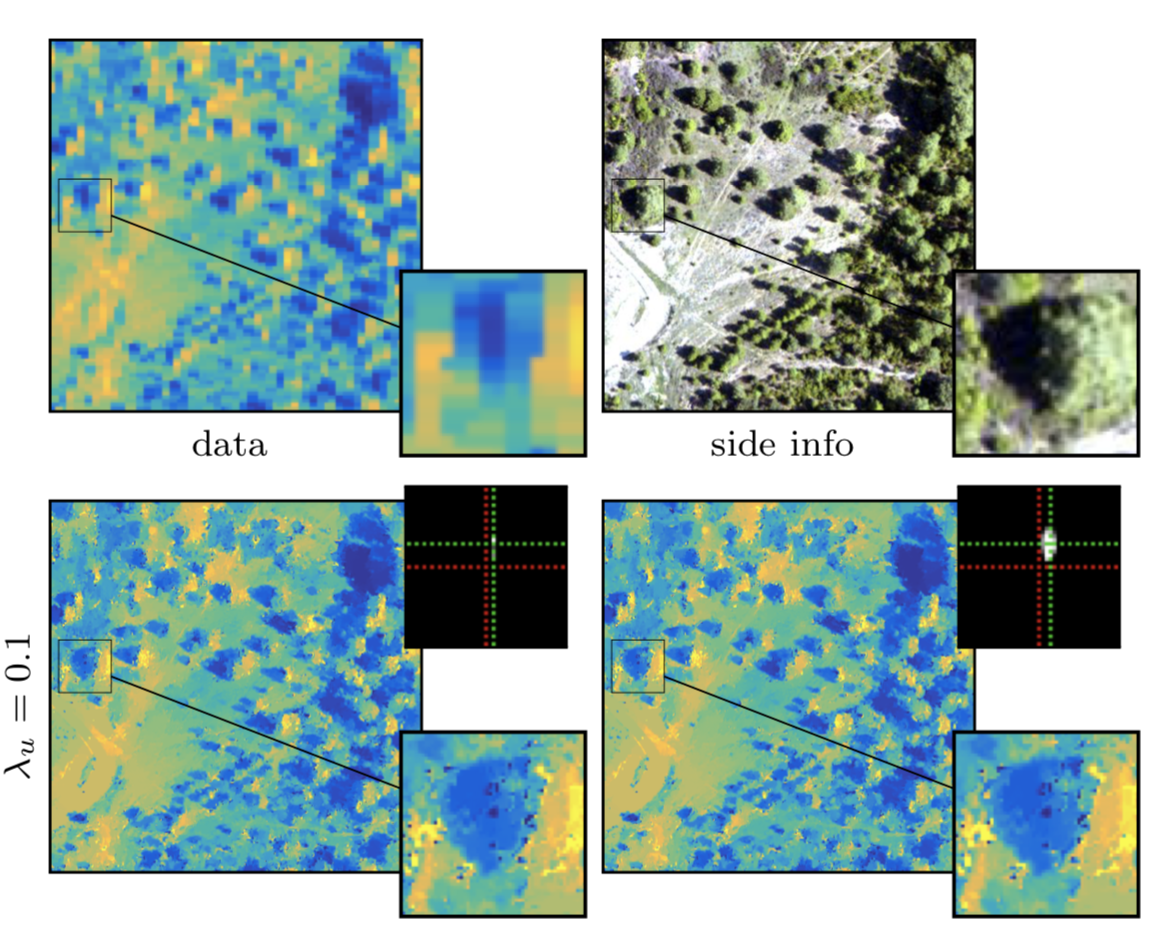} 
        \caption{Multi-modal image fusion: \citep{bungert2018blind}}
        \label{fig:multimodal_fusion}
    \end{subfigure}

    \vspace{1cm} 

    \begin{subfigure}[b]{0.96\textwidth} 
        \centering
        \includegraphics[height=1.3\commonimageheightthr,]{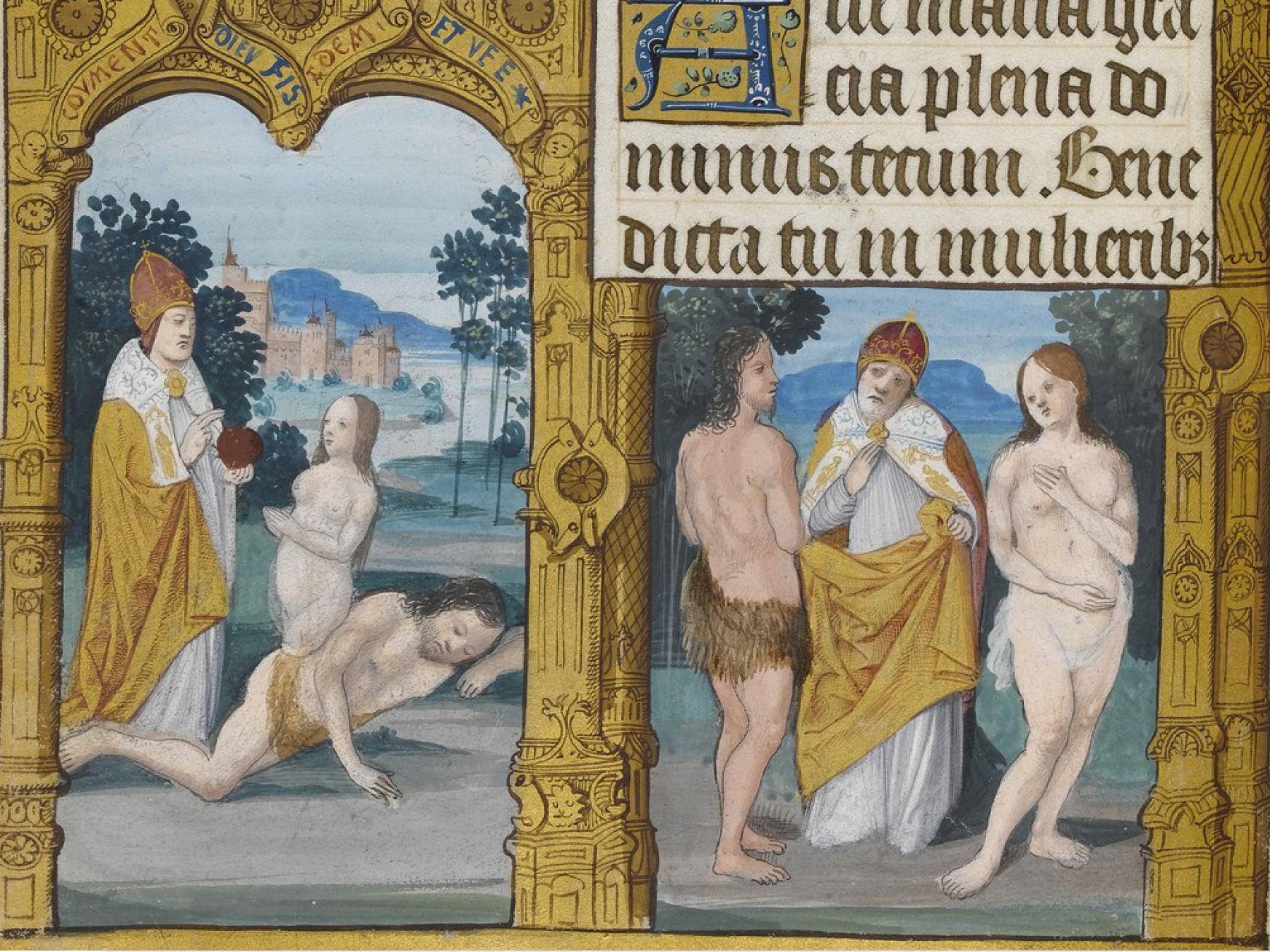}  
        \hfill
        \includegraphics[height=1.3\commonimageheightthr]{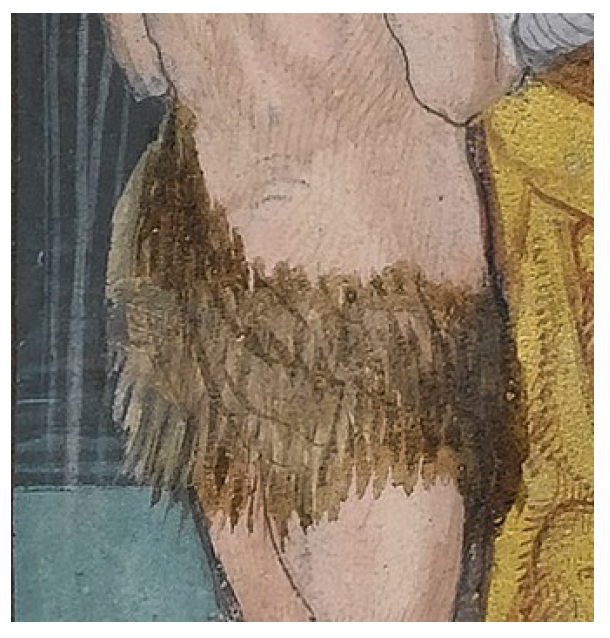} 
        \hfill
        \includegraphics[height=1.3\commonimageheightthr]{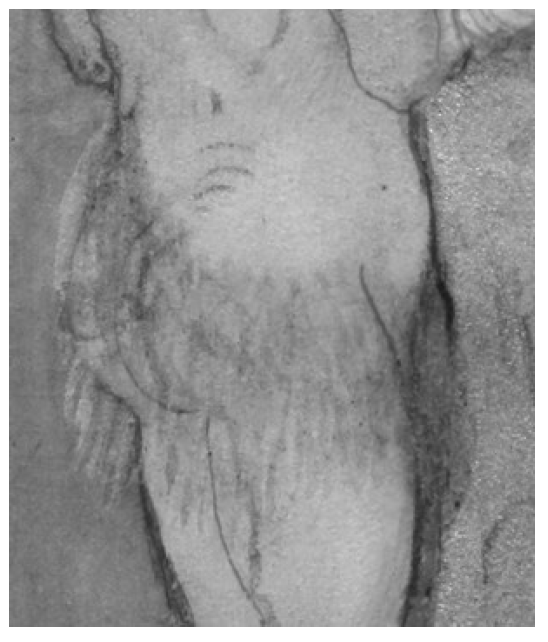} 
        \caption{Virtual art restoration and interpretation: Mathematics for Applications in Cultural Heritage (MACH) \url{https://mach.maths.cam.ac.uk} funded by the Leverhulme Trust; \citep{calatroni2018unveiling,parisotto2019anisotropic,parisotto2020variational}}
        \label{fig:virtual_art_restoration}
    \end{subfigure}
    
    \caption{Applications in conservation, sustainability, and digital humanities, showcasing various remote sensing and image analysis techniques.}
    \label{fig:combined_apps_overview}
    \end{figure} 

    \item \textit{Physical Sciences}:
    Inverse problems are foundational across the physical sciences, enabling researchers to probe and understand phenomena from the vastness of cosmic structures down to the intricacies of material microstructures, as illustrated in \Cref{fig:physical_sciences_overview}. In astrophysics, for instance, inverse problems arise when reconstructing images of distant celestial objects from data collected by telescopes \citep{starck1998image}. This process is complicated by vast distances, interference from various light sources, faint signals, and atmospheric disturbances. 
    A particularly famous example is the first imaging of a black hole (\Cref{fig:black_hole_imaging}) \citep{akiyama2019first}. 
    
    Similarly, geophysics utilizes seismic tomography (\Cref{fig:seismic_tomography}) to create images of the Earth's subsurface by analyzing seismic waves from earthquakes or controlled sources \citep{biegler2011large,biegler2003large,haber2014computational}. The inherent ill-posedness of this problem, due to complex geological layers, noisy and limited view data, necessitates regularization methods. The applications of inverse problem methodologies also extend to material sciences, where they are used to characterize material properties or analyze microstructures from indirect measurements (\Cref{fig:material_sciences}) \citep{tovey2019directional}, and to computational fluid dynamics, for estimating flow parameters or reconstructing complex flow fields from limited sensor data (\Cref{fig:cfd}) \citep{benning2014phase}. 
    
    \begin{figure}[p]
    \centering
    \newlength{\commonimageheightfour}
    \setlength{\commonimageheightfour}{0.225\linewidth}
    \begin{subfigure}[b]{0.48\textwidth}
        \centering
        \includegraphics[height=.92\commonimageheightfour]{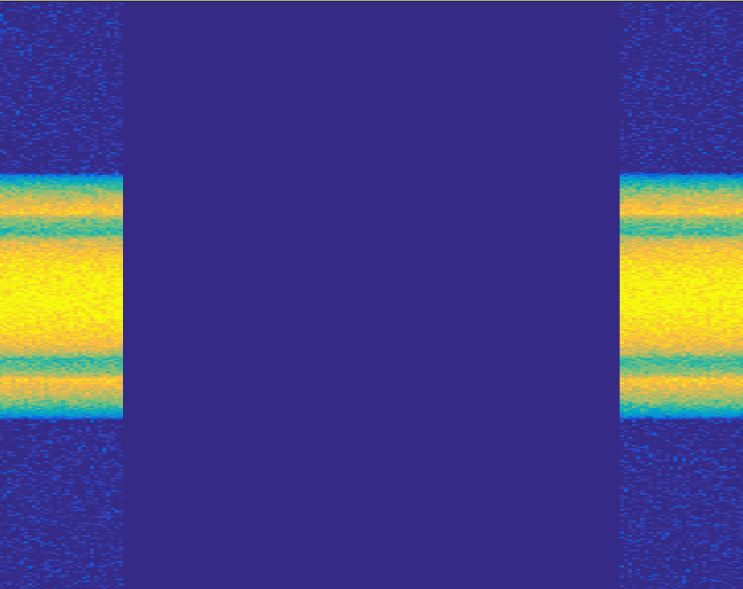} 
        \includegraphics[height=.92\commonimageheightfour]{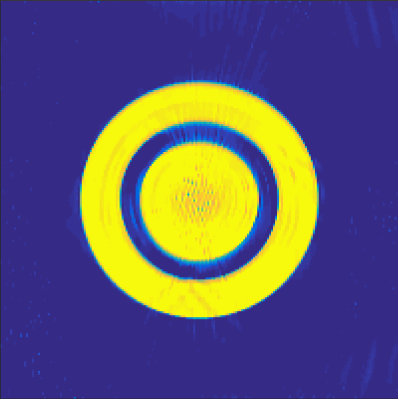} 
        \caption{Material Sciences:  \citep{tovey2019directional}}
        \label{fig:material_sciences}
    \end{subfigure}
    \hfill
    \begin{subfigure}[b]{0.48\textwidth}
        \centering
        \includegraphics[height=.92\commonimageheightfour]{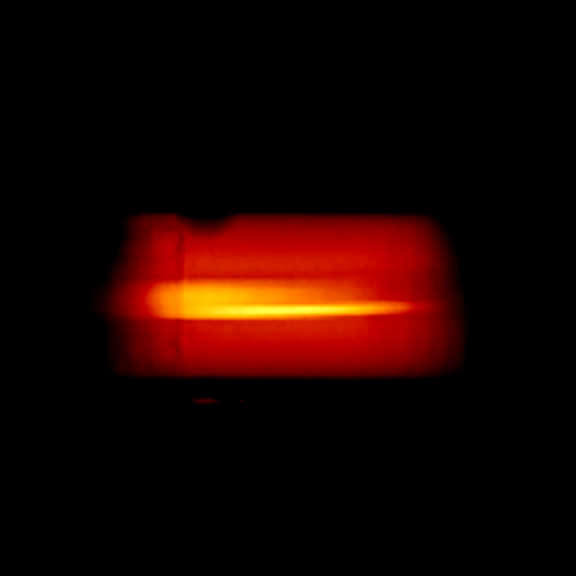} 
        \includegraphics[height=.92\commonimageheightfour]{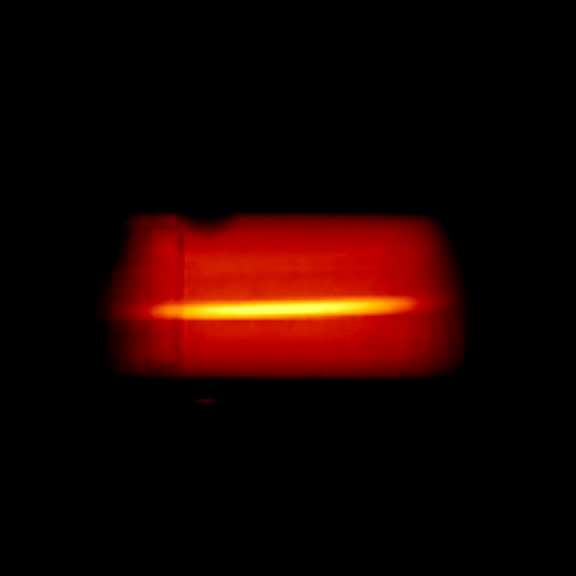} 
        \caption{CFD : \citep{benning2014phase}}
        \label{fig:cfd}
    \end{subfigure}
    \vspace{1cm} 
    \begin{subfigure}[b]{.98\textwidth}
        \centering
\includegraphics[height=2.5\commonimageheightfour]{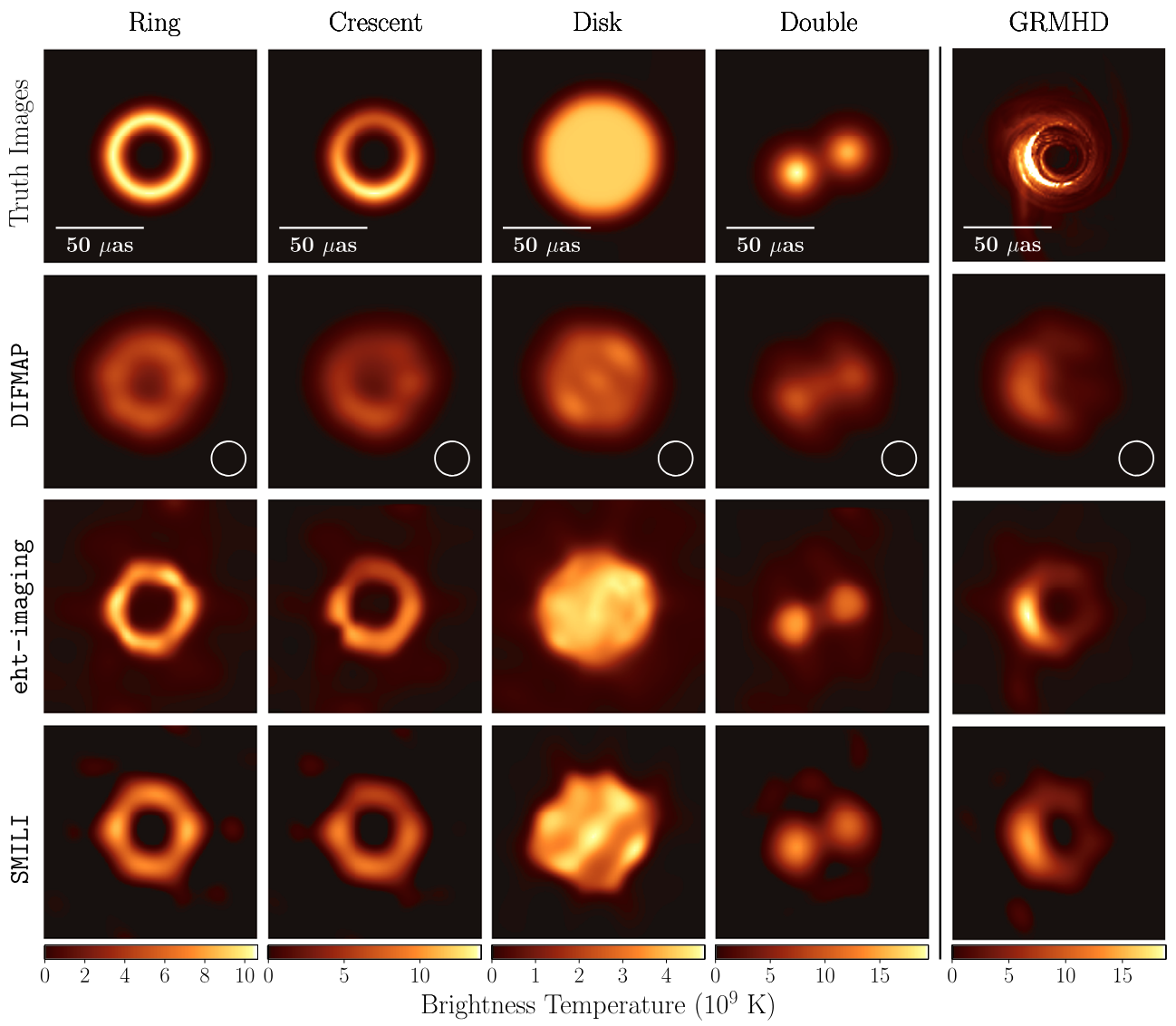} 
        \caption{Black Hole Imaging \citep{akiyama2019first}} 
        \label{fig:black_hole_imaging}
    \end{subfigure}
    \vspace{1cm}
    \begin{subfigure}[b]{.96\textwidth}
        \centering
    \includegraphics[height=1.75\commonimageheightfour]{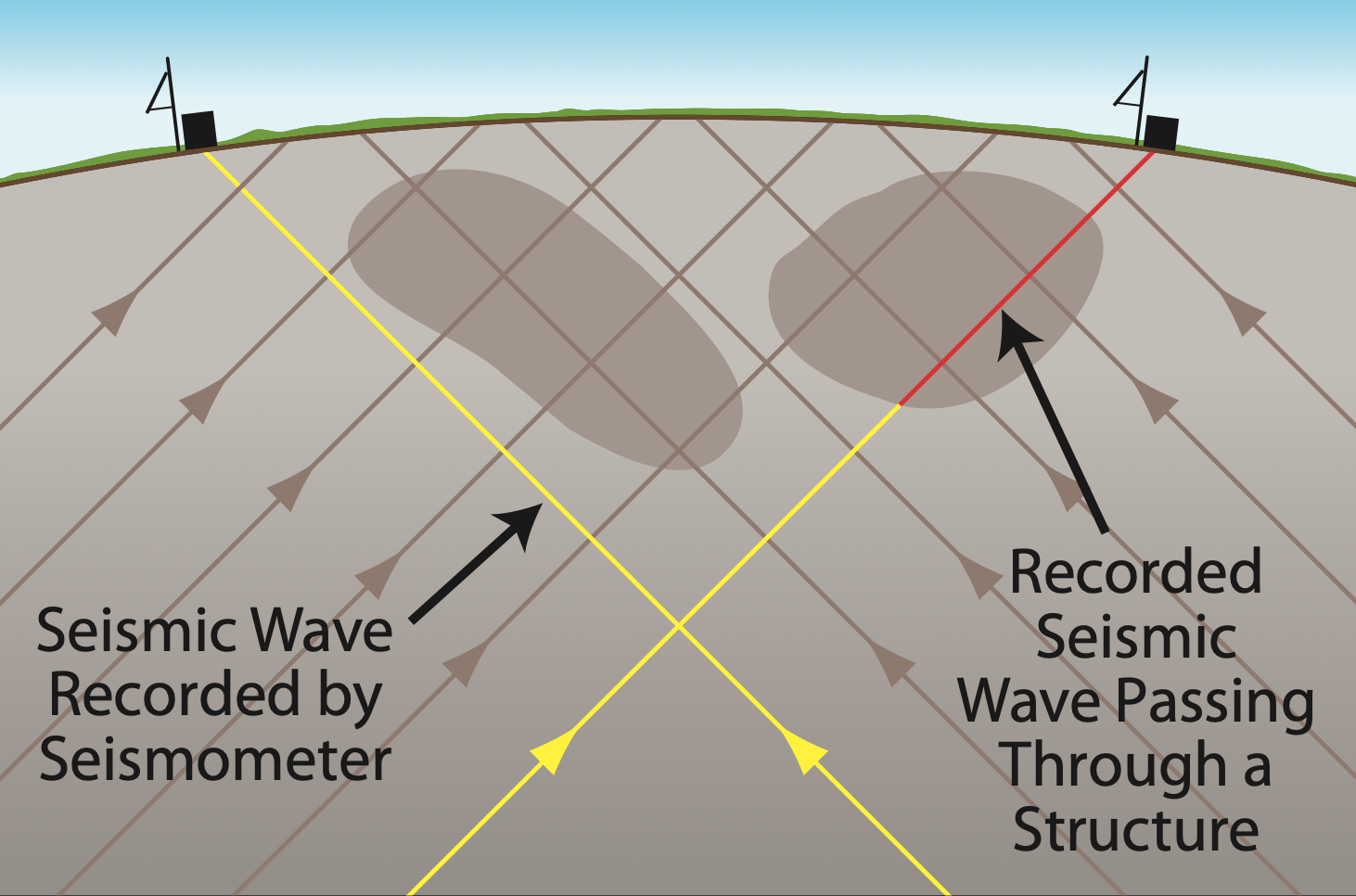} 
        \caption{Seismic Tomography: \href{https://www.iris.edu/hq/inclass/fact-sheet/seismic_tomography}{EarthScope Consortium}}
        \label{fig:seismic_tomography}
    \end{subfigure}
    \caption{Applications in Physical Sciences, including materials science, computational fluid dynamics, astrophysics, and geophysics.}
    \label{fig:physical_sciences_overview}
\end{figure}
\end{itemize}

\clearpage
\section{Well-posed and ill-posed problems}
An \textit{inverse problem} is the problem of finding $u$ satisfying the equation
$$
y=A u,
$$
where $y \in \mathbb{R}^m$ is given. We refer to $y$ as \textit{observed data} or \textit{measurement}, and $u$ as an unknown. The physical phenomena that relates the unknown and the measurement is modeled by a linear operator $A$. The function spaces for $u$ and $A$ will be specified later when particular inverse problems or reconstruction approaches are being discussed. In reality, the perfect data is perturbed by noise. For an additive noise model, we observe measurements of the form
\begin{equation}\label{eq:first_chap_ip}
    y_n=A u+n,
\end{equation}
where $n \in \mathbb{R}^m$ represents the observational noise. The main question that we concern ourselves with is: ``Can we always compute a reliable answer $u$?'' We are often interested in ill-posed inverse problems, where the inverse problem is more difficult to solve than the forward problem of finding $y_n$ when $u$ is given. To explain what we mean by ill-posed, we first need to introduce well-posedness as defined by \citet{hadamard1902problemes}:
\begin{definition}[Well-posed problem]\label{def:wellposed}
A problem is called \textit{well-posed} if:
\begin{itemize}[nosep]
    \item There exists at least one solution. (Existence)
    \item There is at most one solution. (Uniqueness)
    \item The solution depends continuously on data. (Stability)
\end{itemize}
\end{definition}
\noindent While the forward problem is generally assumed to be well-posed, inverse problems are typically \textit{ill-posed}, meaning they violate at least one of the conditions for well-posedness. This focus on ill-posedness arises from physical motivations as opposed to as an abstract concern; it stems directly from the fact that a vast majority of practical problems in science and engineering are indeed ill-posed. The following simple examples illustrate common issues that arise when these conditions are not met. 
\begin{example}
 \begin{itemize}
    \item (\textit{Non-existence}) Assume that $n<m$ and $A: \mathbb{R}^n \rightarrow \operatorname{ran}(A) \subsetneq \mathbb{R}^m$, where the range of $A$ is a proper subset of $\mathbb{R}^m$. Because of the noise in the measurement, it may end up that $y_n \notin \operatorname{ran}(A)$ so that inverting $A$ is not possible. Note that usually only the statistical properties of the noise $n$ are known, so we can not simply subtract it.
    \item (\textit{Non-uniqueness}) Assume next that $n>m$ and $A: \mathbb{R}^d \rightarrow \mathbb{R}^m$, in which case the system is under-determined. We then have more unknowns than equations which means that there may be infinitely many solutions.
    \item (\textit{Instability}) Finally suppose that $n=m$, and that there exists an inverse $A^{-1}: \mathbb{R}^m \rightarrow \mathbb{R}^n$. Suppose further that the condition number $\kappa=\lambda_1 / \lambda_m$ is very large, where $\lambda_1$ and $\lambda_m$ are the biggest and smallest eigenvalues of $A$. Such a matrix is said to be \textit{ill-conditioned}. In this case, the problem is sensitive to even the smallest errors in the measurement. The naive reconstruction $\tilde{u}=A^{-1} y_n=u+A^{-1} n$ would not produce a meaningful solution, but would instead be dominated by the noise component $A^{-1} n$. 
\end{itemize}
The last example illustrates one of the key questions of inverse problem theory: \begin{center}\textbf{How to stabilize the reconstruction while maintaining accuracy?} \end{center}
\end{example}



\begin{example}[Blurring in the continuum]\label{eg:heat}
One common and illustrative inverse problem encountered in image processing is deblurring. Imagine we have an image that has been blurred, perhaps by camera motion or an out-of-focus lens. Our goal is to recover the original, sharp image. This seemingly straightforward task quickly reveals the challenges inherent in many inverse problems.

\noindent Let us consider this in a continuous one-dimensional setting.
Suppose our observed blurred function, $y(x) : \mathbb{R} \rightarrow \mathbb{R}$, results from the convolution of an original sharp function, $u(x)$, with a blurring kernel, i.e. $y(x)=\left(G_\sigma * u\right)(x)$, where
$$
G_\sigma(x):=\frac{1}{2 \pi \sigma^2} e^{-|x|^2 /\left(2 \sigma^2\right)} \quad=\text { Gaussian kernel }
$$
with standard deviation $\sigma$, dictating the extent of the Gaussian blur. Our objective is to reconstruct the original sharp function $u$ from the observed blurred function $y$. This turns out to be equivalent to inverting the heat equation. To be precise, the blurred measurement $y(x)$ can be seen as a solution to the heat equation at a specific time $t=\sigma^2/2$, with $u(x)$ being the initial condition. Therefore, attempting to retrieve $u(x)$ from $y(x)$ is analogous to solving the heat equation backward in time. Ill-posedness in this example arises from a lack of continuous dependence of the solution on the data: small errors in the measurement $y$ can lead to very large errors in the reconstructed $u$. From Fourier theory, we can write the following, where $\mathcal{F}$ and $\mathcal{F}^{-1}$ denote the Fourier transform and inverse Fourier transform respectively:
$$
\begin{aligned}
& y=\sqrt{2 \pi} \mathcal{F}^{-1}\left(\mathcal{F} G_\sigma \mathcal{F} u\right) \\
& u=\frac{1}{\sqrt{2 \pi}} \mathcal{F}^{-1} \frac{\mathcal{F} y}{\mathcal{F} G_\sigma} .
\end{aligned}
$$
Instead of measuring a blurry $y$, suppose now that we measure a blurry and noisy $y_\delta=y+n^\delta$, with deblurred solution $u_\delta$. Then,
$$
\sqrt{2 \pi}\left|u-u_\delta\right|=\left|\mathcal{F}^{-1} \frac{\mathcal{F}\left(y-y_\delta\right)}{\mathcal{F} G_\sigma}\right|=\left|\mathcal{F}^{-1} \frac{\mathcal{F} n_\delta}{\mathcal{F} G_\sigma}\right|
$$
Now, for high-frequencies, $\mathcal{F}\left(n_\delta\right)$ will be large while $\mathcal{F} G_\sigma$ will be small (since $G_\sigma$ is a compact operator). Hence, the high frequency components in the error are amplified!
\end{example}
\noindent In many practical linear inverse problems, the condition of \emph{continuity} is the first to break down. This failure of continuous dependence of the solution on the data leads to extreme amplification of noise. In addition, the \emph{uniqueness} condition often fails in under-sampled inverse problems. While under-sampling has physical advantages, such as faster data acquisition or reduced exposure (e.g., in medical imaging), the trade-off is significant noise amplification, thereby making the ill-posedness even more severe.
\begin{example}[Computed Tomography]
In almost any tomography application, the underlying inverse problem is either the inversion of the Radon transform or of the X-ray transform. Here, we primarily follow \citet{sherry2025part}. 
For $u \in C_0^{\infty}\left(\mathbb{R}^2\right), s \in \mathbb{R}$, and $\theta \in \mathbb{S}^{1}$ the Radon transform $R: C_0^{\infty}\left(\mathbb{R}^2\right) \rightarrow C^{\infty}\left(\mathbb{S}^{1} \times \mathbb{R}\right)$ can be defined as the integral operator
$$
\begin{aligned}
f(\theta, s)=(R u)(\theta, s) & =\int_{x^{\prime} \cdot \theta=s} u(x) \mathrm{d} x =\int_{\theta^{\perp}} u(s \theta+y) \mathrm{d} y
\end{aligned}
$$
for $\theta \in \mathbb{S}^{1}$ and $\theta^{\perp}$ being the vector orthogonal to $\theta$. Effectively the Radon transform in two dimensions integrates the function $u$ over lines in $\mathbb{R}^2$. Since $S^{1}$ is the unit circle $S^1=\left\{\theta \in \mathbb{R}^2 \mid\|\theta\|=1\right\}$, we can choose for instance $\theta=(\cos (\varphi), \sin (\varphi))^{\top}$, for $\varphi \in[0,2 \pi)$, and parameterize the Radon transform in terms of $\varphi$ and $s$, i.e.
$$
f(\varphi, s)=(R u)(\varphi, s)=\int_{\mathbb{R}} u(s \cos (\varphi)-t \sin (\varphi), s \sin (\varphi)+t \cos (\varphi)) \mathrm{d} t
$$
Note that with respect to the origin of the reference coordinate system, $\varphi$ determines the angle of the line along one wants to integrate, while $s$ is the offset from that line from the center of the coordinate system. It can be shown that the Radon transform is linear and continuous, and even compact.
Visually, CT simply turns images into sinograms:
\begin{figure}[h!] 
    \centering 
    \vspace{-0.5em}
    \begin{equation*}
        R  \vcenter{\hbox{\includegraphics[width=2.7cm,height=2.43cm]{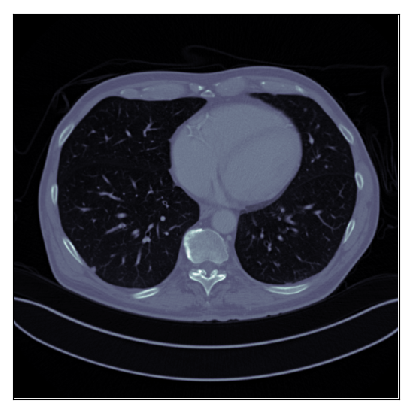}}} 
        + n =
\vcenter{\hbox{\includegraphics[width=2.7cm,height=2.565cm]{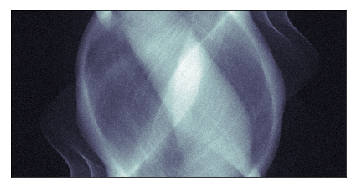}}}
    \end{equation*}
    \vspace{-0.5em}
    \label{fig:visual_equation_t_u_plus_n_equals_y}
\end{figure}
\paragraph{Radon Transform Inversion is Ill-Posed}
While enabling many applications, in practice inversion of the Radon transform can be shown to be ill-posed \citep{hertle1981problem}, by considering singular values of the operator to show unboundedness. A practical illustration of this ill-posedness is shown in \Cref{fig:illposedct}.

\noindent Informally, following \citep{epstein2007introduction,Candes_Math262_Lec10}, we can construct eigenfunctions of the operator $R^*R$, where $R^*$ is the adjoint of the Radon transform. For a formal presentation see \citep{Candes_Math262_Lec11}. For $g(x)=e^{i\left\langle{k}, {x}\right\rangle}$, we have:
$$
{R}^* {R}\left[g\right]({x})=\frac{1}{\left\|{k}\right\|} e^{i\left\langle{k}, {x}\right\rangle}
$$
Thus, $g$ is an eigenfunction of ${R}^* {R}$ with eigenvalue $\frac{1}{\left\|{k}\right\|}$, meaning that singular values of ${R}$ are $\frac{1}{\sqrt{\|{k}\|}}$. As these tend to $0$, the inverse is unbounded and the problem is ill-posed.
\begin{figure}[htb!]
    \centering
    \begin{minipage}{0.45\textwidth}
        \centering
        \includegraphics[width=0.9\textwidth]{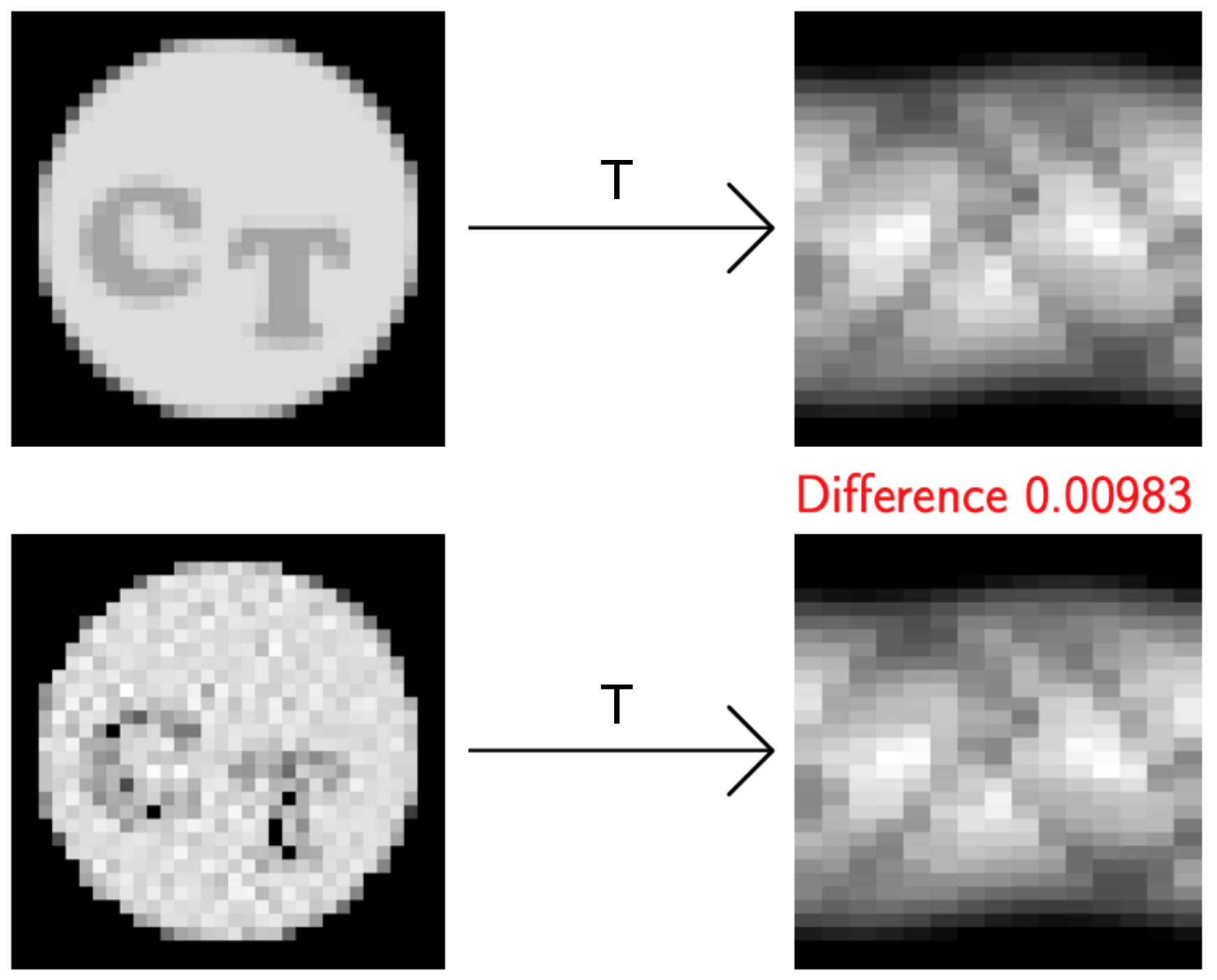} 
    \end{minipage}
    \hspace{0.05\textwidth} 
    \begin{minipage}{0.45\textwidth}
        \small 
        \emph{Ill-posedness:} $A^{-1}$ is not continuously invertible (unbounded or discontinuous)\\[0.5cm]
        \emph{Typical reasons:} noise, undersampling, nonlinearity, \ldots
    \end{minipage}
    \caption{Illustration of non-uniqueness in CT reconstruction (left) and a description of ill-posedness in inverse problems (right). Courtesy of Samuli Siltanen.}
    \label{fig:illposedct}
\end{figure}

\end{example}
\section{Overcoming the ill-posedness}
The primary strategy to overcome the ill-posedness is not to solve the original ill-posed problem directly, but rather to formulate and solve a related well-posed problem whose solution is a good approximation of the true, underlying solution we seek. This is generally called \textit{regularization}. 
\subsection{Regularization in functional analysis} 
More generally considering a linear bounded operator $A\in\mathcal{L}(\mathcal{U}, \mathcal{V})$ between two normed spaces $\mathcal{U}$ and $\mathcal{V}$, the natural choice for the inverse is the Moore-Penrose inverse $A^{\dagger}$ (see e.g. \citep{sherry2025part}), however when $\operatorname{im}(A)$ is not closed, the pseudo-inverse is not bounded: given noisy data $y^\delta$ such that $\left\|y^\delta-y\right\| \leqslant \delta$, we cannot expect convergence $A^{\dagger} y^\delta \rightarrow A^{\dagger} y$ as $\delta \rightarrow 0$. To achieve convergence, $A^{\dagger}$ is replaced with a family of well-posed reconstruction operators operators $\mathcal{R}_\alpha$ with $\alpha=\alpha\left(\delta, y^\delta\right)$ and require that $\mathcal{R}_{\alpha\left(\delta, y^\delta\right)}(y^\delta) \rightarrow A^{\dagger} y$ for all $y \in \operatorname{dom}(A^{\dagger})$ and all $y^\delta \in \mathcal{V}$ s.t. $\left\|y-y^\delta\right\| \leqslant\delta$ as $\delta \rightarrow 0$. 
\begin{definition}\label{def:funan_conv}
Let $A \in \mathcal{L}(\mathcal{U}, \mathcal{V})$ be a bounded operator. A family $\left\{\mathcal{R}_\alpha\right\}_{\alpha>0}$ of continuous operators is called a regularization (or regularization operator) of $A^{\dagger}$ if
$$
\mathcal{R}_\alpha y \rightarrow A^{\dagger} y=u^{\dagger}
$$
for all $f \in \operatorname{dom}(A^{\dagger})$ as $\alpha \rightarrow 0$. In other words, a regularization is a pointwise approximation of the Moore-Penrose inverse with continuous operators.
\end{definition}
\subsection{Variational Regularization}
An influential technique to overcome the ill-posedness of \Cref{eq:first_chap_ip}, ensuring that reconstructions are regularization operators, has been the variational approach pioneered by \citet{Tikh1963} and \citet{phillips1962}. Here, the inverse problem is re-framed as an optimization task of minimizing a carefully constructed ``\textit{variational energy functional}''. The energy functional typically has two components, taking the form:
$$
\min _{u \in U}\left\|A u-y_n\right\|^2+\alpha R(u).
$$
Here, $\alpha>0$ acts as a tuning parameter balancing the effect of the \textit{data fidelity term} $\left\|A u-y_n\right\|^2$, which ensures consistency with the observed measurements. The \textit{regularization term} $R(u)$ aims to incorporate prior knowledge about the reconstruction, penalizing $x$ if it is not ``realistic''. 

The selection of an appropriate regularizer $R$ and the tuning of the parameter $\alpha$ are critical for designing effective regularization methods. By introducing regularization, we aim to achieve more than just any solution: we seek a problem that is \textit{well-posed} and whose solution is a \textit{good approximation of the true solution}. This is normally captured through the properties of \emph{existence}, \emph{uniqueness}, \emph{stability} and \emph{convergence}. Here, we will present this for the simple setting of finite dimensions and strongly convex regularizers, and a proof for this setting can be found in \citet{mukherjee2024data}. 
\begin{theorem}\label{thm:stab}
For a forward operator $A:\mathbb{R}^n\to\mathbb{R}^m$ and noisy measurement $y^\delta=Au^*+n,$ assume that the noise is bounded $\|n\| < \delta$, and let $R$ be a strongly convex regularization term. Let 
\begin{equation}\label{eq:ip_conv}
    u(y^{\delta},\alpha) = \underset{u \in \mathbb{R}^n}{\operatorname{argmin}}\left\|A u-y^{\delta}\right\|^2+\alpha R(u),
\end{equation}
then the following hold:
\begin{itemize}[nosep]\vspace{-0.5em}
    \item \textbf{\upshape Existence and Uniqueness:} 
    For all $\delta \geq 0,\, \alpha > 0,\, y\in\mathbb{R}^m$, the minimizer of \Cref{eq:ip_conv} exists and is unique.
    \item \textbf{\upshape Stability:} For all sequences $y_k \to y, \ u(y_k,\alpha)\to u(y,\alpha)$, i.e. $u$ is continuous in $y$.
    \item \textbf{\upshape Convergent Regularization:} For a certain parameter rule $\alpha(\delta)$ satisfying $(\delta,\alpha(\delta))\to0$, we have $u(y^\delta,\alpha(\delta))\to u^\dagger$, where $$u^\dagger = \underset{\substack{u\in{\mathbb{R}^n} \\ \text{ \upshape s.t. } Au = Ay^0}}{\operatorname{argmin}}\; R(u)$$ is the $R$-minimizing solution, generalizing \Cref{def:funan_conv} to general inverses.
\end{itemize}
\end{theorem}
\noindent This naturally extends to more complex settings \citep{shumaylov2024weakly, poschl2008tikhonov, scherzer2009variational, subho2023}. In words, \emph{Existence}, \emph{Uniqueness} and \emph{Stability} ensure that \Cref{def:wellposed} is satisfied and the regularized inverse problem is well-posed. \emph{Convergent Regularization} on the other hand shows that ``solution is close to the original'', and thus that formulating the inverse problem using a variational formulation is reasonable.
\begin{figure}
    \centering
    \includegraphics[width=.75\textwidth]{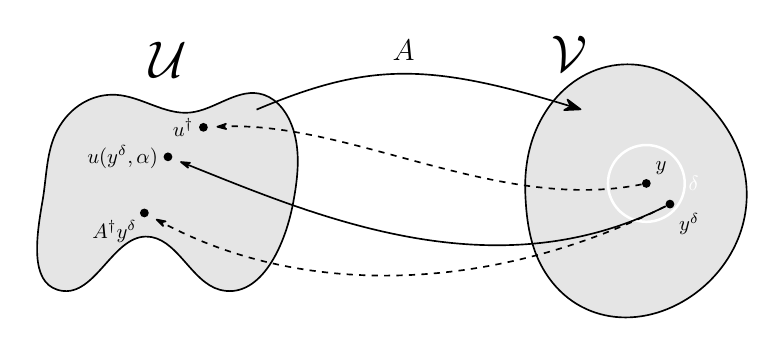}
    \caption{Regularization visualized.}
    \label{fig:reg}
\end{figure}

\subsection{Regularization as maximum a-posteriori estimation}\label{sec:bayesian}
Another way of tackling problems arising from ill-posedness is by adopting the Bayesian perspective for inversion \citep{kaipio2006statistical,stuart2010inverse}. The idea of statistical inversion methods is to rephrase the inverse problem as a question of statistical inference. We consider the problem
$$
Y=A U+N,
$$
where the measurement $Y$, unknown $U$, and noise $N$ are now modeled as random variables. This approach allows us to model the noise through its statistical properties. Through this lens, we can encode a-priori knowledge of the unknown in the form of a \textit{prior probability distribution}, assigning higher probability to those values of $u$ we expect to see. The desired solution then is no longer a vector; instead, it is a so-called \textit{posterior distribution}, which is the conditional probability distribution of $u$ given measurement $y$. This distribution can then be used to obtain estimates that are most likely in some sense. The main distinction from the functional analytic approach is that in the former case, we attain only a single solution; in the Bayesian viewpoint, a distribution of solutions is attained, from which one can either sample, consider point estimates like means or modes, as well as achieve uncertainty quantification by observing the variance of the posterior reconstructions.

\noindent Recall the \textbf{Bayes theorem}, which provides us with a way to statistically invert: for ${y} \in \mathbb{R}^m, {u} \in \mathbb{R}^n$
$$
p({u} \mid {y})=\frac{p({y} \mid {u}) p({u})}{p({y})}.
$$
The likelihood $p({y} \mid {u})$ is determined by the forward model and the statistics of the measurement error, and $p({u})$ encodes our prior knowledge on ${u}$. In practice, $p({y})$ is usually a normalizing factor which may be ignored. The \textit{maximum a-posteriori} (MAP) estimate is the maximizer ${u}^*$ of the posterior distribution: 
$$
p\left({u}^* \mid {y}\right)=\max _{{u}} p({u} \mid {y})=\max _{{u}}\{p({y} \mid {u}) p({u})\}
$$parameter and noise level
This interpretation provides a useful connection between the variational formulation and noise statistics. To be precise, the likelihood can be interpreted as a fidelity term, and the prior as regularization. See e.g. \citep{pereyra2019revisiting,tanunsupervised}. For Gaussian noise, for example, the regularization parameter is given by the variance of the noise model as shown below. 
\begin{example}[MAP for Gaussian noise]
 Consider an inverse problem $y = Au+n$, where $n$ is modeled as a zero-mean Gaussian $n \sim \mathcal{N}(0, \sigma^2I_m)$. Suppose we have a prior distribution on the unknowns $p(u)=e^{-\mathcal{R}(u)}$. The posterior distribution for $u$ knowing $y$ is given by Bayes:
$$
p(u \mid y)=\frac{p(y\mid u) p(u)}{p(y)} \propto e^{-\frac{1}{2 \sigma^2} \left\|A u-y\right\|^2 - \mathcal{R}(u)}
$$
In this case, the MAP image reconstruction is the one that maximizes this probability, or equivalently solves the variational regularization minimization problem
$$
\min _u\left\{\mathcal{R}(u)+\frac{1}{2 \sigma^2} \left\|y-Au\right\|^2\right\} .   
$$

\end{example}

%% file: chapter2.tex
\chapter{Variational Models and PDEs for Inverse Imaging}\label{cap2}
In this chapter we will explore the relationship between variational models and the theory of Partial Differential Equations (PDEs). This connection will equip us with a powerful analytical and computational tools for analyzing and solving inverse problems.
We begin by investigating the impact of various regularization choices, with a particular emphasis on Total Variation (TV) regularization due to its efficacy in preserving edges - a property fundamental in image processing tasks. Because regularizers like TV often result in non-smooth optimization problems, we will then introduce key concepts from classical convex analysis.  This provides the essential numerical tools for minimizing such variational energies. These tools are fundamental to many established methods and have led to a diverse ``zoo'' of regularizers designed for various imaging tasks. Finally, we will finish the section by highlighting the inherent limitations of purely model-driven approaches, thereby motivating the exploration of data-driven and hybrid techniques in subsequent sections.

\noindent In the introduction, we primarily considered finite-dimensional Euclidean spaces. However, as we shift our attention to modeling inverse problems in imaging, it becomes necessary to return to the continuous setting. In this framework, we typically consider the unknown image $u$ as a function in $L^2(\Omega)$, where $\Omega$ is an open and bounded domain with a Lipschitz boundary (often a rectangle in $\mathbb{R}^2$). The transformation $A$, which maps the true image to the observed data, is taken to be a bounded linear operator mapping from $L^2(\Omega)$ into itself.

\section{Variational Models}
\begin{figure}[ht!]
    \newlength{\commonDHeight}
    \setlength{\commonDHeight}{4.4cm} %
    \newlength{\commonDWidth}
    \setlength{\commonDWidth}{4.4cm} %
    \centering
    \begin{tabular}{c c c}
    \textbf{Gaussian} & \textbf{Poisson} & \textbf{Impulse}\\ %
    $\mathcal{D}(Au,f)=\|Au-f\|_2^2$ & $\mathcal{D}(Au,f)=\int Au-f\log(Au)~dx$ & $\mathcal{D}(Au,f)=\|Au-f\|_1$ \vspace{0.2cm}\\ %
    
    \includegraphics[height=\commonDHeight, width=\commonDWidth]{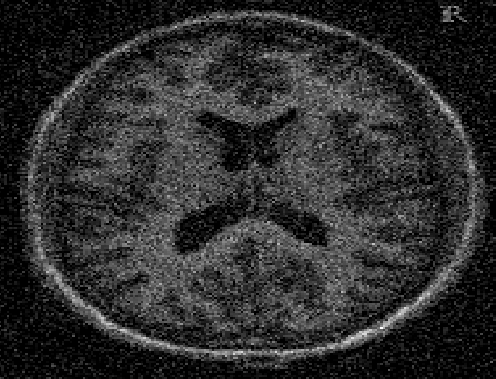} & %
    \includegraphics[height=\commonDHeight, width=\commonDWidth]{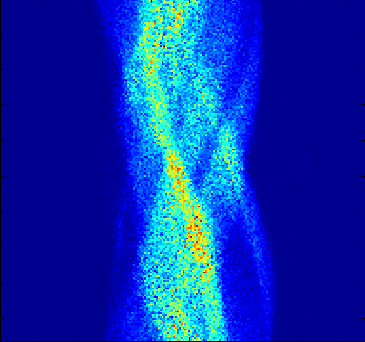} & %
    \includegraphics[height=\commonDHeight, width=\commonDWidth]{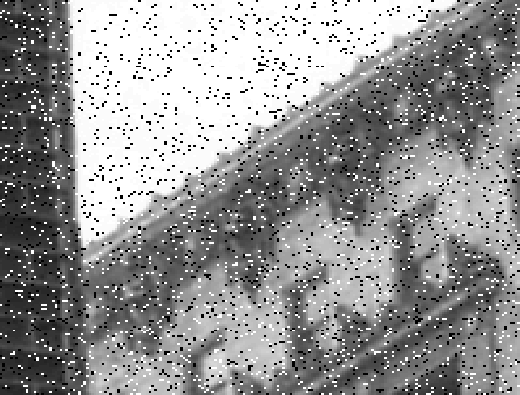} \\ %
    \includegraphics[height=\commonDHeight, width=\commonDWidth]{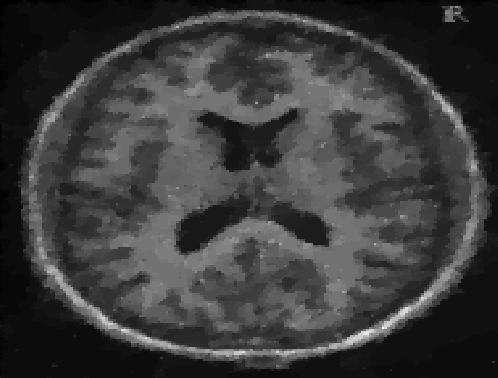} & %
    \includegraphics[height=\commonDHeight, width=\commonDWidth]{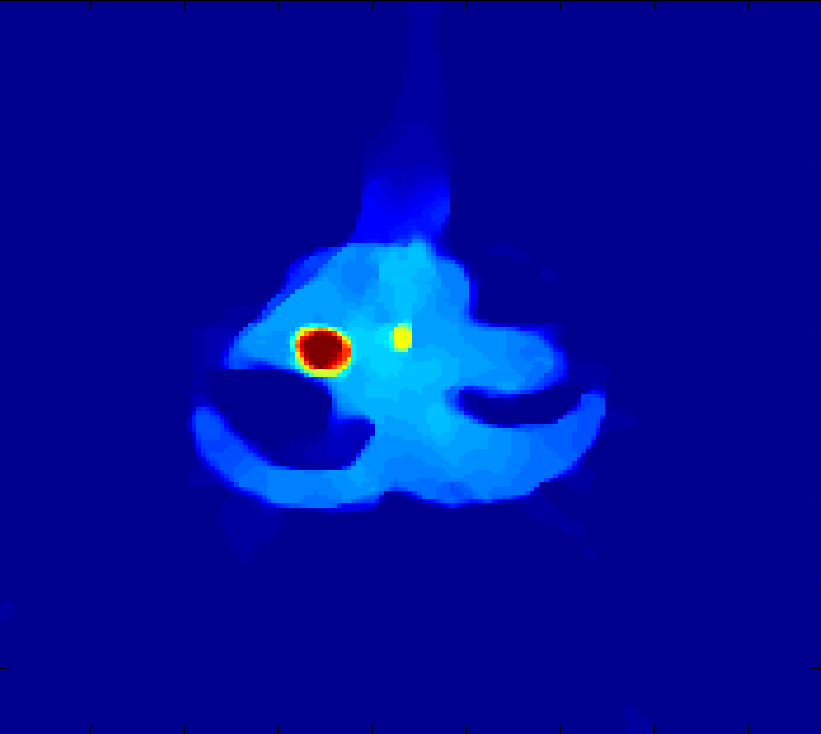} & %
    \includegraphics[height=\commonDHeight, width=\commonDWidth]{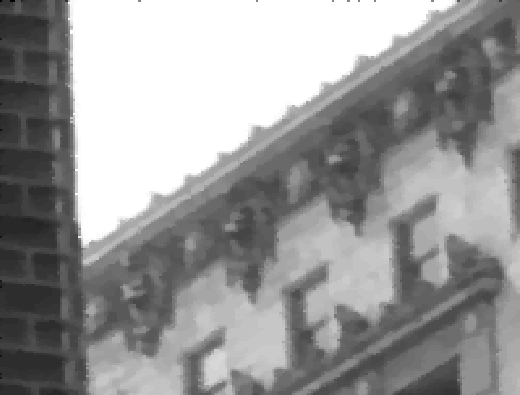}\\ %
    MRI & PET\footnote{Data courtesy of EIMI, M\"unster.} & Sparse noise
    \end{tabular}
    \caption{Examples of different noise models and corresponding data fidelity terms, with example images. See works by \citet{werner2012convergence,hohage2014convergence}.}
    \label{fig:noise_models_examples} %
\end{figure}

The reconstruction process of recovering the unknown image $u$ from the observed data $y$ in many inverse imaging problems is formulated as a minimization problem. This is typically  formulated as (see e.g. \citet{engl1996regularization,benning2018modern}):
\begin{equation}\label{eq:ip_main}
    \min_u \left\{ \alpha \mathcal{R}(u) + \mathcal{D}(Au, y) \right\}.
\end{equation}
Here, the data fidelity $\mathcal{D}(Au, y)$ enforces alignment between the forward model applied to $u$ and the observed data $y$. A common choice for $\mathcal{D}$ is a least-squares distance measure, such as $\frac{1}{2}\|A u-y\|^2$, however generally choice of $\mathcal{D}$ depends on data statistics (see \Cref{sec:bayesian}), and some examples are shown in \Cref{fig:noise_models_examples}. 
The second term $\mathcal{R}(u)$ is a functional that incorporates \textbf{a-priori information} about the image, acting as a regularizer, and $\alpha > 0$ is a weighting parameter that balances the influence of this prior information against fidelity to the data. Some basic examples of forward operators include:
\begin{itemize}[nosep]
    \item $A=\operatorname{Id}$ for image denoising,
    \item $A=1_{\Omega \backslash D}$ for image inpainting,
    \item $A=* k$ for image deconvolution,
    \item $A$ is (a possibly undersampled) Radon/Fourier transform for tomography.
\end{itemize}

\noindent A critical question then arises: \textbf{how do we choose an appropriate regularizer $\mathcal{R}(u)$ for images?} From a Bayesian perspective, the regularizer should be the data prior, penalizing characteristics likely introduced by noise. Beyond this, the choice of $\mathcal{R}(u)$ is guided by the specific properties we wish to enforce on the reconstructed image, such as smoothness or preservation of edges. The following examples illustrate this fundamental concept.

\begin{example}[1D Tikhonov]
Classical Tikhonov regularization often employs simple quadratic regularizers like $\mathcal{R}(u) = \frac{1}{2} \int_{\Omega} u^2 dx$ or, more commonly for images, $\mathcal{R}(u) = \frac{1}{2} \int_{\Omega} |\nabla u|^2 dx$. The latter term penalizes large gradients, encouraging smoothness in the solution. However, this choice implies that the reconstructed image $u$ possesses a certain degree of regularity. Most crucially for imaging, \textbf{the reconstruction cannot exhibit sharp discontinuities} like object boundaries or fine edges within an image. 

\noindent To see this, consider a one-dimensional scenario where $u:[0,1] \rightarrow \mathbb{R}$ and $u \in H^1(0,1)$, i.e. is $L^2$ with $L^2$ derivative. For any $0 < s < t < 1$, we have:
$$
u(t)-u(s)=\int_s^t u^{\prime}(r) dr \leq \sqrt{t-s} \sqrt{\int_s^t\left|u^{\prime}(r)\right|^2 dr} \leq \sqrt{t-s}\|u\|_{H^1(0,1)}
$$
This inequality shows that $u$ must be Hölder continuous with exponent $1/2$ (i.e., $u \in C^{1/2}(0,1)$), precluding jump discontinuities. 
\end{example}

\begin{example}[2D Tikhonov]
Extending this to a two-dimensional image $u \in H^1((0,1)^2)$, one can show that for almost every $y \in (0,1)$, the function $x \mapsto u(x,y)$ (a horizontal slice of the image) belongs to $H^1(0,1)$. This is because:
$$
\int_0^1\left(\int_0^1\left|\frac{\partial u(x, y)}{\partial x}\right|^2 dx\right) d y \leq \|u\|_{H^1}^2 < \infty
$$
This implies that $u$ cannot have jumps across vertical lines in the image (and similarly for horizontal lines). 
\end{example}

\section{Total Variation (TV) regularization}\label{sec:tvreg}

\begin{figure}
    \centering
    \begin{subfigure}{0.48\textwidth}
        \centering
        \includegraphics[width=\linewidth]{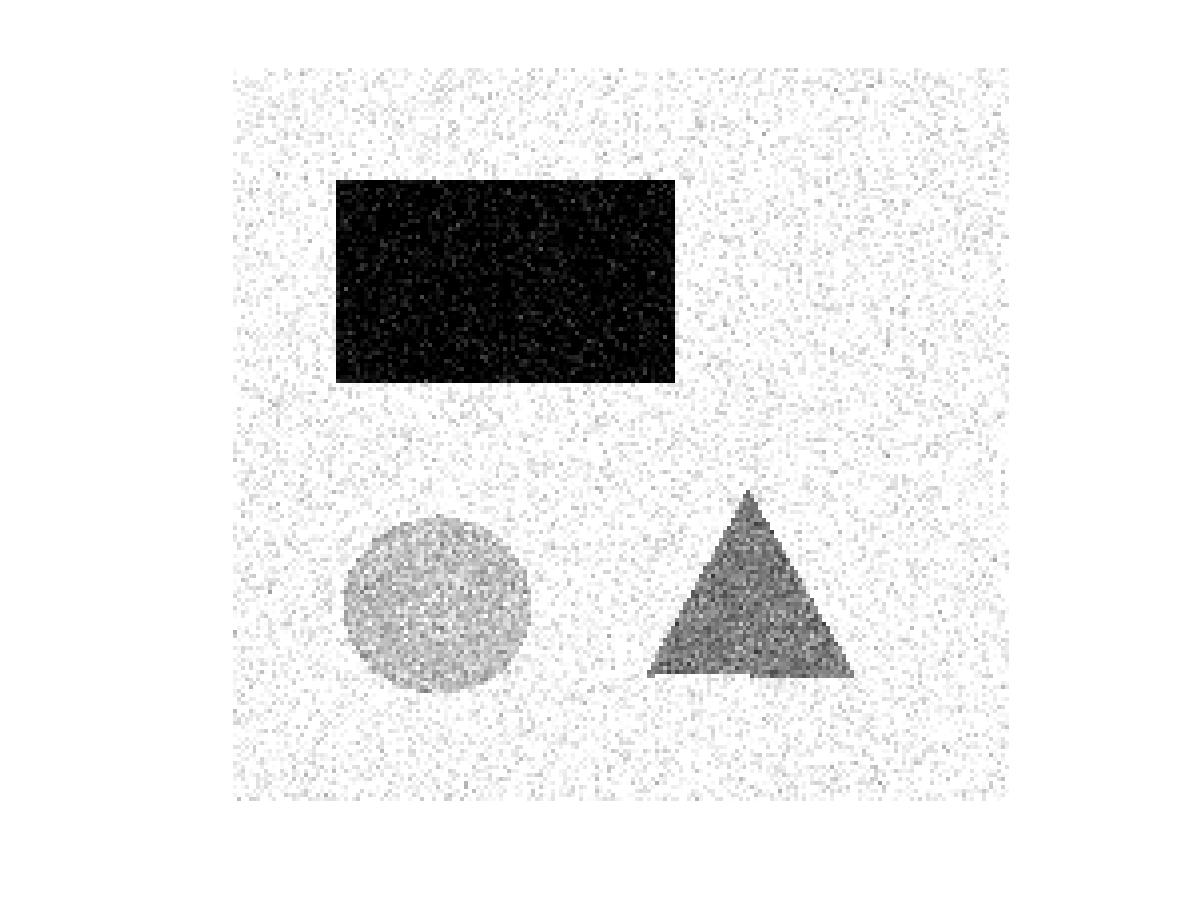}
    \end{subfigure}%
    \hfill
    \begin{subfigure}{0.48\textwidth}
        \centering
        \includegraphics[width=\linewidth]{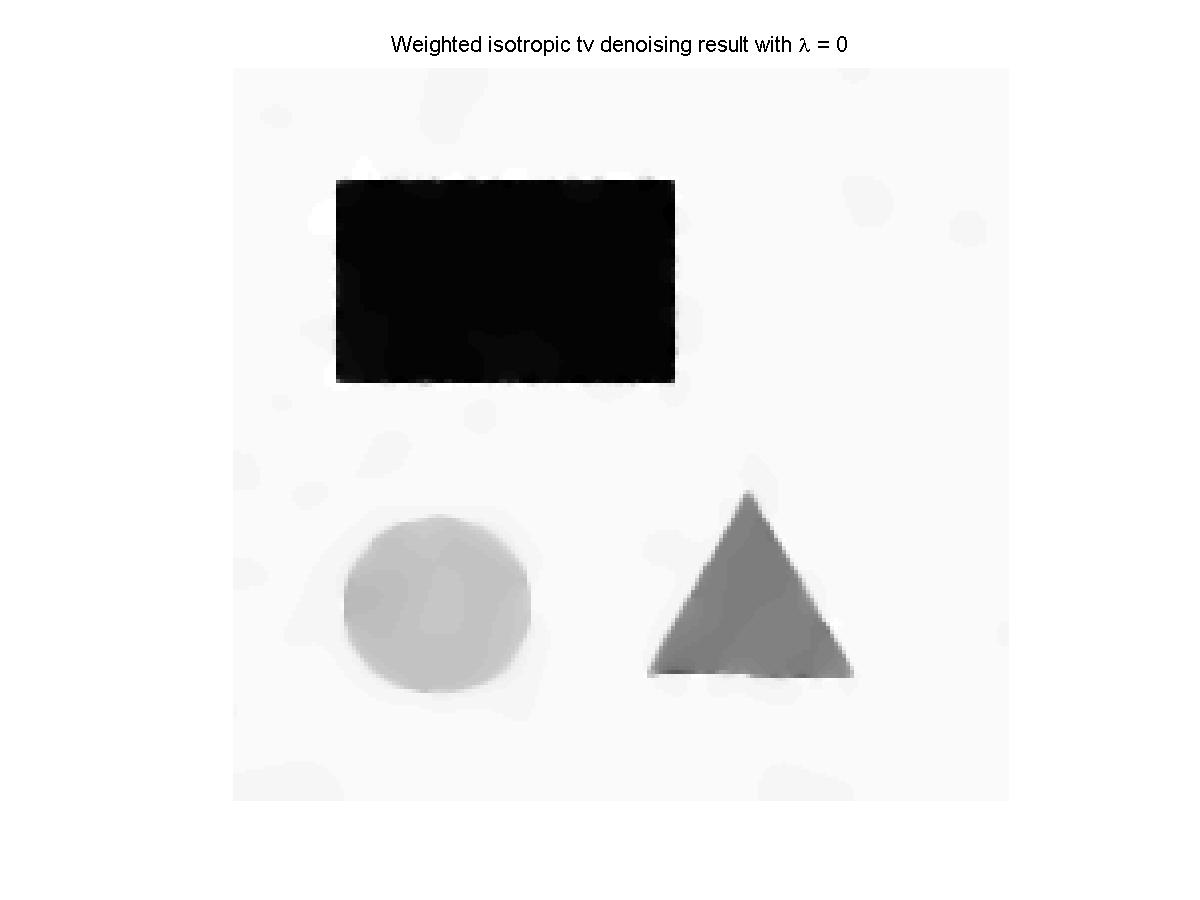}
    \end{subfigure}
    \caption{Example of TV denoised image of rectangles. 
 The total variation penalizes small irregularities/oscillations while respecting intrinsic image features such as edges.}
    \label{fig:combined}
\end{figure}

The examples above demonstrate that traditional Sobolev spaces like $H^1$ are often too restrictive for image processing because they penalize the very edges that define objects. To overcome this, we seek to relax the smoothness constraint by employing a weaker notion of the derivative, leading us to \textbf{Total Variation (TV) regularization}.

\begin{example}[Total Variation regularization]

In TV regularization, the functional used is $\mathcal{R}(u) = |Du|(\Omega)$, which represents the total variation of $u$ over the domain $\Omega$. For a locally integrable function $u \in L_{loc}^1(\Omega)$, its variation is defined as:
$$V(u, \Omega) := \sup \left\{ \int_{\Omega} u \nabla \cdot \varphi dx : \varphi \in C_c^1(\Omega; \mathbb{R}^n), \|\varphi\|_{\infty} \leq 1 \right\}.$$
A function $u$ belongs to the space of \textbf{functions of Bounded Variation}, denoted $BV(\Omega)$, if and only if its variation $V(u, \Omega)$ is finite. For such functions, the total variation is given by:
$$|Du|(\Omega) = V(u, \Omega)$$
where $|Du|(\Omega)$ is the total mass of the Radon measure $Du$, which is the derivative of $u$ in the sense of distributions.

\end{example}

\begin{figure}[ht] %
    \centering
    \newlength{\commonTVHeight}
    \setlength{\commonTVHeight}{3.4cm} %

    \begin{subfigure}[b]{0.24\textwidth} %
        \centering
        \includegraphics[height=\commonTVHeight]{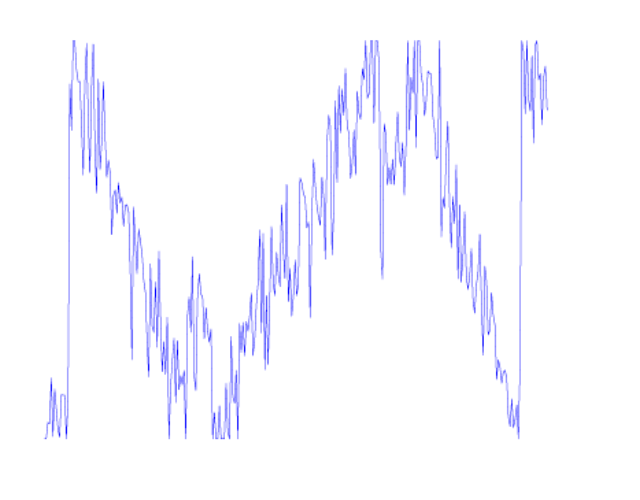} %
        \caption{Large TV} %
        \label{sfig:tv_expl1}
    \end{subfigure}
    \hfill %
    \begin{subfigure}[b]{0.24\textwidth}
        \centering
        \includegraphics[height=\commonTVHeight]{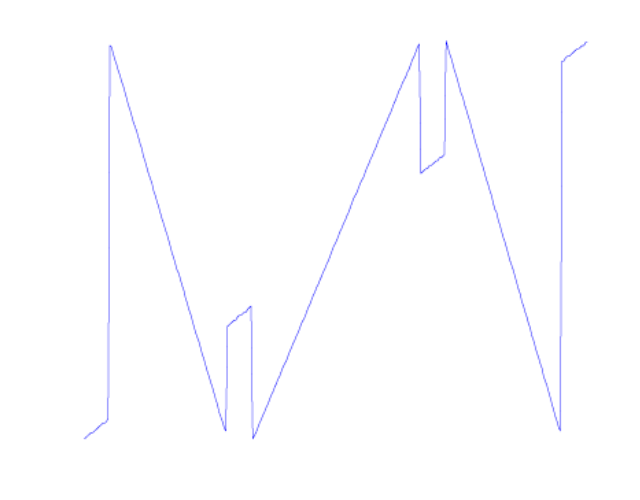} %
        \caption{Small TV}
        \label{sfig:tv_expl2}
    \end{subfigure}
    \hfill
    \begin{subfigure}[b]{0.48\textwidth}
        \centering
        \includegraphics[height=\commonTVHeight]{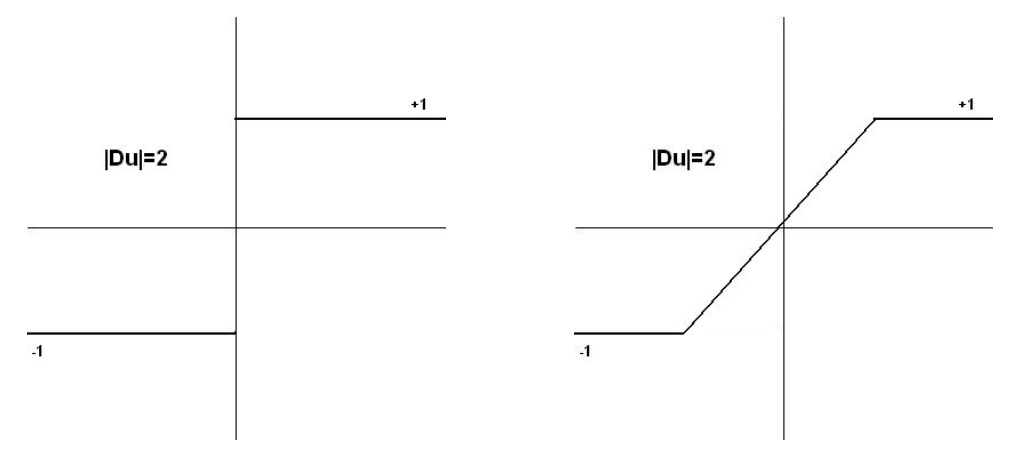} %
        \caption{Examples of functions with equal TV}
        \label{sfig:tv_expl_general}
    \end{subfigure}

    \caption[Properties of Total Variation Smoothing]{Properties of Total Variation (TV) smoothing.
    (\subref{sfig:tv_expl1}-\subref{sfig:tv_expl2}) 
    TV penalizes small irregularities and oscillations, and tends to preserve edges.
    (\subref{sfig:tv_expl_general}) The total variation measures the size of the jump discontinuity.
    Overall, the total variation penalizes small irregularities/oscillations while respecting intrinsic image features such as edges \citep{rudin1992nonlinear}.
    }
    \label{fig:tv_properties_combined}
\end{figure}

\begin{figure}[ht] %
    \centering
    \begin{subfigure}{0.25\textwidth}
        \includegraphics[width=\linewidth]{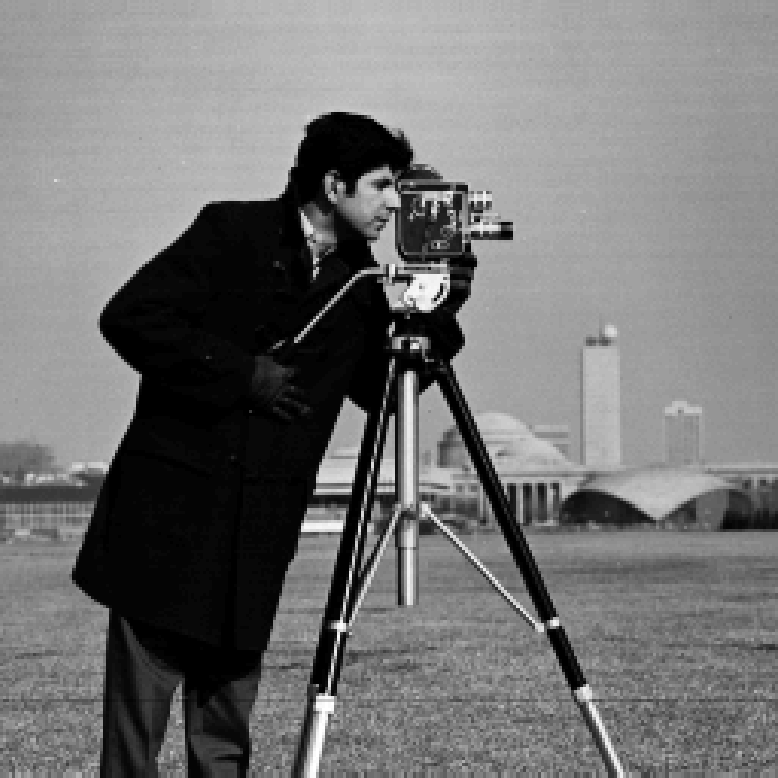}
        \caption{original}
    \end{subfigure}
    \hspace{1cm}
    \begin{subfigure}{0.25\textwidth}
        \includegraphics[width=\linewidth]{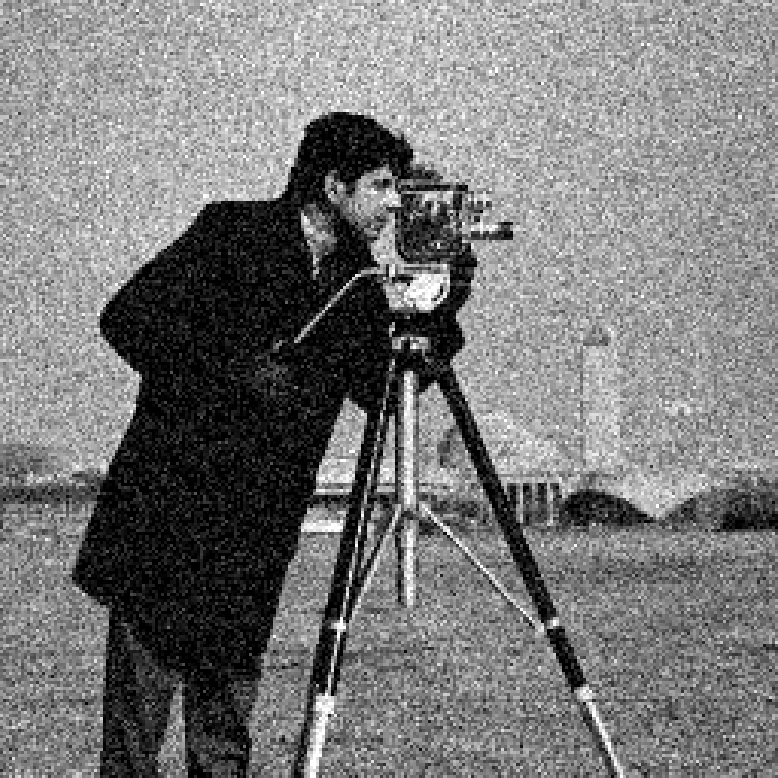}
        \caption{noisy}
    \end{subfigure}
    \\ %
    \vspace{0.5cm} %
    \begin{subfigure}{0.25\textwidth}
        \includegraphics[width=\linewidth]{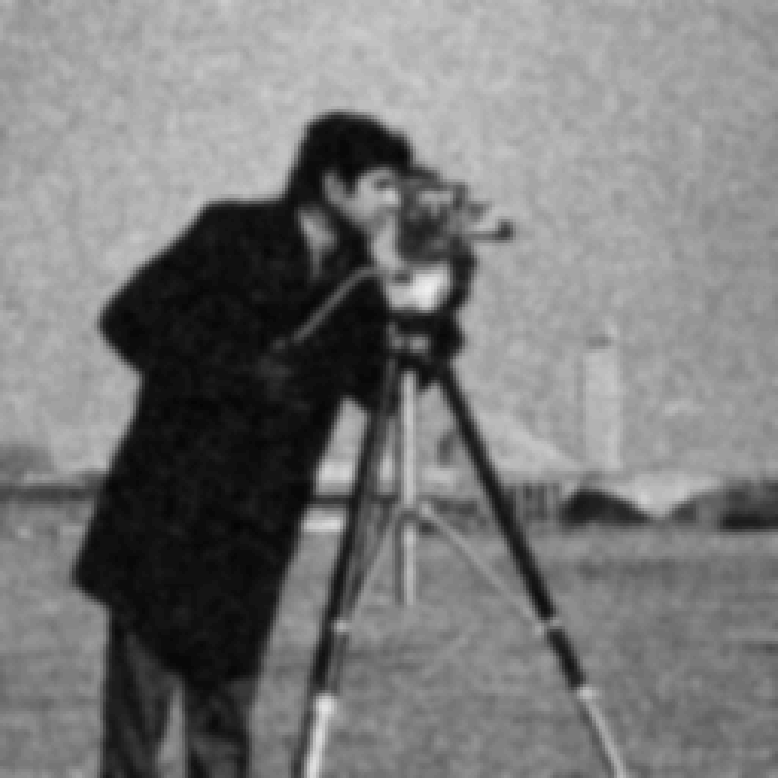}
        \caption{$\|\nabla u\|_2^2$}
    \end{subfigure}
    \hspace{1cm}
    \begin{subfigure}{0.25\textwidth}
        \includegraphics[width=\linewidth]{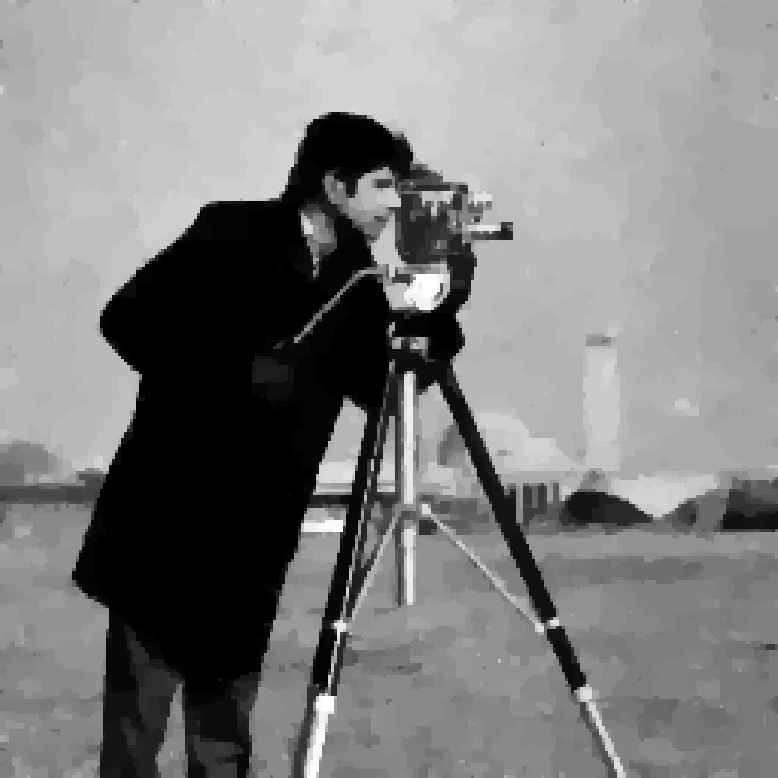}
        \caption{$\|\nabla u\|_1$}\label{sfig:tvl1}
    \end{subfigure}
    \caption{Comparison of regularization methods. This showcase that convex relaxation of $l_0$ of TV with $l_1$ successfully achieves sparsity and is a more natural prior for denoising.\label{fig:cameraman}
    } %
\end{figure}

The space $BV(\Omega)$ is particularly well-suited for images because, unlike $H^1(\Omega)$, $BV$ functions \textbf{can have jump discontinuities (edges)}. Minimizing the total variation penalizes small irregularities and oscillations while respecting instrinsic image features such as edges. See \Cref{fig:tv_properties_combined} for a visualisation of properties of TV. Heuristically, the total variation of a function quantifies the ``amount of jumps'' or oscillations it contains; thus, noisy images, which typically have many rapid oscillations, have a large TV value. 
\noindent Owing to these desirable properties, TV regularization has become a widely used technique in image processing and inverse problems. It promotes solutions that are \textbf{piecewise smooth} yet retain sharp edges, a property that is crucial in image processing, see e.g. \Cref{fig:tv_properties_combined,sfig:tvl1}. The non-differentiability of the TV term, however, necessitates specialized optimization algorithms, such as primal-dual methods or general splittings \citep{lions1979splitting,combettes2005signal,hintermuller2003semi}.

\begin{example}[Compressed sensing]
In compressed sensing \citep{candes2006stable,poon2015role}, TV regularization plays a vital role. Images often exhibit sparse gradients (large areas of constant intensity), a key assumption in compressed sensing. For $u \in W^{1,1}(\Omega)$, the total variation coincides with the $L^1$ norm of its gradient:
$$
|Du|(\Omega) = \|\nabla u\|_{L^1(\Omega)}
$$
The $L^1$ norm is well-known for promoting \textbf{sparsity}. While the $L^0$ norm (counting non-zero gradient values) would be ideal for enforcing sparse gradients, it leads to computationally intractable (NP-hard) problems. The $L^1$ norm serves as a convex relaxation, making optimization feasible while still encouraging solutions with few non-zero gradient values, characteristic of piecewise constant regions. Remarkably, if the underlying data is indeed sparse, TV regularization enables near-perfect reconstruction from significantly undersampled data, for example \Cref{fig:mri_reconstruction_example,fig:cameraman}.

\end{example}
\begin{example}[MRI]
\textit{Magnetic Resonance Imaging} ( \citep{lustig2008sparse,fessler2008image}) is a medical imaging technique that measures the response of atomic nuclei in a strong magnetic field. The measured data in MRI is essentially a sampled Fourier transform of the object being imaged. In many MRI applications, acquiring a full set of Fourier measurements is time-consuming and can be uncomfortable for the patient.  Compressed sensing offers a way to speed up the process by acquiring only a subset of the Fourier data. This is known as undersampled Fourier acquisition, see \Cref{sfig:mri_undersampled}.
$$
y=(\mathcal{F} u)_{\mid \Lambda}+n,
$$
where $\mathcal{F}$ denotes the Fourier transform operator,  $\Lambda$ is the set of (undersampled) Fourier coefficients and $n$ represents noise in the measurement process.
We often seek a piecewise constant image, which means the image has distinct regions with constant intensities (like different tissues in the body). The goal, then, is to identify a piecewise constant function $u$ consistent to the datum $g$. As before, minimizing $\|\nabla u\|_0$ under data consistency is NP-hard, and $\ell_1$ can be used as convex relaxation, leading to total variation minimization:
$$
 \min _u\quad  \alpha\|\nabla u\|_1+\frac{1}{2}\left\|(\mathcal{F} u)_{\mid \Lambda}-y\right\|^2. 
$$
\begin{figure}[ht]
    \centering
    \begin{subfigure}[b]{0.24\textwidth}
        \centering
        \includegraphics[width=.96\linewidth]{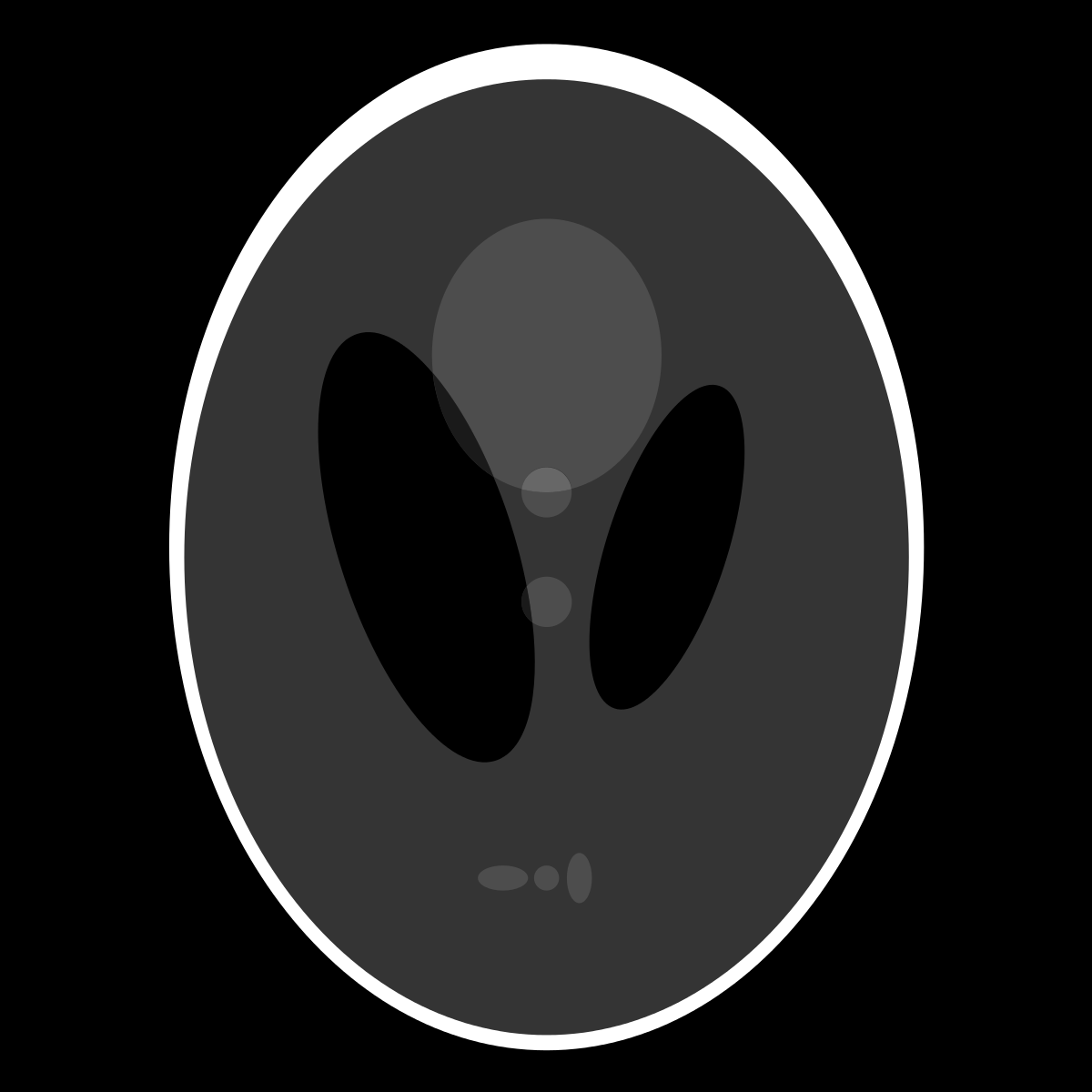} %
        \caption{Ground truth}
        \label{sfig:shepplogan_gt}
    \end{subfigure}
    \hfill
    \begin{subfigure}[b]{0.24\textwidth}
        \centering
        \includegraphics[width=\linewidth]{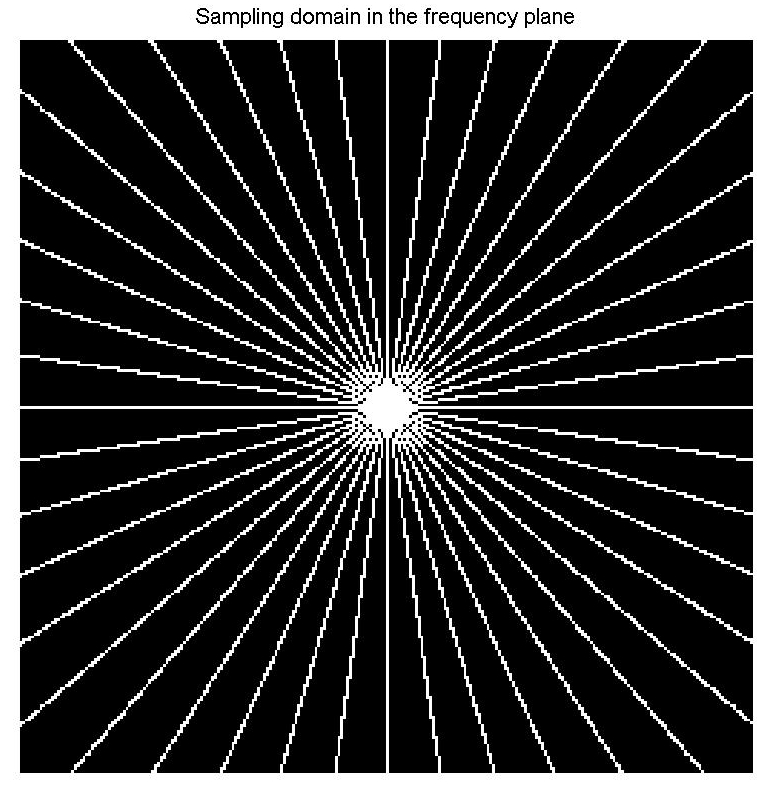} %
        \caption{Undersampled Fourier}
        \label{sfig:mri_undersampled}
    \end{subfigure}
    \hfill
    \begin{subfigure}[b]{0.24\textwidth}
        \centering
        \includegraphics[width=\linewidth]{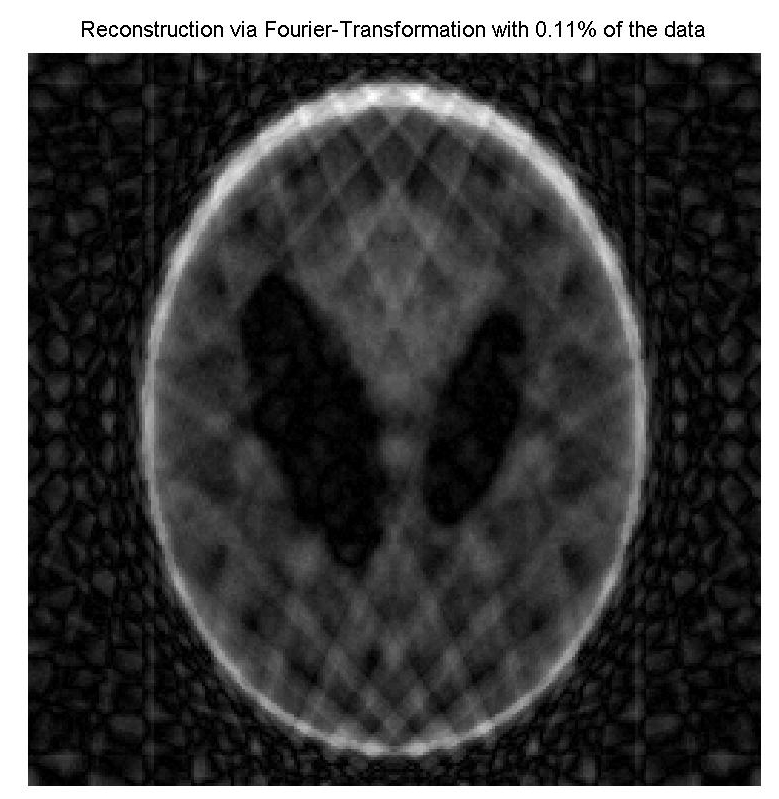} %
        \caption{Zero-filling}
        \label{sfig:mri_zerofill}
    \end{subfigure}
    \hfill
    \begin{subfigure}[b]{0.24\textwidth}
        \centering
        \includegraphics[width=\linewidth]{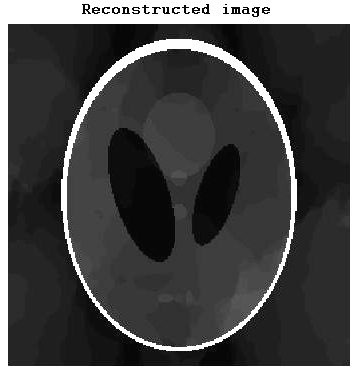} %
        \caption{TV solution}
        \label{sfig:mri_tv_solution}
    \end{subfigure}

    \caption{MRI reconstruction example:
    (\subref{sfig:shepplogan_gt}) Ground truth Shepp-Logan phantom.
    (\subref{sfig:mri_undersampled}) Undersampled k-space (Fourier) data.
    (\subref{sfig:mri_zerofill}) Reconstruction via zero-filling the undersampled k-space and inverse Fourier transform.
    (\subref{sfig:mri_tv_solution}) Reconstruction using a Total Variation (TV) regularized approach.
    }
    \label{fig:mri_reconstruction_example}
\end{figure}
\end{example}
Another key insight is that the TV measure can be interpreted as an accumulation of the perimeters of all level sets. Thus, penalizing the TV encourages the stretching or smoothing of these level sets, leading to a reduction in their overall length. This property makes TV regularization particularly well-suited for segmentation problems, where the goal is to partition an image into distinct regions with well-defined boundaries.

\begin{example}[Sets of Finite Perimeter and the Co-area Formula]

For instance, if $\Omega \subset \mathbb{R}^2$ is an open set and $D$ is a subset with a $C^{1,1}$ boundary, the total variation of its characteristic function $u = \chi_D$ (1 inside $D$, 0 outside) is simply the perimeter of $D$ within $\Omega$: $|Du|(\Omega) = \mathcal{H}^1(\partial D \cap \Omega)$. See \citet{ambrosio2000functions}. More generally, the \textbf{co-area formula} states that for any $u \in BV(\Omega)$:
$$
|Du|(\Omega) = \int_{-\infty}^{+\infty} \operatorname{Per}(\{u>s\} ; \Omega) ds
$$
where $\operatorname{Per}(\{u > s\}; \Omega) = \|D\chi_{\{u > s\}}\|(\Omega)$ is the perimeter of the superlevel set of $u$ at level $s$. This formula reveals that the total variation of $u$ is the integral of the perimeters of all its level sets.

\end{example}
In the context of image segmentation, minimizing the TV of an image encourages the boundaries between different regions to be smooth and well-defined. This is because minimizing the perimeters of the level sets leads to a reduction in the overall length and complexity of the boundaries. This has made TV a fundamental tool in numerous segmentation algorithms, for example through the Chan--Vese model.
\begin{example}[Chan--Vese Segmentation]
The Chan--Vese model \citep{chan2001active} is a popular variational approach for image segmentation that leverages the TV regularization. It stems from the Mumford--Shah functional \citep{mumford1989optimal}, which aims to find an optimal piecewise smooth approximation of a given image. The Chan--Vese model simplifies this by assuming that the image can be segmented into regions with constant intensities.

\noindent Let $\Omega\subset\mathbb{R}^2$ represent the image domain, and let $y:\Omega\to\mathbb{R}$ denote the given image. The Chan--Vese model seeks to partition $\Omega$ into two regions, represented by a binary function $\chi:\Omega\to\{0,1\}$. The objective functional to be minimized is:
$$
\min _{\chi, c_1, c_2} \quad\alpha|D \chi|(\Omega)+\int_{\Omega}\left(y-c_1\right)^2 \chi+\int_{\Omega}\left(y-c_2\right)^2(1-\chi),
$$
where $\chi$ is the binary segmentation function, $c_1$ and $c_2$ are the average intensities within the regions where $\chi=1$ and $\chi=0$, respectively.
\begin{figure}[hbt!] %
    \centering
    \includegraphics[width=0.8\textwidth]{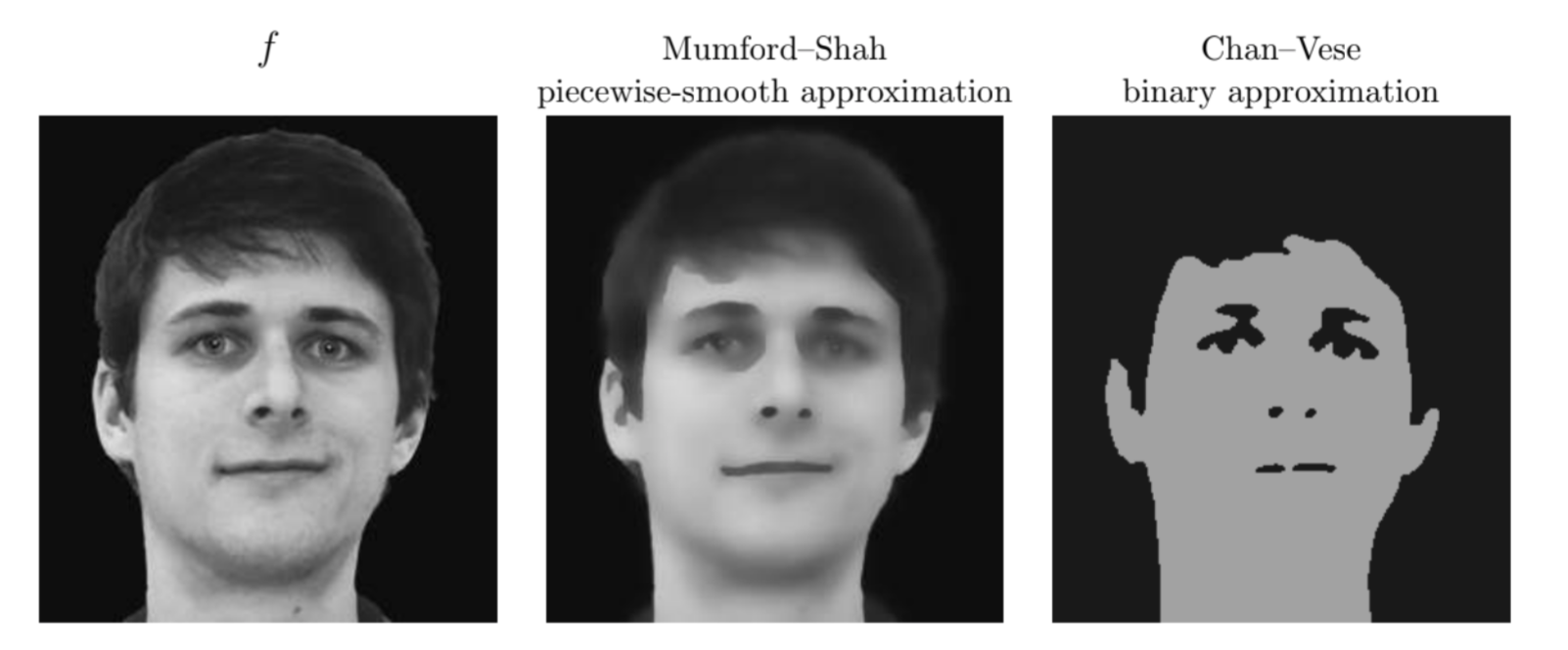} %
    \caption{Example of binary Chan--Vese segmentation compared to Mumford-Shah segmentation. \citep{mumford1989optimal,pock2009algorithm,getreuer2012chan}.
    }
    \label{fig:chan_vese_example}
\end{figure}
\noindent Solving the original Chan--Vese formulation with the binary constraint is computationally challenging.  A common approach \citep{cai2013two,cai2019linkage} to address this is to relax the binary constraint. This leads to the following convex optimization problem (with given $c_1$ and $c_2$):
$$
\min _v \quad\alpha|D v|(\Omega)+\int_{\Omega}\left(y-c_1\right)^2 v+\int_{\Omega}\left(y-c_2\right)^2(1-v),
$$
with the relaxed segmentation function $v \in[0,1]$. The final binary segmentation is then typically obtained by thresholding the resulting $v$. See \Cref{fig:chan_vese_example} for a visual example.
\end{example}

\section{From Total Variation Regularization to Nonlinear PDEs}

The study of variational methods such as TV regularization is deeply enriched by its connection to the theory of Partial Differential Equations (PDEs).
PDEs frequently arise as mathematical descriptions of the gradient flow minimizing some energy functional. This dynamic perspective is instrumental not only for characterizing important properties of solutions, such as scale-space \citep{perona1990scale,burger2006nonlinear,florack2000topological} and image decomposition \citep{rudin1992nonlinear,ambrosio2001connected,caselles1993geometric,aujol2005image,alvarez1992image,caselles1998introduction,chambolle2010introduction}, but also for understanding the continuous analog of iterative optimization algorithms, thereby informing the analysis of their discrete counterparts. Moreover, the PDE framework itself can motivate the development of novel regularizers or reconstruction paradigms, occasionally leading to methods like Cahn--Hilliard inpainting \citep{doi:10.1137/080728548} that may not possess an explicit variational formulation. The extensive analytical toolkit of PDE theory further allows for rigorous investigation into the qualitative aspects of solutions, such as their regularity. This PDE-centric view enriches the understanding of TV-based methods and is pivotal in guiding the design of efficient numerical solvers, a concept we will illustrate with the Rudin--Osher--Fatemi (ROF) model \citep{rudin1992nonlinear}.

\begin{example}[Nonlinear Image Smoothing: The ROF Model]\label{eg:ROF}
The ROF model seeks to determine an image $u$ that remains close to a noisy observation $y$ while also possessing a minimal total variation. The corresponding optimization problem is:
$$
\min_u \left( \alpha |Du|(\Omega) + \frac{1}{2} \|u - y\|^2 \right).
$$
To understand the process by which $u$ evolves to minimize this energy functional, we can examine its (sub)gradient flow. This flow describes the path of steepest descent for the functional. It is given by the differential inclusion \citep{bellettini2002total}:
$$u_t = \alpha p + (u-y), \quad \text{where } p \in \partial(|Du|), \quad \text{in } \Omega.$$
In this expression, $u_t$ denotes the derivative of $u$ with respect to an artificial time variable $t$ (representing the evolution of the flow), and $p$ is an element from the subdifferential of the TV term.

\noindent In regions of the image where the gradient is component-wise non-zero, the subdifferential reduces to a singleton, and the gradient flow to be expressed more explicitly as the following nonlinear PDE:
$$u_t = \alpha \operatorname{div}\left(\frac{D u}{|D u|}\right) + (u-y), \quad \text{in } \Omega.$$
This equation is a nonlinear diffusion equation. Its key characteristic is that the effective diffusion coefficient is inversely proportional to the magnitude of the image gradient, $|D u|$. 
This property leads to a highly desirable selective smoothing behavior: in relatively flat regions of the image where $|D u|$ is small (often dominated by noise), the diffusion is strong, leading to significant smoothing. Conversely, near sharp edges where $|Du|$ is large, the diffusion is weak, which helps to preserve these important structural features of the image while reducing noise elsewhere.
\end{example}

\section{Numerical Aspects}
Having explored various regularization functionals and their impact on reconstructed image properties, we now turn to the crucial aspect of their practical implementation: numerical optimization. This section discusses established algorithms for finding minimizers of the energy functionals that arise in variational image reconstruction, largely drawing from the framework presented in \citet{chambolle2016introduction}. We typically consider optimization problems in a finite-dimensional setting, $u \in \mathbb{R}^n = X$. The general form of the minimization problem is:
$$
\min_{u\in X} \mathcal{J}(u)+\mathcal{H}(u),
$$
where $\mathcal{J}$ and $\mathcal{H}$ are proper and convex functions, and one or both may be Lipschitz differentiable. Note that for this section the notation has changed from the usual, as will become clear from \Cref{eg:ROF_opt} - both $\mathcal{J}$ and $\mathcal{H}$ can take on the rule of the data fidelity.

\begin{example}[ROF problem]\label{eg:ROF_opt}
For a given noisy image $y \in \mathbb{R}^n$, the ROF model from \Cref{eg:ROF} seeks to find an image $u$ by solving:
$$
\min _u\; \alpha\|\nabla u\|_{2,1}+\frac{1}{2}\|u-y\|_2^2,
$$
where the TV term is  $\|\nabla u\|_{2,1}=\sum _{ij}\left|(\nabla u) _{ij}\right|_2=\sum _{ij} \sqrt{\left(u_x\right)_{i j}^2+\left(u_y\right)_{ij}^2}$.

\noindent Several algorithms have been developed to compute minimizers of such functionals. One approach involves regularising the TV term to make it differentiable. For instance, one might consider instead solving the following regularized ROF problem:
$$
\min _u\left\{\alpha \sum \sqrt{u_x^2+u_y^2+\epsilon}+\frac{1}{2}\|u-g\|_2^2\right\}
$$
for a small $0<\epsilon \ll 1$. The regularized TV is differentiable in the classical sense, therefore we can apply classical numerical algorithms to compute a minimizer, e.g. gradient descent, conjugate gradient methods etc.   
\end{example}
However, to address the original, unregularized problem, we need to employ techniques from convex analysis, specifically the concept of subgradients. We will consider the following basic properties \citep{chambolle2016introduction}.
\begin{definition}[Lower Semi-Continuity (lsc.)] \label{def:lsc}
A function $f: X \to \mathbb{R} \cup \{+\infty\}$ is {lower semi-continuous} (lsc) at a point $x_0 \in X$ if for any sequence $\{x_k\}_{k=1}^\infty$ in $X$ such that $x_k \to x_0$, we have:
$$ \liminf_{k \to \infty} f(x_k) \ge f(x_0) $$
A function $f$ is said to be lsc on $X$ if it is lsc at every point in $X$.
\end{definition}
\begin{definition}[Proper] \label{def:proper_convex}
An extended real-valued convex function $f: X \to \mathbb{R} \cup \{+\infty\}$ is called a proper if $\text{dom}(f) = \{x \in X : f(x) < +\infty\}$, is non-empty, and $f(x) > -\infty \;,\forall x \in X$.
\end{definition}
\begin{definition}[Subdifferential]
For a convex function $J: X \to \mathbb{R}$, we define the subdifferential of $J$ at $x \in X$, as $\partial J(x)=\emptyset$ if $J(x)=\infty$, otherwise
$$
\partial J(x):=\left\{p \in X^{\prime}:\langle p, y-x\rangle+J(x) \leq J(y) \quad \forall y \in X\right\},
$$
where $X^{\prime}$ denotes the dual space of $X$. It is obvious from this definition that $0 \in \partial J(x)$ if and only if $x$ is a minimizer of $J$.
   
\end{definition}
\begin{example}[Subdifferential of the $\ell_1$ norm]
To illustrate the concept of the subdifferential for a common non-smooth function in imaging, consider the $\ell_1$-norm. Let $X=\ell_1(\Lambda)$ and $J(x):=\|x\|_1$, with $\Lambda=[1,\dots,n]$ or $\mathbb{N}$. The subdifferential is given by:
$$
\partial\|\cdot\|_1(x)=\left\{\xi \in \ell_{\infty}(\Lambda): \xi_\lambda \in \partial|\cdot|\left(x_\lambda\right), \lambda \in \Lambda\right\},
$$
where $\partial|\cdot|(z)=\{\operatorname{sign}(z)\}$ if $z \neq 0$, and $\partial|\cdot|(0)=[-1,1]$.
\end{example}
\paragraph{The Legendre--Fenchel transform}
In order to approach the problem in \Cref{eg:ROF_opt}, it will turn out to be useful to write down the dual formulation, for which we need to define the Legendre--Fenchel transform. 
\begin{definition}[Legendre-Fenchel Transform] \label{def:lf_transform}
The {Legendre-Fenchel transform} (also known as the convex conjugate) of a function $f: X \to \mathbb{R} \cup \{\pm\infty\}$ is a function  $f^*: X' \rightarrow \mathbb{R} \cup\{-\infty,+\infty\}$ defined as
$
f^*\left(y\right):=\sup_{x \in X} \left\{\left\langle y, x\right\rangle-f(x)\right\}
$
where $\langle x, y \rangle$ denotes the inner product of $x$ and $y$. In particular, $f^*$ is convex and lsc.
\end{definition}
\noindent What is more, when $f$ is a proper, convex, lsc function, applying the Legendre-Fenchel transform twice returns the original function:
$\forall x \in \mathbb{R}^n: \;  f^{**}(x) = f(x) $.

\begin{example}[One homogeneous functions]
For example, for a function $J$ that is one-homogeneous (i.e., $J(\lambda u)=\lambda J(u)$ for every $u$ and $\lambda>0$), its Legendre--Fenchel transform is the characteristic function of a closed convex set $K$:
$$
J^*(v)=\chi_K(v)= \begin{cases}0 & \text { if } v \in K, \\ +\infty & \text { otherwise. }\end{cases}
$$
Since $J^{**}=J$ (as $J$ is proper, convex, and lsc), we can recover $J(u)$ from its transform:
$$
J(u)=\sup _{v \in K}\langle u, v\rangle_X.
$$    
\end{example}

\paragraph{Proximal Map}
To effectively handle non-differentiable terms in optimization problems and to construct iterative algorithms, the concept of the proximal map is essential. For a convex, proper, and lower semi-continuous function $J: X \rightarrow \mathbb{R} \cup \{+\infty\}$, its proximal map at a point $y \in X$, with parameter $\tau > 0$, is defined as the unique minimizer of the following problem:
$$
\operatorname{prox}_{\tau J}(y) = \underset{u\in X}{\operatorname{argmin}}\left( J(u)+\frac{1}{2 \tau}\|u-y\|_2^2 \right).
$$
Let $u^* = \operatorname{prox}_{\tau J}(y)$. The optimality condition for this minimization is:
$$
0 \in \partial J\left(u^*\right)+\frac{u^*-y}{\tau},
$$
which can be rewritten as $u^*=(I+\tau \partial J)^{-1} y$.
Furthermore, Moreau's identity provides a relationship between the proximal map of $J$ and its convex conjugate $J^*$:
$$
y=\operatorname{prox}_{\tau J}(y)+\tau \operatorname{prox}_{\frac{1}{\tau} J^*}\left(\frac{y}{\tau}\right).
$$
This identity implies that if $\operatorname{prox}_{\tau J}$ is known, $\operatorname{prox}_{\frac{1}{\tau} J^*}$ can also be computed.

\subsection{Dual Problem}
Transforming the original (primal) optimization problem into its dual counterpart can often lead to more tractable problems or enable the development of efficient algorithms, especially when dealing with complex coupling terms. Consider the primal problem:
$$
\min _{u \in X} \mathcal{J}(A u)+\mathcal{H}(u),
$$
where $\mathcal{J}: Y \rightarrow(-\infty,+\infty]$ and $\mathcal{H}: X \rightarrow(-\infty,+\infty]$ are convex, lower semi-continuous (l.s.c.) functions, and $A: X \rightarrow Y$ is a bounded linear operator. Using the definition of the convex conjugate, we have $\mathcal{J}(Au) = \sup_{p \in Y} (\langle p, Au \rangle - \mathcal{J}^*(p))$. Substituting this into the primal problem leads to:
$$
\adjustlimits \min _{u \in X}\sup _{p \in Y} \left( \langle p, A u\rangle-\mathcal{J}^*(p)+\mathcal{H}(u) \right).
$$
Then, under certain mild assumptions on $\mathcal{J},\, \mathcal{H}$, e.g. in finite dimensions it is sufficient to have a point with $A x$ in the relative interior of $\operatorname{dom} \mathcal{J}$ and $x$ in the relative interior of $\operatorname{dom} \mathcal{H}$ see e.g. \citep[Section 3.1]{chambolle2016introduction} or \citep[Section 5.1]{sherry2025part}, one can swap the min and sup to have:
$$
\begin{aligned}
& \min _{u \in X} \mathcal{J}(A u)+\mathcal{H}(u) \\
& \underbrace{=}_{\mathcal{J}^{* *}=\mathcal{J}} \adjustlimits\min _{u \in X} \sup _{p \in Y}\langle p, A u\rangle-\mathcal{J}^*(p)+\mathcal{H}(u) \\
& =\max _p \inf _u\langle p, A u\rangle-\mathcal{J}^*(p)+\mathcal{H}(u) \\
& =\max _p-\mathcal{J}^*(p)-\mathcal{H}^*\left(-A^* p\right), 
\end{aligned}
$$
where the latter is the dual problem, $\mathcal{H}^*$ is the convex conjugate of $\mathcal{H}$, and $A^*$ is the adjoint of $A$. Under the above assumptions, there exists at least one solution $p^*$ (see, e.g., \citet{ekeland1999convex,borwein2015duality}).
If $u^*$ solves the primal problem and $p^*$ solves the dual problem, then $(u^*, p^*)$ is a saddle-point of the Lagrangian $\mathcal{L}(u,p)$ defined as follows, which provides a link between primal and dual solutions:
$$
\mathcal{L}(u, p):=\langle p, A u\rangle-\mathcal{J}^*(p)+\mathcal{H}(u),
$$
such that for all $(u, p) \in X \times Y$, we have $\mathcal{L}(u^*, p) \leq \mathcal{L}(u^*, p^*) \leq \mathcal{L}(u, p^*)$. Moreover, we can define the primal-dual gap, a measure of suboptimality, is defined as
$$
\mathcal{G}(u, p) := \left( \mathcal{J}(A u)+\mathcal{H}(u) \right) + \left( \mathcal{J}^*(p)+\mathcal{H}^*(-A^* p) \right),
$$
which vanishes if and only if $(u, p)$ is a saddle point of $\mathcal{L}$.

\begin{example}[Dual ROF]
To show how duality can simplify or offer new perspectives, we can derive the dual of the ROF problem from \Cref{eg:ROF_opt}. Let $A=\nabla$, $\mathcal{J}(z)=\alpha\|z\|_{2,1}$, and $\mathcal{H}(u)=\frac{1}{2}\|u-y\|_2^2$.
The convex conjugate of $\mathcal{J}$ is $$
\mathcal{J}^*(p)=\chi_{\left\{\|\cdot\|_{2, \infty} \leq \alpha\right\}}(p)= \begin{cases}0 & \text { if }\left|p_{i, j}\right|_2 \leq \alpha \quad \forall i, j, \\ +\infty & \text { otherwise }\end{cases}.
$$
The conjugate of $\mathcal{H}$ is $\mathcal{H}^*(q) = \frac{1}{2}\|q+y\|_2^2 - \frac{1}{2}\|y\|_2^2$.
Substituting into the dual formulation, we get:
$$
\begin{aligned}
& \max _p-\mathcal{J}^*(p)-\left(\frac{1}{2}\left\|\nabla^* p\right\|_2^2-\left\langle\nabla^* p, y\right\rangle\right) \\
& =-\min _p\left(\mathcal{J}^*(p)+\frac{1}{2}\left\|\nabla^* p-y\right\|_2^2\right)+\frac{1}{2}\|y\|^2.
\end{aligned}
$$
So the dual ROF problem is equivalent to solving:
$$
\min _p\left\{\frac{1}{2}\left\|\nabla^* p-y\right\|_2^2 : \left\|p_{i, j}\right\|_2 \leq \alpha \text{ for all } i, j\right\}.
$$
This dual problem is a constrained least-squares problem, which can be easier to solve than the primal non-smooth problem. From the optimality conditions of the saddle-point problem, we also have the relationship $u=y-\nabla^* p$ connecting the primal and dual solutions.

\end{example}
With these tools, we can now introduce several iterative algorithms designed for non-smooth convex optimization.
\subsection{Proximal descent}
The proximal descent algorithm is a foundational iterative method for minimizing non-smooth convex functions, generalizing the idea of gradient descent by leveraging the proximal operator. Starting from an initial guess $u_0$, it generates a sequence of iterates $u^{(k)}$ according to:
$$
u^{k+1} = \operatorname{prox}_{\tau J}(u^k) = (I + \tau \partial J)^{-1}(u^k).
$$
If $J$ is differentiable, this corresponds to an implicit gradient step $$u^{k+1}=u^k-\tau \nabla J\left(u^{k+1}\right).$$
The iterate $u^{k+1}$ is the unique minimizer of the proximal subproblem $J(v) + \frac{1}{2\tau}\|v-u^k\|_2^2$.
This algorithm can be interpreted as an explicit gradient descent step on the Moreau--Yosida regularization of $J$. The \textit{Moreau--Yosida regularization} of $J$ (or \textit{envelope}) with parameter $\tau >0$ is:
$$
J_\tau(\bar{u}):=\min _v \left( J(v)+\frac{\|v-\bar{u}\|_2^2}{2 \tau} \right).
$$
It can be shown that $J_\tau$ is continuously differentiable (even if $J$ is not) with gradient:
$$
\nabla J_\tau(\bar{u})=\frac{\bar{u}-\operatorname{prox}_{\tau J}(\bar{u})}{\tau}.
$$
Thus, the proximal descent update $u^{k+1} = \operatorname{prox}_{\tau J}(u^k)$ can be rewritten as $u^{k+1} = u^k - \tau \nabla J_\tau(u^k)$, which is an explicit gradient descent step on the smoothed function $J_\tau$.

\subsection{Forward-Backward Splitting}
To solve optimization problems structured as the sum of a smooth function and a non-smooth (but proximally tractable) function, the Forward-Backward Splitting algorithm offers an effective iterative approach. Consider the problem of minimizing the sum of two convex, proper, and lsc. functions:
$$
\min _u \mathcal{J}(u)+\mathcal{H}(u),
$$
where $\mathcal{J}$ is ``simple'' (its prox is easily computable) and $\mathcal{H}$ is differentiable with a Lipschitz continuous gradient (Lipschitz constant $L_{\mathcal{H}}$). The Forward-Backward splitting algorithm, also known as the proximal gradient algorithm, combines an explicit gradient descent step on $\mathcal{H}$ (forward step) and an implicit proximal step on $\mathcal{J}$ (backward step):
$$
u^{k+1} = \operatorname{prox}_{\tau \mathcal{J}}(u^k - \tau \nabla \mathcal{H}(u^k)).
$$
A point $u$ is a minimizer of the composite objective if and only if it is a fixed point of this iteration, which corresponds to the optimality condition $0 \in \nabla \mathcal{H}(u)+\partial \mathcal{J}(u)$. If the step size $\tau$ satisfies $0 < \tau \leq 1/L_{\mathcal{H}}$, the iterates $u^k$ converge to a minimizer.
\subsection{Primal-Dual Hybrid Gradient descent}\label{sec:pdhg}
For more complex structured problems, particularly those involving a linear operator $A$ coupling terms (e.g., $\mathcal{J}(Au)$), primal-dual algorithms iteratively seek a saddle point of an associated Lagrangian, which corresponds to a solution of the original problem. Consider now problems of the form 
$$
\min _u \mathcal{J}(A u)+\mathcal{H}(u),
$$
where $\mathcal{J}, \mathcal{H}$ are convex, lsc and simple, and $A$ bounded and linear. Primal-dual hybrid (PDHG) gradient descent alternates between proximal descent in the primal variable $u$ and ascent in the dual variable $p$ for the corresponding saddle-point problem
$$
\max _p \inf _u\langle p, A u\rangle-\mathcal{J}^*(p)+\mathcal{H}(u)
$$
via updating the primal and dual variables in an  alternating fashion:
$$
\begin{aligned}
u^{k+1} & =\operatorname{prox}_{\tau \mathcal{H}}\left(u^k-\tau A^* p^k\right) \\
p^{k+1} & =\operatorname{prox}_{\sigma \mathcal{J}^*}\left(p^k+\sigma A\left(2 u^{k+1}-u^k\right)\right)
\end{aligned}
$$
This algorithm is closely related to other optimization techniques like the augmented Lagrangian method and the alternating direction method of multipliers (ADMM) \citep{arrow1958studies,pock2009algorithm,esser2010general}.

\section{Regularizer Zoo}
Beyond the ROF model, TV regularization has found applications in a multitude of image processing tasks such as deblurring, inpainting, and segmentation. Each distinct application typically gives rise to a unique PDE, thereby presenting specific analytical and numerical challenges. Since the inception of these TV-based approaches, considerable research effort has been dedicated to the comprehensive analysis of these models, encompassing both their numerical and analytical properties. The analytical and numerical challenges posed by these problems have driven significant advancements in optimization and PDE theory, leading to the development of sophisticated algorithms and a deeper understanding of the underlying mathematical structures. Here we summarize key aspects of this extensive body of work, though it represents only a selection from the broader literature.
 
\paragraph{Analytical Properties:} 
\begin{itemize}[nosep]
    \item Function Space: 
    The natural function space for TV regularization is the space of functions of bounded variation, $B V(\Omega)$. This space is non-reflexive, necessitating the use of specialized compactness properties for analysis 
\citep{ambrosio2000functions,ambrosio1990metric,de1988nuovo,MR1783032}.
    \item Stability: 
    Novel metrics like Bregman divergences can be employed for deriving stability estimates in TV-regularized problems \citep{burger2007error,hofmann2007convergence,schonlieb2009image}.

    \item  Non-differentiability: 
    The non-differentiability of the TV term requires tools from convex analysis, such as subgradients, and leads to the study of TV flow via differential inclusions or viscosity solutions \citep{chen1999uniqueness,ambrosio1996level,caselles2007discontinuity,alter2005evolution,caselles2006anisotropic,bellettini2002total,bellettini2002total_2,paolini2003relaxed,novaga2005regularity,bellettini2006global}.
    Analysis often draws upon geometric measure theory 
    \citep{federer1996applications,federer2014geometric,allard1986geometric,allard2008total}.
\end{itemize}
\paragraph{Numerical Properties:}
\begin{itemize}[nosep]
    \item Non-smooth Optimization: 
    Standard optimization frameworks for smooth, strictly convex problems are not directly applicable to TV regularization. New analysis is based on algorithms involving generalized Lagrange multipliers, Douglas--Rachford splitting \citep{lions1979splitting}, iterative thresholding algorithms \citep{chambolle2015convergence, chambolle2004algorithm, combettes2008proximal, combettes2005signal, daubechies2004iterative}
    and semi-smooth Newton methods \citep{hintermuller2010semismooth, hintermuller2003semi}.
    \item Scalability:
    Large-scale non-smooth convex minimization problems involving TV regularization often scale poorly. Accelerations are achieved through preconditioning, splitting approaches, partial smoothness, and stochastic optimization \citep{chan1990circulant,chan1989toeplitz, fornasier2009subspace, fornasier2010convergent, afonso2010fast, bioucas2007new, figueiredo2007gradient, beck2009fast,cevher2014convex,liang2014local,bredies2015preconditioned,chambolle2018stochastic}.
    \item Non-smooth and Nonlinear Problems:
    Further challenges arise when using non-smooth regularization for nonlinear inverse problems. From an optimization perspective, the combination of non-smoothness and non-convexity opens yet another chapter in numerical analysis and optimization with only partly satisfying results. Seminal contributions include the work of \citet{kaltenbacher2008iterative} and \citet{bachmayr2009iterative} on computational solutions for nonlinear inverse problems, and several works on non-smooth and non-convex optimization \citep{attouch2010proximal, valkonen2014primal, bolte2014proximal,chizat2018global, driggs2021stochastic, benning2021choose}.
\end{itemize}

\subsection{Regularizer Zoo: Electric boogaloo}

Real-world inverse problems often present unique challenges, such as diverse noise, intricate structures, and specific modeling requirements, necessitating tailored regularizers for optimal reconstruction.  Classically, this has been addressed by designing and refining handcrafted models encouraging desirable properties in the reconstructed image, such as smoothness, sparsity, or adherence to specific geometric features. Over time, this has lead to a diverse array of specialized regularizers, based on:

\begin{itemize}[nosep]
    \item Multi-resolution analysis, wavelets \citep{mallat1999wavelet,daubechies1992ten,vonesch2007generalized,unser2000fractional,dragotti2003wavelet,kutyniok2012introduction,foucart2013invitation,fornasier2012wavelet};
    \item Other Banach-space norms, e.g. Sobolev norms, Besov norms, etc. \citep{saksman2009discretization,lassas2009discretization};
    \item  Higher-order total variation regularization \citep{osher2003image, chambolle1997image, setzer2011infimal, bredies2010total}, 
    as well as higher-order PDEs, Euler elastica \citep{masnou1998level,shen2003euler,bertozzi2006inpainting}; 
    \item Non-local regularization \citep{gilboa2009nonlocal,buades2011non};
    \item Anisotropic regularization \citep{weickert1998anisotropic};
    \item Free-discontinuity problems; \citep{mumford1989optimal, carriero1992plastic}.
    \item and mixtures of the above and others\dots
\end{itemize}
\noindent This ``regularizer zoo'' highlights the diversity of approaches developed to address the specific challenges posed by different image reconstruction problems. The selection of an appropriate regularizer requires careful consideration of the image properties, the degradation process, and the desired characteristics of the reconstruction.
\section{Limitations and move towards data-driven}
While knowledge-driven regularization has significantly advanced image reconstruction, its effectiveness is inherently limited by our ability to accurately model the complex structures present in real-world images. Even with sophisticated mathematical models and physical principles from non-linear PDEs, variational formulations, and multi-resolution or sparsity constraints, our understanding of image formation and degradation remains incomplete. Consider the remarkable ability of the human brain to denoise and interpret images under challenging conditions — a feat far surpassing current knowledge-driven approaches.

This limitation highlights the need for a paradigm shift towards data-driven reconstruction methods, which leverage the power of overparameterized models like support vector machines and neural networks \citep{Goodfellow-et-al-2016}. These models, trained on vast amounts of data, can learn intricate patterns and relationships that may be difficult to capture through explicit mathematical formulations. The next section will delve into existing paradigms and challenges of data-driven reconstruction, exploring how it can complement and even surpass knowledge-driven methods in the quest for accurate and robust image recovery.  Furthermore, we will investigate the emerging field of hybrid approaches that combine the strengths of both data-driven and knowledge-driven techniques, potentially leading to a new generation of image reconstruction algorithms that push the boundaries of performance and applicability.

%% file: chapter3.tex
\chapter{Data-Driven Approaches to Inverse Problems}\label{cap3}

This chapter explores the paradigm shift from knowledge-driven to data-driven approaches for solving inverse problems in imaging. As we established in previous chapters, traditional knowledge-driven methods, while powerful, are fundamentally limited by our ability to accurately model the complexities of real-world images.

\noindent It is important to preface this section by acknowledging the exceptionally rapid evolution of the field of deep learning, particularly in imaging. The content presented herein largely reflects the understanding and prominent methods as of 2023, the time of the original CIME lectures. Consequently, while the foundational concepts discussed remain relevant, the field has likely seen further advancements since.
\subsection{Knowledge-Driven vs. Data-Driven Models}\label{sec:knowledge_vs_data}

\begin{figure}[!h]
    \centering
    \begin{subfigure}[t]{0.24\textwidth} %
        \centering
        \begin{tikzpicture}[spy using outlines={
            rectangle, 
            red, 
            magnification=3.5,
            size=1.2cm, 
            connect spies}]
            \node {\includegraphics[width=1.2in]{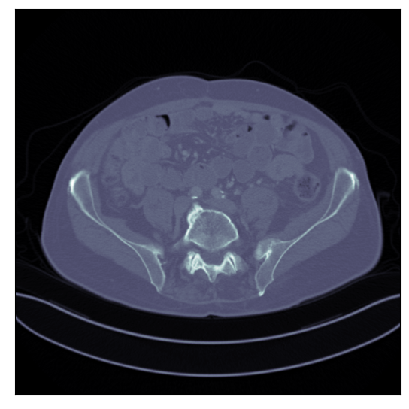}}; %
            \spy on (0.05,-0.50) in node [left] at (2.00,-1.0);
            \spy on (-0.3,0.65) in node [left] at (2.00,1.0);  
        \end{tikzpicture}
        \caption{Ground-truth}
        \label{sfig:kgp_truth_subset}
    \end{subfigure}
    \hfill %
    \begin{subfigure}[t]{0.24\textwidth}
        \centering
        \begin{tikzpicture}[spy using outlines={
            rectangle, 
            red, 
            magnification=3.5,
            size=1.2cm, 
            connect spies}]
            \node {\includegraphics[width=1.2in]{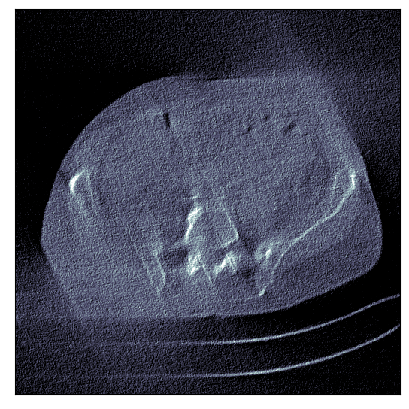}}; %
            \spy on (0.05,-0.50) in node [left] at (2.00,-1.0);
            \spy on (-0.3,0.65) in node [left] at (2.00,1.0);  
        \end{tikzpicture}   
        \caption{FBP: 21.61 dB, 0.17}
        \label{sfig:kgp_fbp_subset}
    \end{subfigure}
    \hfill
    \begin{subfigure}[t]{0.24\textwidth}
        \centering
        \begin{tikzpicture}[spy using outlines={
            rectangle, 
            red, 
            magnification=3.5,
            size=1.2cm, 
            connect spies}]
            \node {\includegraphics[width=1.2in]{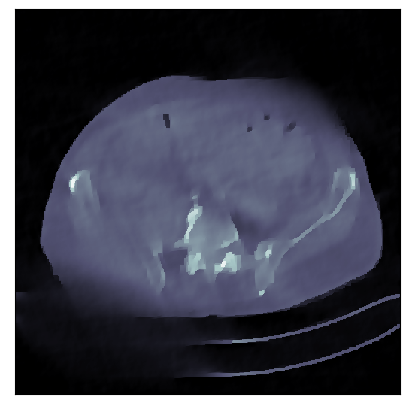}}; %
            \spy on (0.05,-0.50) in node [left] at (2.00,-1.0);
            \spy on (-0.3,0.65) in node [left] at (2.00,1.0);  
        \end{tikzpicture}    
        \caption{TV: 25.74 dB, 0.80}
        \label{sfig:kgp_tv_subset}
    \end{subfigure}
    \hfill %
    \begin{subfigure}[t]{0.24\textwidth}
        \centering
        \begin{tikzpicture}[spy using outlines={
            rectangle, 
            red, 
            magnification=3.5,
            size=1.2cm, 
            connect spies}]
            \node {\includegraphics[width=1.2in]{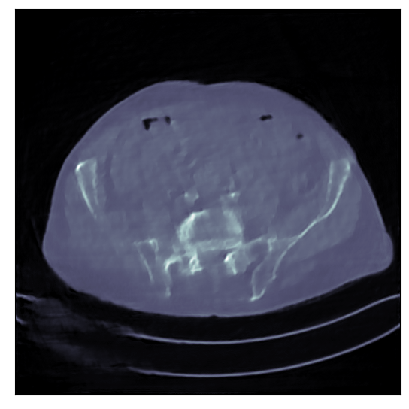}}; %
            \spy on (0.05,-0.50) in node [left] at (2.00,-1.0);
            \spy on (-0.3,0.65) in node [left] at (2.00,1.0);  
        \end{tikzpicture}
        \caption{LPD: 29.51 dB, 0.85}
        \label{sfig:kgp_lpd_subset}
    \end{subfigure}
    \caption{
    Limited angle CT reconstruction: Heavily ill-posed problem. Deep Learning cannot do magic and also hits boundaries of what is mathematically possible. A fully learned method LPD (\Cref{sec:unroll}) in \subref{sfig:kgp_lpd_subset} begins hallucinating, as highlighted in red boxes, despite resulting in better performance metrics (here PSNR and SSIM).}
\label{fig:badexample}
\end{figure}

\begin{figure}[htb!]
    \centering
    \includegraphics[width=0.7\textwidth]{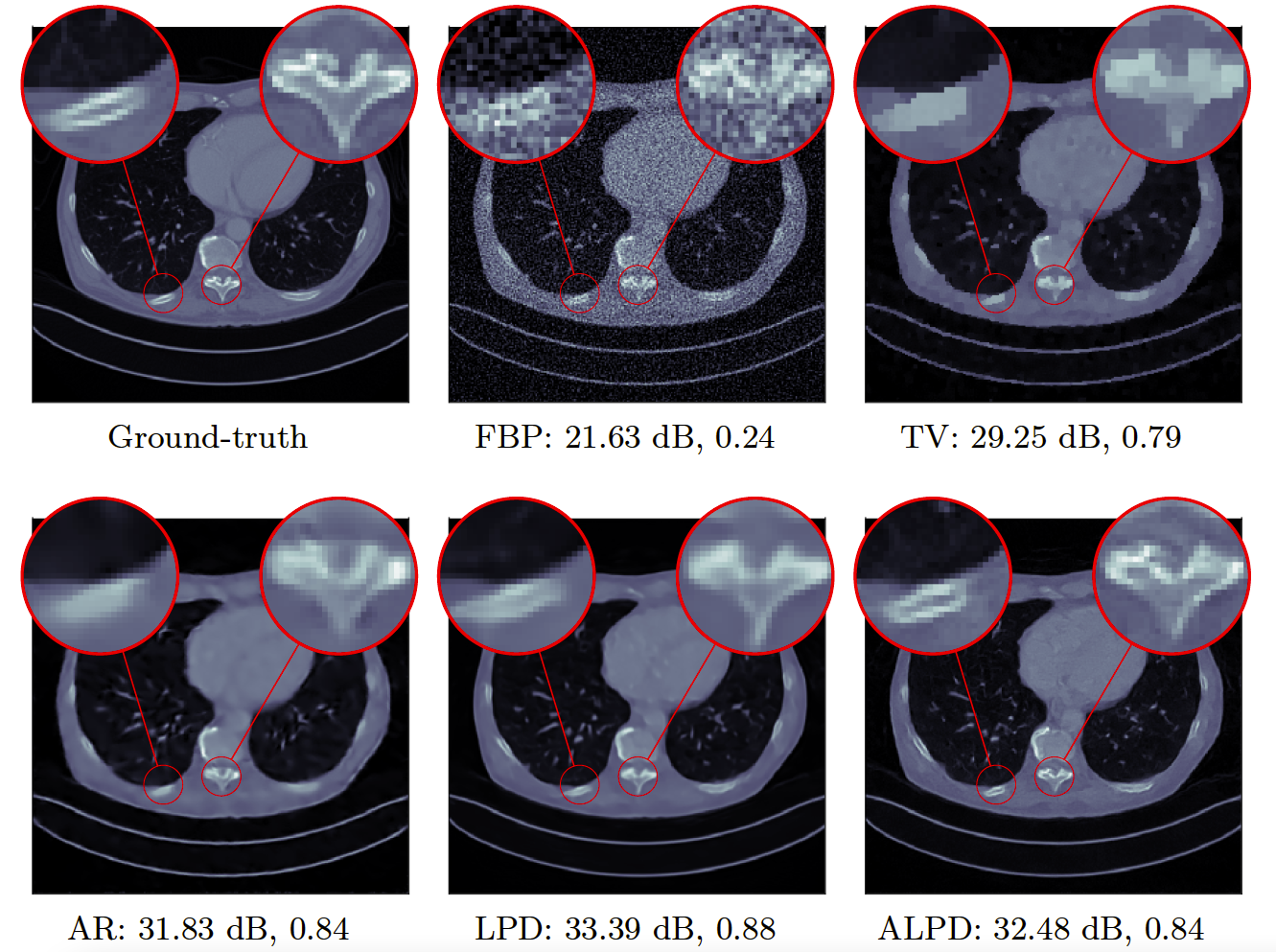}
    \caption{Sparse view CT reconstruction: top row is based on mathematical/handcrafted models; bottom row is using novel deep learning based models. For this problem, deep learning methods result in both improved metrics (here PSNR and SSIM) and visually better reconstructions. Photo courtesy of Subhadip Mukherjee \citep{subho2023}.}
    \label{fig:goodexample}
\end{figure}

\paragraph{Knowledge-Driven Models}
Knowledge-driven models are rooted in mathematical principles and domain expertise. These models are often interpretable and offer theoretical guarantees, but their performance is fundamentally limited by two critical factors: the accuracy of the underlying mathematical model, and the designer's ability to capture complex system behaviors. Such models often struggle with complex or unknown noise patterns and complex image structures that defy straightforward mathematical representation. Even if such complexities could be modeled, practical limitations arise due to computational constraints: while we might develop highly detailed models, the sheer complexity may render them computationally infeasible to use effectively.

The deep learning revolution of the 2010s fundamentally transformed our approach to complex imaging tasks by challenging traditional modeling paradigms. As computational power and data availability dramatically increased, neural networks demonstrated their ability to learn intricate representations directly from massive datasets, often outperforming carefully crafted mathematical models. Consequently, given the abundance of image data available today, a natural question emerges: why meticulously handcraft models when we can potentially derive effective priors simply by providing sufficient data to overparameterized models?

\paragraph{Data-Driven Models}
Data-driven models, in contrast to knowledge driven approaches, directly extract information from data. While traditional knowledge-driven methods often involve some degree of learning (e.g., parameter estimation), they typically rely on models with a limited number of parameters. Unlike knowledge-driven models, data-driven approaches have a significant number of parameters and leverage large datasets to learn complex patterns and relationships without explicit mathematical modeling.

\noindent Deep learning exemplifies this paradigm, using extensive computational resources to train highly flexible, over-parameterized neural networks that can adapt to diverse imaging tasks and datasets, especially in high dimensions, while remaining computationally efficient. An example in the context of inverse problems is shown in \Cref{fig:goodexample} -- learned methods consistently and significantly outperform knowledge driven methods like TV regularization. Despite their power, these models often sacrifice interpretability and demand substantial training data to achieve good performance. An example of this in the context of inverse problems is shown in \Cref{fig:badexample} - for a significantly ill-posed problem, fully learned methods begin hallucinating, despite resulting in better performance metrics.

\noindent We note however, that derivation of models directly from data is by no means exclusive to deep neural networks; in fact, such approaches predate and extend beyond them, constituting a rich methodological landscape within machine learning and signal processing. To illustrate, classical learning techniques have long been employed to explore data-driven regularization models, for example:

\begin{itemize} 
    \item \textbf{Sparse Coding and Dictionary Learning}: These methods aim to find a sparse representations of signals as linear combinations of a few elementary atoms from a dictionary, which itself can be learned from data. Approaches optimize signal representation through a minimization problem that balances data fidelity and sparsity:
    $$\min_{\gamma, \phi}\left\|A\left(\sum_i \gamma_i \phi_i\right)-y\right\|_2^2+\|\gamma\|_1.$$
    Some examples include 
    \citet{elad2006image, aharon2006k,mairal2009online,rubinstein2010dictionaries, moreau2016understanding, chandrasekaran2011rank, devore2007deterministic, fadili2009image, mallat1993matching,elad2006image,rubinstein2009double,papyan2017convolutional, peyre2009sparse}.
    
\item \textbf{Black-Box Denoiser Methods}: These techniques integrate powerful, often pre-existing, denoising algorithms as implicit priors within iterative reconstruction schemes, without requiring explicit knowledge of the denoiser's internal structure. Some examples include the Plug-and-Play Prior method \citet{venkatakrishnan2013plug, wei2020tuning} and Regularization by Denoising \citet{romano2017little,terris2020building}:
$$\min_u\; \mathcal{D}(A(u), y)+\alpha R(u), \quad \text{with } R(u)=\langle u, u-\Lambda(u)\rangle, \quad \Lambda: X \rightarrow X \text{ denoiser}.$$
\item \textbf{Bilevel Optimization (since early 2000s)}: This class of methods addresses the challenge of learning model parameters by formulating a nested optimization problem, where an outer problem optimizes parameters used in an inner image reconstruction or processing task. 
    $$\min_\lambda F(u_\lambda) \quad \text{s.t. } u_\lambda=\underset{u}{\operatorname{argmin}}\; R(\lambda, u)+\mathcal{D}(A(u), y).$$ 
Some examples include \citet{calatroni2017bilevel,kunisch2013bilevel, de2016structure, haber2009numerical,langer2017automated,horesh2010optimal}.
\end{itemize}
\noindent The key distinctions between the paradigms emerge not just in methodology, but in their philosophical approach: knowledge-driven models seek to understand through explicit modeling assumptions, while data-driven models pursue understanding through statistical learning and pattern recognition.

The rest of this section will discuss deep learning more generally, will present a number of recent approaches within the data-driven paradigm, progressively advancing towards methodological frameworks that intersect the two paradigms, exemplified through methods using deep neural networks as regularizers.

\subsection{The Black Box of Deep Learning}
Deep learning has shown remarkable success in various fields, but comes with significant  challenges and limitations associated with these approaches. We refer the interested reader to the review by \citet{grohs2022mathematical}. State of the art deep neural networks have `too many' degrees of freedom:
\begin{itemize}[nosep]
    \item millions of free parameters, i.e., the parameter space $\Theta$ is super high-dimensional;
    \item complex concatenation of diverse mathematical constructs (convolutions, activations, attention, skip connections, normalization, dropout, ...);
    \item high-dimensional and non-convex optimization problems.
\end{itemize}
The usefulness of the resulting model is influenced by all of these model design parameters, as well as the quality of the training set, and optimization approach used. This in turn, makes it difficult to understand their internal workings and interpret their outputs. This makes deep learning models effectively a ``\textit{black box}''. More precisely, the resulting issues are: 
\begin{itemize}[nosep]
    \item \textbf{Lack of Interpretability}: It is difficult to understand why a deep learning model produces a particular output. This can make it challenging to identify biases, errors, or limitations in the model.   

    \item \textbf{Safety Concerns}: In applications where safety is critical, such as medical imaging, or autonomous driving, the lack of interpretability can raise concerns about the reliability and trustworthiness of deep learning models.

    Particularly, for the problem of CT reconstruction: how can we reliably say that a deep learning model did not introduce or obscure cancerous tissue during reconstruction?
    
    \item \textbf{Limited Design Principles}: Due to their complexity, there is no systematic way to design deep learning models. Their development often involves trial and error, making it difficult to guarantee optimal performance or generalize to new tasks.  
\end{itemize}
Despite these, deep learning offers interesting opportunities for inverse problems thanks to its ability to produce highly accurate and computationally fast solutions. However, to fully leverage the potential of deep learning,  it is essential to integrate these techniques with established mathematical modeling principles. This synergy is crucial to ensure predictable and reliable (in a certain sense) behavior in the resulting solutions, as oftentimes interpretability goes hand in hand with mathematical guarantees, but often comes at an expense in either reduced performance or computational intractability. The current main goal of the field thus is in trying to find the sweet spot between computational power and mathematical guarantees.

\begin{example}[Neural Networks in a Nutshell]

A neural network can be formally defined as a mapping:
$$
\begin{aligned}
\Psi: X \times \Theta & \rightarrow Y \\
(x, \theta) & \mapsto z^K ,
\end{aligned}
$$
where in above, $X$ is the input space, $Y$ is the output space, $x \in X$ is the input data, $z^K \in Y$ is the output of the network, $\Theta = (\Theta^0, ..., \Theta^{K-1})$ represents the parameter space, and $\Theta^k$ denoting the parameter space of the $k$-th layer. The network's internal operations are characterized by a sequence of layer-wise transformations:
$$
\begin{aligned}
z^0 & =x \in X \\
z^{k+1} & =f^k\left(z^k, \theta^k\right), \quad k=0, \ldots, K-1
\end{aligned}
$$
where $z^k \in X^k$ represents the feature vector at the $k$-th layer, with $X^k$ being the corresponding feature space,
and $f^k: X^k \times \Theta^k \rightarrow X^{k+1}$ is the non-linear transformation at the $k$-th layer, parameterized by $\theta^k$.
A common choice for $f^k$ is an affine transformation followed by an element-wise non-linear activation function:
$$
f^k(z)=\sigma\left(W^k z+b^k\right),
$$
where $W^k$ is a weight matrix (for imaging tasks often represented by a convolution operator), $b^k$ is a bias vector, $\sigma$ is an element-wise non-linear activation function (e.g., \texttt{ReLU}, \texttt{tanh}).

\noindent The training process aims to optimize the network parameters $\theta$ by minimizing a loss function $L_n$ over a given dataset $\left\{\left(x_n, c_n\right)\right\}_n$, often with an added regularization term $R(\theta)$ (this time to regularize the \emph{training}):
$$
\min _{\theta \in \Theta} \frac{1}{N} \sum_{n=1}^N L_n\left(\Psi\left(x_n, \theta\right), c_n\right)+R(\theta)
$$
This generic framework can be adapted and applied to various mathematical imaging tasks, such as image classification, segmentation, and reconstruction, by appropriately defining the network architecture, loss function, and training data.
\end{example}

\section{Learned Iterative Reconstruction Schemes}\label{sec:unroll}

This section explores learned iterative reconstruction schemes, a class of end-to-end deep learning methods for solving inverse problems. 
Some good reviews on the topic include \citet{arridge2019solving,mccann2017convolutional}, with many works in the literature proposing similar approaches \citep{gregor2010learning,sun2016deep,meinhardt2017learning,putzky2017recurrent,adler2017solving,adler2018learned,hammernik2018learning,hauptmann2018model,de2022deep,gilton2019neumann,bubba2021deep}. These methods draw inspiration from classical iterative algorithms, where individual iterative steps are replaced or augmented with neural networks. The core idea is to ``unroll'' a fixed number of iterations of an optimization algorithm and learn parts of this unrolled scheme from data. This approach often begins by considering a gradient descent update rule for a variational problem and then introducing parameterized blocks within these iterations. The parameters of these blocks are subsequently optimized over a fixed number of steps using supervised data.
\begin{figure}
    \centering
\includegraphics[width=0.8\textwidth]{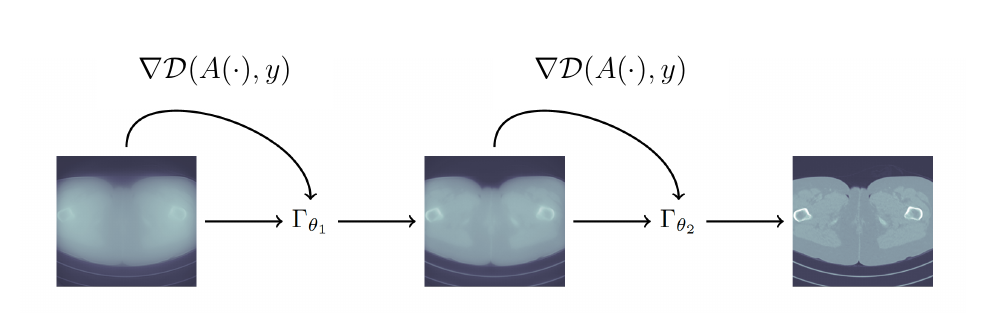}
    \caption{Learned Iterative Schemes Schematic.
    }
    \label{fig:lgd}
\end{figure}
\noindent The general concept can be illustrated by comparing it with standard gradient descent. Given an initial guess $u^0 \in X$, gradient descent follows a sequence of steps as 
\[
u:=\left(\operatorname{Id}-\eta\nabla \mathcal{D}(A (\cdot), y)\right)^N u_0.
\]
Learned iterative schemes generalize this by parameterizing these steps individually:
 $$
u:=\left(\Lambda_{\Theta_N} \circ \cdots \circ \Lambda_{\Theta_1}\right)\left(u^0\right) .
$$
Each $\Lambda_{\Theta_k}$ can be viewed as a residual layer in a neural network $\Psi_{\Theta}(y)$ with $N$ layers, which reconstructs $u$ from $y$. These schemes are typically derived from an iterative method designed to solve a variational regularization problem and several variations of learned iterative schemes exist, each with a different formulation for the layers $\Lambda_{\Theta_k}$. Writing the steps as 
$$
u^{k+1}=\Lambda_{\theta_k}\left(u^k, {A}^*\left({A} u^k-y\right)\right) \quad \text { for } k=0, \ldots, N-1,
$$
for some neural networks $\Lambda_{\theta_k}: X \times X \rightarrow X$, the main examples (described in more detail below) are
\begin{alignat}{2}
\Lambda_\theta(u, h)&:=u+\Gamma_\theta(h) &\quad \text { (original learned gradient) } \\
\Lambda_\theta(u, h)&:=u-h+\Gamma_\theta(u) &\quad \text { (variational networks) } \\
\Lambda_\theta(u, h)&:=\Gamma_\theta(u-h) &\quad \text { (plug-and-play; learned proximal) } 
\end{alignat}
for some neural network $\Gamma_\theta: X \rightarrow X$ with an architecture that does not involve data or the forward operator (or its adjoint), which only enter into the evaluation of $h = \mathcal{D}(A(\cdot), y)$. An illustration is shown in \Cref{fig:lgd}. Parameters are then learned from supervised data.

\paragraph{Variational Networks.} First proposed by \citet{hammernik2018learning, kobler2017variational}, they represent one of the earliest learned iterative schemes and are notable for their explicit connection to variational regularization. These networks are inspired by variational regularization models where the regularization functional incorporates parameterizations that extend beyond traditional handcrafted regularisers like TV. More specifically, their development draws from the Fields of Experts (FoE) model  and, by extension, conditional shrinkage fields \citep{schmidt2014shrinkage}, which allow for adaptive parameters across iterations. Variational Networks are defined by unrolling an iterative scheme designed to minimize a functional comprising both a data discrepancy component and a regularizer. Such networks can be interpreted as performing block incremental gradient descent on a learned variational energy or as learned non-linear diffusion.

\begin{figure}
    \centering
    \includegraphics[width=\textwidth]{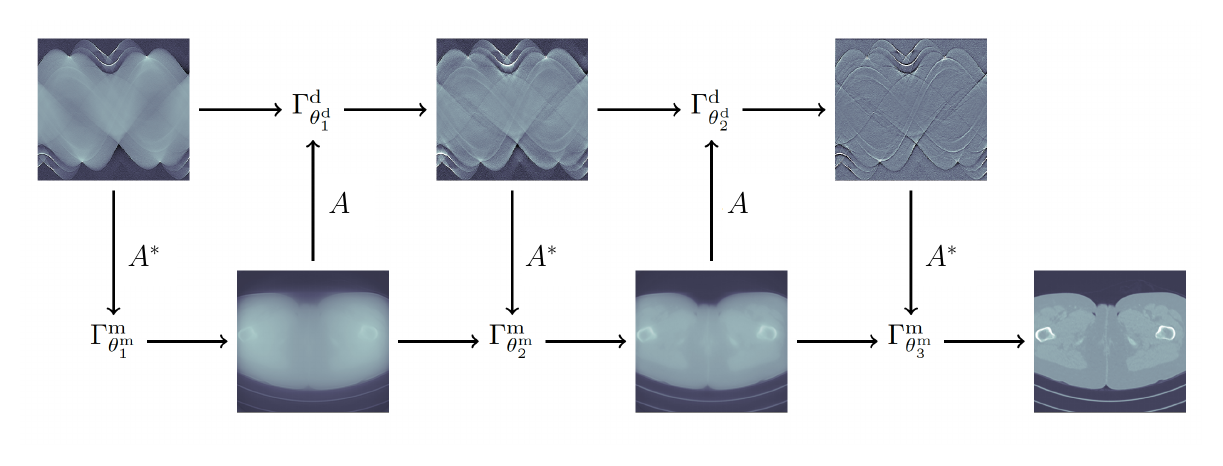}
    \caption{Learned Primal Dual Schematic.}
    \label{fig:lpd}
\end{figure}

\paragraph{Learned Gradient, Proximal, Primal Dual.} These are all further extensions of the idea introducing increasingly more freedom to parameters \citep{adler2017solving}, with \citet{adler2018learned} generalizing the steps even further to include steps in both primal (image) and dual (measurement) spaces, inspired by the primal-dual hybrid gradient from \Cref{sec:pdhg}. This can be summarized with a more general parameterization of the gradient steps, potentially including some $\mathcal{R}$ regularization information, e.g. using TV \citep{kiss2025}. An illustration is shown in \Cref{fig:lpd}.
$$
\Lambda_{\Theta_k}:=\Gamma_{\Theta_k}\left(u, y, A^* y, A u, \nabla \mathcal{R}(u)\right) \quad \text { for }(u, \tilde{y}) \in X \times Y,
$$
Empirical evidence suggests such approaches result in models that are easier to train and demonstrate excellent \textit{reconstruction quality} for mildly ill-posed inverse problems and offer considerable versatility, for instance in task-adapted reconstruction \citep{adler2022task,lunz2018adversarial}. While intuitively appealing, this approach can lose connection to the original variational problem, leading to a lack of theoretical guarantees. Despite this, this generic recipe has been used to propose numerous methods, leveraging various iterative algorithms such as gradient descent, proximal descent, and primal-dual methods.

\subsection{Limitations and Challenges}
Despite their promising performance, learned iterative schemes are still often perceived as ``black boxes'' due to several outstanding challenges.
\begin{itemize}[nosep]
    \item \textbf{Limited Theoretical Understanding:} There is a general lack of rigorous analysis regarding their well-posedness and regularization properties, with few exceptions \citep{hertrich2021convolutional,sun2021scalable,gilton2021model}. This theoretical gap contributes to an ongoing debate within the field concerning the applicability of deep learning to inverse problems arising in critical fields due to a lack of robustness \citep{genzel2022solving}. 
    
    \item \textbf{Interpretability:} 
    The learned operators within these schemes frequently lack a clear mechanistic explanation, making it difficult to fully understand their behavior. While some asymptotic properties are known, such as the convergence of learned iterative schemes with an $\ell_2$ loss to the conditional mean under infinite data conditions \citep{adler2017solving}, a deeper, more general understanding remains elusive.
    \item \textbf{Data Consistency} with the measurements is not always guaranteed in the final reconstructions, although specific approaches exist, attempting to address this limitation, e.g. \citep{schwab2019deep}.
    \item \textbf{Supervised Training:} These methods typically require large amounts of supervised data, which can be challenging to obtain in practice.
    \item \textbf{Convergence:} Iterating beyond the number of training steps may not guarantee convergence.  Learned iterative schemes are trained for a fixed number of iterations (typically $\leq 20$) due to computational constraints, and the reconstruction deteriorates if more iterations are performed at test time. An example of such deterioration can be seen in bottom row of \Cref{fig:deq}. 
    \item \textbf{Computational Cost:}  Evaluating the forward and adjoint operators in each layer can be computationally memory expensive, hindering scalability.
\end{itemize}
While this may seem gloomy, many current research efforts are focused on addressing these limitations, including
\begin{itemize}[nosep]
    \item \textbf{Theoretical Analyses:} Investigating the convergence properties and approximation capabilities of learned iterative schemes, a particular example of which is the usage of equilibrium models discussed in the next subsection.
    \item \textbf{Efficient Implementations:} Exploring techniques like invertible neural networks \citep{rudzusika2021invertible,rudzusika20243d}, stochastic subsampling fo the forward operator \citep{tang2025iterative,tang2021stochastic}, and greedy training \citep{hauptmann2018model} to reduce computational cost.
\end{itemize}
\subsection{Deep Equilibrium Networks}

\begin{figure}
\includegraphics[width=\textwidth]{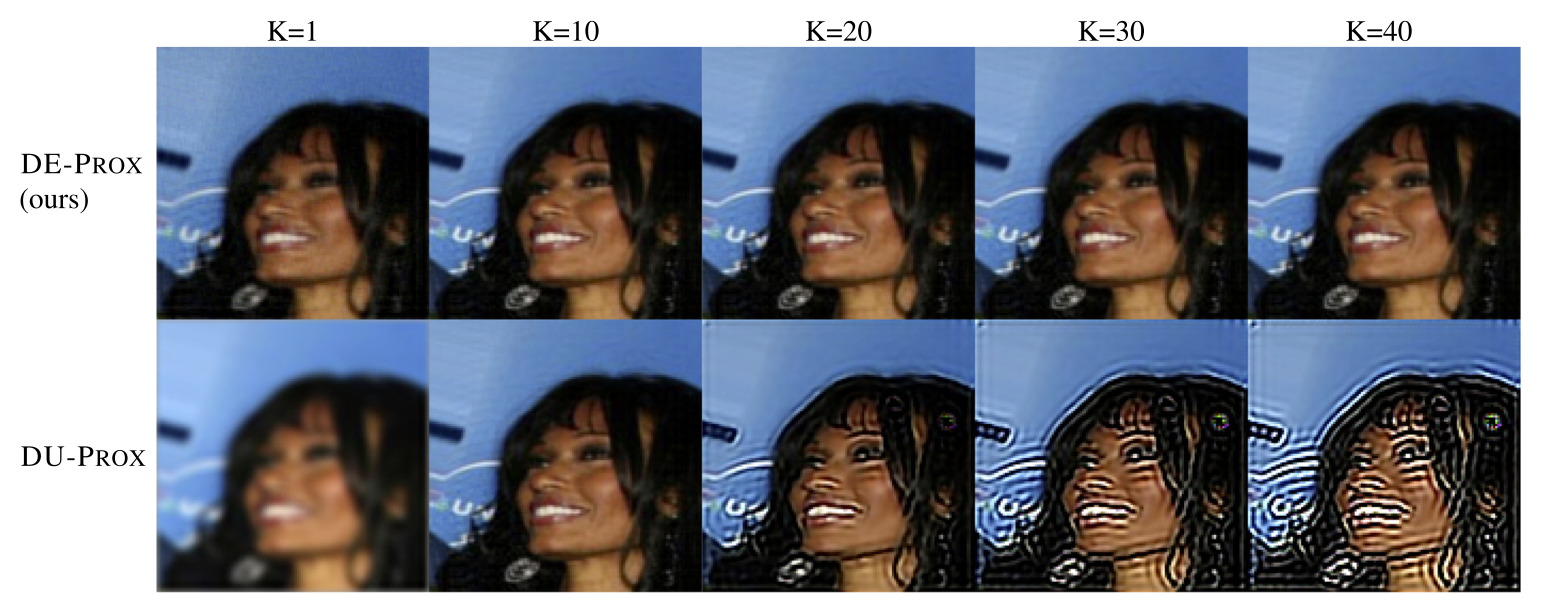}
    \centering
    \caption{Illustration of what artifacts appearing whenever learned operators are applied repeatedly without convergence guarantees. Example borrowed from \citep{gilton2021deep}.}
    \label{fig:deq}
\end{figure}
One promising avenue is the use of learned fixed point iterations, exemplified by deep equilibrium networks (DEQs).  These networks are designed such that the desired reconstruction is a fixed point of a learned operator 
$$
u=\Gamma_{\Theta}(u ; y) .
$$
This formulation naturally leads to iterative schemes that provably converge (under certain assumptions) to a fixed point as the number of iterations approaches infinity. For instance, consider a deep equilibrium gradient descent scheme where:
$$
\Gamma_{\Theta}(u ; y)=u+\eta A^*(y-A u)-\eta {R}_{\Theta}(u).
$$
Here, $A$ is the forward operator, $A^*$ is its adjoint, $\eta$ is a step size, and ${R}_{\Theta}$ is a trainable neural network representing a gradient of a learned regularizer. Note that these models are more general than gradient descent, as convergence can be ensured even when ${R}_{\Theta}$ is not a gradient of a function. To ensure convergence, $\Gamma_{\Theta}$ can be constrained to be a contraction mapping. This, once again, is not simply an academic exercise and has a significant effect in practice, guaranteeing convergence to a fixed point as showcased in \Cref{fig:deq} on top row, compared to iterate divergence for a non-constrained model, showcased on the bottom row. The following theorem provides sufficient conditions for convergence in the context of deep equilibrium gradient descent:
\begin{theorem}[\citep{gilton2021deep}]
Assume ${R}_{\Theta} - \operatorname{Id}$ is $\epsilon$-Lipschitz continuous and let $L=$ $\lambda_{\max }\left(A^* A\right)$ and $\mu=\lambda_{\min }\left(A^* A\right)$. If $0<\eta<1 /(L+1)$, then $\Gamma_{\Theta}$ fulfills 
$$\left\|\Gamma_{\Theta}(u ; y)-\Gamma_{\Theta}(\tilde{u} ; y)\right\| \leq \underbrace{(1-\eta(1+\mu)+\eta \epsilon)}_{=\gamma}\|u-\tilde{u}\|, \quad \forall u, \tilde{u} \in X.$$
Therefore, $\Gamma_{\Theta}$ is a contraction if $\epsilon<1+\mu$, and hence the iterates converge.
\end{theorem}
\begin{remark}
An interesting additional avenue with convergent learned iterative schemes is the ability to accelerate their convergence by increasing the memory of the iterations, i.e. by introducing dependency of each iteration from just the previous iterate only to a couple of previous iterates, for instance via Anderson acceleration as in \citet{gilton2021model}. While convergence to a fixed point is a desirable property, further investigation is needed to characterize the properties of this fixed point and its relationship to the solution of the underlying inverse problem, e.g. in analogy with \citet{obmann2023convergence}.
\end{remark}

\section{Learned Variational Models}

This section delves into learned variational models for inverse problems. While such models have a history in signal processing discussed in \Cref{sec:knowledge_vs_data}, application of overparameterized models represents a more recent development. 

\noindent The central idea is to leverage the well-established mathematical framework of variational methods, integrating deep learning while retaining theoretical guarantees and enforcing desirable structural properties on the neural networks. This approach allows us to combine the expressive power of deep learning with the stability and interpretability of variational methods.

\subsection{Learning the regularizer}
A particularly appealing strategy within the learned variational model paradigm focuses on \textbf{learning the regularizer} itself, while maintaining the classical variational framework. Consider the general form of a variational problem:
\begin{equation}\label{eq:varprob_learnreg}
\underset{u\in X}{\arg \min }\;\|A u-y\|_2^2+\alpha R(u)    
\end{equation}
where $y\in Y$ is the measured data, $A$ is a linear and bounded forward operator, $X$ and $Y$ are Banach spaces, $\alpha>0$ is a regularization parameter, and $R(u)$ is the regularizer.

\noindent Instead of learning the entire reconstruction mapping from $y$ to $u$ as in \Cref{sec:unroll}, this approach concentrates on learning a \textbf{data-adaptive regularizer} $R(u)$. The goal is for $R(u)$ to effectively capture prior knowledge about the desired solution, promoting reconstructions with desirable characteristics (e.g., ``good-looking'' images) while penalizing undesirable features. Learning the regulariser offers several advantages:
\begin{itemize}[nosep]
    \item \textbf{Interpretability:} The learned regulariser provides an explicit prior on the solution space.
    \item \textbf{Stability and Convergence:} Existing variational theory can be applied to analyze stability and convergence, e.g. through \Cref{thm:stab}.
    \item \textbf{Adaptability:} The regularization parameter $\alpha$ can be adjusted to accommodate different noise levels.
    \item \textbf{Incorporation of Forward Model and Noise Statistics:} The variational framework explicitly incorporates the forward model and noise statistics, e.g. ensuring data consistency.
\end{itemize}
In essence, this approach seeks to learn an ``image prior'' that is both data-driven and amenable to mathematical analysis. 

\paragraph{How to learn?} Given some parametric model $R_\theta$ for the regularizer, its parameters $\theta$ still need to be learned from data. Over the past decades, a variety of paradigms have been introduced for learning the regularizer given an image distribution. 
\noindent Analogous to \Cref{sec:unroll}, the direct approach is to train parameters such that the optimal solution in \Cref{eq:varprob_learnreg} minimizes the $\ell_2$ loss over training data. This results in the so called bilevel learning discussed in \Cref{sec:knowledge_vs_data}. While parameter hypergradients can be computed via implicit differentiation or unrolling, they can quickly become computationally infeasible, necessitating approximations. Alternatively, interpreting the problem as maximum-a-posteriori estimation (\Cref{sec:bayesian}), the prior can be learned directly from data. See \citet{habring2024neural} or \citet{dimakis2022deep} for an overview.

In what follows, we consider one specific approach for learning $\mathcal{R}_\theta$, which relies neither on the bilevel structure, nor on learning the whole prior. Instead, we view the regulariser as a ``classifier'' that distinguishes between desirable and undesirable solutions. This decouples the problem of learning the regulariser from the underlying variational problem.

\subsection{In-Depth: Adversarial regularization}

\begin{figure}[ht] %
    \centering
    \newlength{\commonFigHeight}
    \setlength{\commonFigHeight}{3.5cm} %

    \begin{subfigure}[b]{0.32\textwidth} %
        \centering
        \includegraphics[trim={2.5cm 6cm 15cm 6.5cm},clip, height=\commonFigHeight]{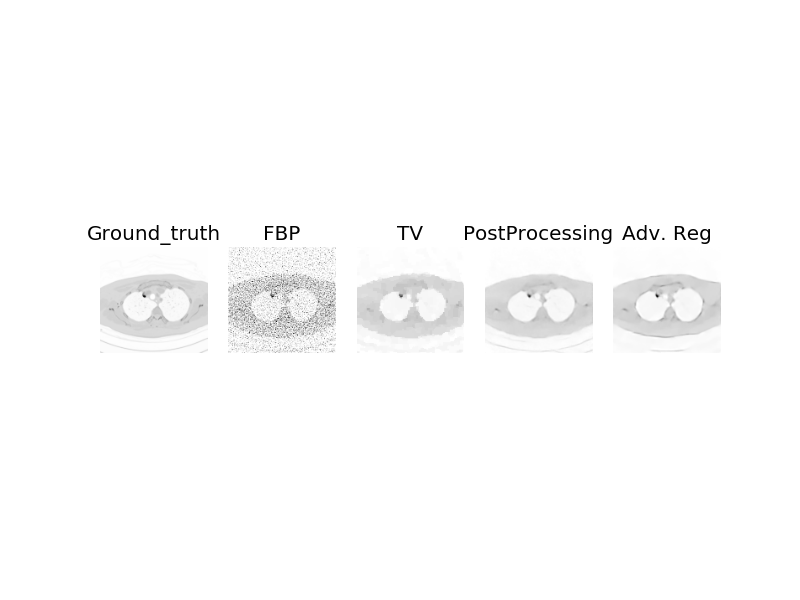} %
        \caption{Good}
        \label{sfig:good_recon}
    \end{subfigure}
    \hfill
    \begin{subfigure}[b]{0.32\textwidth} %
        \centering
        \includegraphics[trim={0cm 0cm 0cm 0cm},clip, height=\commonFigHeight]{Sinogram1_turned_scaled.png} %
        \caption{Sinogram}
        \label{sfig:sinogram_data}
    \end{subfigure}
    \hfill
    \begin{subfigure}[b]{0.32\textwidth} %
        \centering
        \includegraphics[trim={5.5cm 6cm 11.5cm 6.5cm},clip, height=\commonFigHeight]{fe_ct_4.png} %
        \caption{Bad}
        \label{sfig:bad_recon}
    \end{subfigure}

    \caption{Comparison of CT reconstructions: (\subref{sfig:good_recon}) a good quality reconstruction, (\subref{sfig:sinogram_data}) the corresponding sinogram data, and (\subref{sfig:bad_recon}) a poor quality reconstruction.}
    \label{fig:good_sinogram_bad_example} %
\end{figure}

In this section we will present the method for learning the regularizer
This section delves into the concept of \textbf{adversarial regularization} for learning priors in inverse problems. The core idea is to train the regularizer $R$ such that it assigns low values to samples from a target distribution of ``good'' images, denoted $\mathbb{P}_U$, and high values to samples from a distribution of ``bad'' or undesirable images, $\mathbb{P}_n$. \Cref{fig:good_sinogram_bad_example} conceptually illustrates this distinction, where clean ground truth images are ``good'', while noisy corruptions are ``bad''.

\noindent Inspired by the Wasserstein GAN framework \citep{arjovsky2017wasserstein}, the 1-Wasserstein distance between the clean and noisy distributions is employed as a weakly supervised loss for the regulariser. The 1-Wasserstein distance between $\mathbb{P}_n$ and $\mathbb{P}_U$ is given by:
$$
\operatorname{Wass}_1\left(\mathbb{P}_n, \mathbb{P}_U\right)=\underset{R \in 1\text{-Lip}}{\sup} \mathbb{E}_{U \sim \mathbb{P}_n}\left[R(U)\right]-\mathbb{E}_{U \sim \mathbb{P}_U}\left[R(U)\right],
$$
where the supremum is taken over all 1-Lipschitz functions $R$. By finding an appropriate $R$, this formulation allows to train the regulariser to effectively capture image statistics without requiring paired examples of ``good'' and ``bad'' images. In practice, the regulariser is parameterized as:
$$
R_{\Theta}(u)=\Psi_{\Theta}(u)+\rho_0|u|_2^2,
$$
where $\Psi_{\Theta}(u)$ is a (potentially convex \citep{mukherjee2024data} or weakly convex \citep{shumaylov2024weakly}) convolutional neural network (CNN) and an additional $l_2$ regularization term that enhances analysis and ensures coercivity. Ensuring exact 1-Lipschitzness turns out to be rather complicated, and the network is trained by minimizing the following loss function \citep{lunz2018adversarial}:
$$
\min_ {\Theta}\; \mathbb{E}_{U \sim \mathbb{P}_U}\left[\Psi_{\Theta}(U)\right]-\mathbb{E}_{U \sim \mathbb{P}_n}\left[\Psi_{\Theta}(U)\right]+\mu \cdot \mathbb{E}\left[\left(\left\|\nabla_u \Psi_{\Theta}(U)\right\|-1\right)_{+}^2\right]
$$
The first two terms encourage the regularizer to distinguish between ``good'' and ``bad'' images, while the third term softly enforces a Lipschitz constraint on the regularizer.
\paragraph{Analysis.}
Once trained (with parameters $\Theta^*$), the learned (convex/weakly-convex) adversarial regulariser (AR/ACR/AWCR) is incorporated into a variational problem:
$$
\underset{u}{\arg \min }\; \frac12\|A u-y\|_2^2+\alpha\left(\Psi_{\Theta^*}(u)+\rho_0\|u\|_2^2\right)
$$
This variational problem can be solved using standard optimization techniques such as (sub)gradient methods, or proximal descent. The resulting model benefits from the theoretical properties of variational methods, including well-posedness and stability \citep{shumaylov2023provably,lunz2018adversarial,mukherjee2024data,shumaylov2024weakly}.

A theoretical justification for using such a loss can be seen by analyzing the gradient descent flow of the regularizer with respect to the Wasserstein distance. Under certain assumptions, it can be shown that this flow leads to the fastest decrease in Wasserstein distance for any regularization functional with normalized gradients.
\begin{theorem}[Theorem 1 in \citet{lunz2018adversarial}]
Let $$
\begin{gathered}
g_\eta(u):=u-\eta \cdot \nabla_u \Psi_{\Theta^*}(u) . \\
\mathbb{P}_\eta:=\left(g_\eta\right) \# \mathbb{P}_n
\end{gathered}
$$ and assume that $\eta \mapsto \operatorname{Wass}\left(\mathbb{P}_u, \mathbb{P}_\eta\right)$ admits a left and a right derivative at $\eta=0$, and that they are equal.
Then,
$$
\left.\frac{\mathrm{d}}{\mathrm{d} \eta} \operatorname{Wass}\left(\mathbb{P}_u, \mathbb{P}_\eta\right)\right|_{\eta=0}=-\mathcal{E}_{U \sim \mathbb{P}_n}\left[\left\|\nabla_u \Psi_{\Theta^*}(U)\right\|^2\right]=-1 .
$$
This is the fastest decrease in Wasserstein distance for any regularization functional with normalized gradients!
\end{theorem}
Furthermore, under a data manifold assumption (DMA) and a low-noise assumption, the distance function to the data manifold is a maximizer of the Wasserstein loss. This provides further intuition for the effectiveness of adversarial training in learning a regularizer that captures the underlying data distribution.

\begin{assumption}[Data Manifold Assumption (DMA)]
   The measure $\mathbb{P}_u$ is supported on the weakly compact set $\mathcal{M}$, i.e. $\mathbb{P}_u\left(\mathcal{M}^c\right)=0$. 
\end{assumption}
\noindent Denote by $P_{\mathcal{M}}: X \rightarrow \mathcal{M}, u \mapsto \arg \min _{v \in \mathcal{M}}\|u-v\|$ the projection onto the data manifold.
\begin{assumption}[Low Noise Assumption (LNA)] The pushforward of the noisy distribution under the projection equals the clean distribution,
$\left(P_{\mathcal{M}}\right) _\#\left(\mathbb{P}_n\right)=\mathbb{P}_r$. This corresponds to an assumption that the noise level is low in comparison to manifold curvature.
\end{assumption}
\begin{theorem}
Assume DMA and LNA. Then, the distance function to the data manifold
$$
u \mapsto \min _{v \in \mathcal{M}}\|u-v\|_2
$$
is a maximizer to the Wasserstein Loss
\begin{equation}\label{eq:wasloss}
\sup _{R \in 1\text{-Lip}} \mathbb{E}_{U \sim \mathbb{P}_n} R(U)-\mathbb{E}_{U \sim \mathbb{P}_r} R(U)    
\end{equation}
\end{theorem}
\begin{remark}
The functional in \Cref{eq:wasloss} does not necessarily have a unique maximizer. However, in certain settings it can be shown to be unique almost everywhere, see \citet{staudt2022uniqueness,milne2022new}.
\end{remark}
\noindent It is worth mentioning that while theoretically appealing, recent work \citep{stanczuk2021wasserstein} has shown that the practical success of Wasserstein GANs may not be solely attributed to their ability to approximate the Wasserstein distance.
\subsubsection{Extensions}
The development of learned regularizers, particularly adversarial ones, is an active research area with several important possible extensions:
\begin{itemize}[nosep]
    \item \textbf{Generalization:} The generalization capabilities of machine-learned regularisers to out-of-distribution data remain a subject of ongoing empirical investigation, as explored in, e.g. \citep{lunz2022machine}.
    \item \textbf{Stronger Constraints:} Imposing stronger structural constraints on the learned regulariser, such as those related to source conditions or optimization landscapes, can lead to improved theoretical guarantees \citep{mukherjee2021learning} and practical implementations \citep{shumaylov2023provably,shumaylov2024weakly}, but these need to be designed with a problem in mind \citep{shumaylov2025lie}. %
    \item \textbf{Qualitative Properties:} Ensuring that learned regularisers fulfill desired qualitative properties, such as invariance or equivariance to certain transformations (e.g., affine transformations), can be achieved by designing specialized network architectures like equivariant NNs \citep{celledoni2021structure, celledoni2021equivariant}. %
    \item \textbf{Choice of Optimality Criteria:} The definition of optimality for a regulariser should be task-dependent. Task-adapted inversion strategies aim to learn reconstructions that are optimal for a specific end-goal or metric \citep{escudero2023integrating, adler2022task}. %
    \item \textbf{Scalability:} Training regularisers for large-scale and high-dimensional inverse problems could necessitate efficient network architectures, such as invertible networks, to manage computational complexity \citep{etmann2020iunets}. %
    \item \textbf{Uncertainty Quantification:} Learned convex regularisers could be integrated into frameworks for uncertainty quantification, for example, using proximal Markov Chain Monte Carlo (MCMC) methods \citep{pereyra2016proximal}.
\end{itemize}
\section{Plug-and-Play (PnP) Methods}
This section explores Plug-and-Play Prior (PnP) methods, a class of iterative reconstruction algorithms that utilize black-box denoisers, which can range from traditional algorithms to powerful deep learning-based denoisers. The core idea behind PnP is to consider operator splitting techniques for optimizing the variational objective in \Cref{eq:ip_main}, decoupling the regularization step from the data fidelity term, and replacing regularization steps with sophisticated denoisers.
PnP methods are rooted in operator splitting techniques, such as the Alternating Direction Method of Multipliers (ADMM) \citep{setzer2011operator}. Consider a reformulation of \Cref{eq:ip_main}, introducing an auxiliary variable $v$:
$$
\min _{u, v}\{\mathcal{D}(A u, y)+\alpha R(v)\} \quad \text { s.t. } u=v .
$$
The augmented Lagrangian associated with this constrained problem is:
$$
L_\lambda(u, v, h)=\mathcal{D}(A u, y)+\alpha R(v)+\frac{\lambda}{2}\|u-v+h\|_2^2-\frac{\lambda}{2}\|h\|_2^2,
$$
where $h$ is the Lagrange dual variable and $\lambda>0$ a penalty parameter. ADMM consists of approximating a solution to the saddle point problem for $L_\lambda$ by iterating
\begin{align*}
u^{k+1} & =\underset{u}{\arg \min } \;L_\lambda(u,v^{k},h^{k}),\\
v^{k+1} & = \underset{v}{\arg \min } \;L_\lambda(u^{k+1},v^k,h^{k}),\\
h^{k+1} & = h^{k} + (u^{k+1} - v^{k+1}).
\end{align*}
In particular, the updates in $u$ and $v$ read
\begin{align}
u^{k+1} & = \underset{u}{\arg \min } \;\mathcal{D}(Au,y) + \frac{\lambda}{2} \|u-v^{k} + u^{k}\|_2^2, \label{eq:pnp1}\\
v^{k+1} & = \operatorname{prox}_{\tau \frac{\alpha}{\lambda} R}\left(u^{k+1}+h^k \right).\label{eq:pnp2}
\end{align}
\noindent The crucial insight for PnP methods is that this decouples the measurement fidelity step from the reconstruction regularization, done by denoising $v$. The (regularizing) $v$-update step in \Cref{eq:pnp2} can be recognized as the proximal operator of the regularizer $R$ scaled by $\alpha/\lambda$. This allows for the replacement of the proximal operator with any effective denoising algorithm $D$, such as BM3D \citep{dabov2009bm3d}, non-local means \citep{buades2005non}, or deep learning-based denoisers. This flexibility gives rise to the name ``Plug-and-Play''.

\subsection{Theoretical Properties}
Despite their widespread empirical success and practical utility, PnP methods present several theoretical challenges that have been the subject of ongoing research. A primary concern is that, due to the black-box nature of the denoiser $D$, it is often unclear to what objective function, if any, the PnP iterations converge, or even if they converge at all without restrictive assumptions. What is worse however, is that even under significant restrictions, it remains unclear what properties the limit points satisfy, unless there exists a corresponding regularizer.

\begin{itemize}[nosep]
    \item \textbf{Convergence}: In general, PnP is not provably convergent. 
Early results typically required strong conditions, such as non-expansive denoisers, or the data-fidelity term $\mathcal{D}(Au,y)$ being strongly convex in $u$. Without access to the explicit form of the regularizer (implicitly defined by the denoiser), characterizing the fixed points of the iteration or the properties of the limit can be challenging.

\noindent More recent variants, like the gradient step denoisers have shown promise in achieving convergence under milder conditions \citep{hurault2021gradient}. These approaches often involve specific parameterizations or interpretations of the denoising operator.

    \item \textbf{Regularization}: A fundamental characteristic of many PnP schemes is the lack of an explicit representation for the regularizer $R(u)$ that the denoiser $D$ implicitly implements. This poses difficulties for a direct Bayesian interpretation, where the regularizer would typically correspond to a prior probability distribution. Without an explicit $R(u)$, it is hard to ascertain the precise nature of the prior being enforced or to analyze its properties.
    
    Some works exist, e.g. \textbf{Regularization by Denoising (RED)}, proposed by \citet{romano2017little}. RED defines an explicit variational regularizer based on a given denoiser ${D}(\cdot)$:
$$
R_{RED}(u) = \frac{1}{2} \langle u, u - {D}(u) \rangle. %
$$
The iterative schemes based on this $R_{RED}$ can be shown to seek stationary points of an explicit objective function $\frac{1}{2}\|Au-y\|^2 + \alpha R_{RED}(u)$. 
However, a critical and restrictive condition for this gradient formulation to hold is that the Jacobian of the denoiser ${D}(u)$ must be symmetric \citep{reehorst2018regularization}, a property not generally satisfied by many advanced denoisers, especially deep neural networks.

   \item Gradient-step (GS) denoisers \citep{hurault2021gradient} model the denoiser ${D}_{\theta}$ as an explicit gradient step of a potential function $R_{\theta}$:
$$
{D}_{\theta}(u) = u - \nabla R_{\theta}(u).
$$
A common choice for the potential is $R_{\theta}(u) = \frac{1}{2}\|u - \Psi_{\theta}(u)\|_2^2$, where $\Psi_{\theta}$ is any differentiable (deep) network. This formulation similarly directly provides an explicit representation of the regulariser $R_{\theta}(u)$ being enforced. What is more, the denoiser in this case is exactly a proximal operator $ D_\theta(x)=\operatorname{prox}_{\phi_\theta}(x)$, where $\phi_\theta$ is defined by
$$
\phi_\sigma(x)=R_\theta\left(D_\theta^{-1}(x)\right)-\frac{1}{2}\left\|D_\theta^{-1}(x)-x\right\|^2,
$$
admitting various desirable properties like Lipshitz smoothness and weak convexity \citep{hurault2022proximal,tan2024provably}, often resulting in provable convergence.  
\end{itemize}

\subsection{In-depth: Linear Denoiser Plug and Play}
A critical aspect of robust inverse problem solving is ensuring that the chosen regularization strategy is \textit{convergent} as in \Cref{thm:stab}. 
While PnP schemes using learned denoisers can achieve convergence of iterates under certain conditions, explicitly controlling the regularization strength to ensure \textit{convergent regularization} (i.e., convergence to a true solution as data noise $\delta \to 0$) is crucial. This subsection details a principled approach for such control when the denoiser is linear, based on spectral filtering, as introduced by \citet{hauptmann2024convergent}. Insights from the linear case may also inform strategies for nonlinear denoisers \citep{khelifa2025enhanced}.

A significant challenge in this case, particularly when utilizing learned denoisers, is the adjustment of regularization strength. These denoisers are often trained for a specific noise level $\sigma$, yet the effective noise within PnP iterations can vary, and the overall regularization must be adapted to the noise present in the measurements $y^\delta$. Empirically, this has been approached by tuning regularization strength by denoiser scaling  \citep{xu2020boostingperformanceplugandplaypriors}.

\begin{figure}[!ht]
    \centering
    \includegraphics[width=0.5\linewidth]{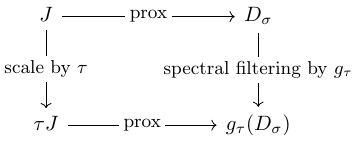} %
    \caption{Diagram illustrating the concept of spectral filtering. From \citep{hauptmann2024convergent}.} %
    \label{fig:spectral_filtering_pnp}
\end{figure}

\begin{figure}[!htb] %
    \centering
    \begin{subfigure}[b]{.48\textwidth} 
        \centering
        \includegraphics[height=0.67\linewidth]{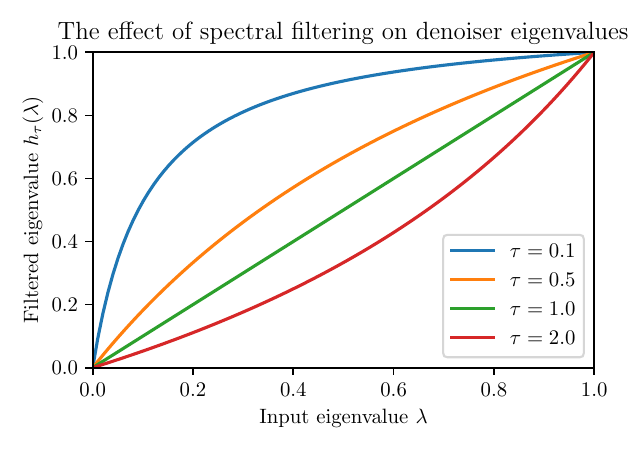} %
        \caption{Eigenvalue spectral filtering.} 
        \label{sfig:eigenvals_filtering_combined}
    \end{subfigure}
    \vspace{0.5cm}
    \hfill
    \begin{subfigure}[b]{.48\textwidth}
        \centering
        \includegraphics[height=0.67\linewidth]{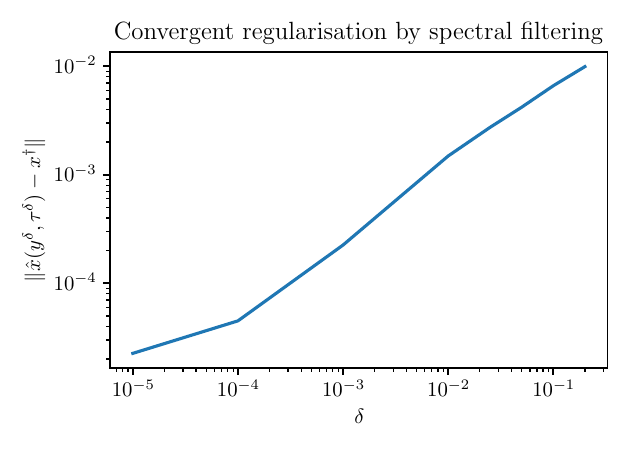} %
        \caption{Illustration of convergent regularization in practice.}
        \label{sfig:spectral_convergence_combined}
    \end{subfigure}
    
    \begin{subfigure}[b]{.98\textwidth}
        \centering
        \includegraphics[width=0.9\linewidth]{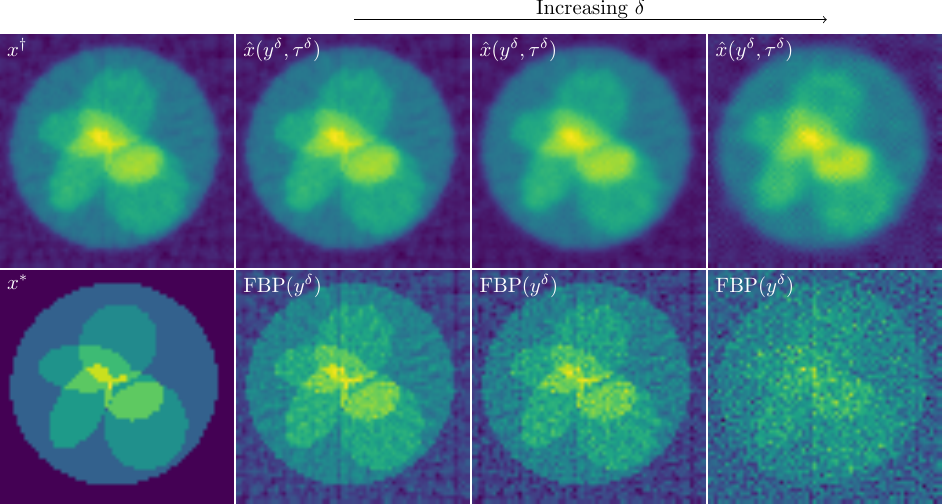} %
        \caption{Illustration of convergent regularization via a selection of snapshots from the plot in \subref{sfig:spectral_convergence_combined}.}
        \label{sfig:spectral_images_combined}
    \end{subfigure}

    \caption[Spectral Filtering for Convergent Regularisation]{Further concepts in spectral filtering and application to CT reconstruction.
    (\subref{sfig:eigenvals_filtering_combined}) Eigenvalue spectral filtering.
    (\subref{sfig:spectral_convergence_combined}) Spectral filtering to control regularisation strength for convergent regularisation.
    (\subref{sfig:spectral_images_combined}) Resulting images from the CT reconstruction.
    The linear denoiser filter is $g_\tau(\lambda) = \lambda / (\tau - \lambda (\tau - 1))$. All illustrations from \citep{hauptmann2024convergent}.}
    \label{fig:spectral_filtering_combined_main} 
\end{figure}

Consider a linear denoiser $D_\sigma: X \rightarrow X$. For $D_\sigma$ to be the proximal operator $\operatorname{prox}_J$ of some convex functional $J: X \rightarrow \mathbb{R} \cup \{\infty\}$, $D_\sigma$ must satisfy specific conditions: it must be symmetric and positive semi-definite \citep{moreau1965proximite, gribonval2020characterization}. For the resulting functional to be convex as well, $D_\sigma$ must be non-expansive, i.e. we will assume that its eigenvalues live in the interval $[0, 1]$ for the sake of contractivity. Lastly, we will assume that the operator norm is bounded from below, such that the inverse is well defined and is a bounded operator. If these conditions hold, the functional $J$ is uniquely determined by $D_\sigma$ (up to an additive constant). The objective then becomes controlling the regularization strength by effectively scaling this underlying functional $J$. The difficulty is that one typically only has access to the denoiser $D_\sigma$, not $J$ itself. However, when the denoiser is linear, it turns out to be possible to appropriately modify the denoiser based on the following observations. By definition of a proximal operator
$$
D_\sigma=\operatorname{prox}_J=(\mathrm{id}+\partial J)^{-1}.
$$
On the other hand, since $D_\sigma$ is linear, $D_\sigma^{-1}$ is linear, and by above $\partial J=: W$ is also linear. As a result, $J(x)=\frac{1}{2}\langle x, W x\rangle$ up to an additive constant. Inverting, we have $W=D_\sigma^{-1}-\operatorname{Id}$. Therefore,
$$
J(x)=\frac{1}{2}\left\langle x,\left(D_\sigma^{-1}-\operatorname{Id}\right) x\right\rangle.
$$
We can control the regularization strength, by scaling the regularization functional $J$:
$$
\tau J(x)=\frac{1}{2}\left\langle x,\left(\left[\tau D_\sigma^{-1}-(\tau-1) \operatorname{Id}\right]-\operatorname{Id}\right) x\right\rangle,
$$
resulting in 
$$
\operatorname{prox}_{\tau J}=\left(\tau D_\sigma^{-1}-(\tau-1) \operatorname{Id}\right)^{-1}=g_\tau\left(D_\sigma\right) .
$$
Here $g_\tau: \mathbb{R} \rightarrow \mathbb{R}$, given by $g_\tau(\lambda)=\lambda /(\tau-\lambda(\tau-1))$ is applied to $D_\sigma$ using the functional calculus. This implies that applying this filter (as illustrated conceptually in \Cref{fig:spectral_filtering_pnp}) effectively transforms the original denoiser $D_\sigma = \operatorname{prox}_J$ into $g_\tau(D_\sigma) = \operatorname{prox}_{\tau J}$. Furthermore \Cref{sfig:eigenvals_filtering_combined} illustrates the effect on eigenvalues of the resulting linear denoiser as a function of original eigenvalues.

\noindent This approach differs from traditional spectral regularization methods (e.g., Tikhonov regularization, Landweber iteration) \citep{engl1996regularization}, where filtering is typically applied to the singular values of the forward operator $A$. In contrast, here the denoiser (and thus the implicit prior) is modified.

\noindent It turns out to be possible to show convergent regularization in general for other spectral filters satisfying technical conditions that $\left(1-g_\tau(\lambda)\right) /\left(\tau g_\tau(\lambda)\right)$ is bounded above and below by positive values and converges as $\tau \rightarrow 0$.
Under such conditions one achieves convergent PnP regularization by Theorem 5 in \citet{hauptmann2024convergent}, and an example of this is illustrated in \Cref{sfig:spectral_convergence_combined,sfig:spectral_images_combined}.

\section{Outlooks}

The preceding sections have largely focused on methodologies that adapt existing mathematical frameworks to incorporate deep learning techniques. Examples include learned iterative schemes, where neural networks replace components of classical algorithms, and Plug-and-Play (PnP) methods, where denoisers replace splitting steps. While these approaches have demonstrated considerable empirical success, they often represent incremental adaptations rather than fundamental redesigns. Consequently, they can sometimes lack a deep theoretical grounding or may appear as ad-hoc integrations rather than solutions derived from first principles tailored to the unique characteristics of deep learning.

\noindent A central challenge and a key direction for future research is to move beyond merely ``plugging in'' deep learning components into pre-existing structures. To fully harness the potential of high-capacity, overparameterized models, the development of new frameworks that are \textit{fundamentally designed} with deep learning in mind is essential. Such a paradigm shift would likely involve concerted efforts in several interconnected areas:

\begin{itemize}[nosep]
    \item \textbf{Rethinking optimization:} Can we design optimization algorithms specifically tailored to the properties of deep neural networks, moving beyond simple gradient descent?
    \item \textbf{Embracing inductive biases:} How can we incorporate domain-specific knowledge and structure into the architecture and training of deep networks, moving beyond generic black boxes?
    \item \textbf{Developing new theoretical tools:} Can we create new mathematical tools and theoretical frameworks that can better explain the generalization capabilities of deep learning in the context of inverse problems, and provide useful guarantees guarantees on stability and convergence for learned solution maps?
\end{itemize}

A persistent theme in this endeavor is navigating the trade-off between capacity and guarantees.  Constraints are necessary for interpretability and reliability, but excessive constraints limit the expressive power of deep learning. Finding the sweet spot is key. This might involve:
\begin{itemize}[nosep]
    \item \textbf{Developing more flexible constraints:} Can we design constraints that are less restrictive but still ensure desirable properties like stability and convergence?
    \item  \textbf{Learning constraints from data:} Can we use data to learn optimal constraints that balance capacity and guarantees?
\end{itemize}
Ultimately, the future of data-driven inverse problems lies in integrating deep learning more deeply with the underlying mathematical principles. This will require both theoretical and practical innovations, but the potential rewards are immense.

%% file: chapter4.tex
\chapter{Perspectives}\label{cap4}
\section{On Task Adaptation}
\begin{figure}[!htb]
\centering
\includegraphics[width=.8\textwidth]{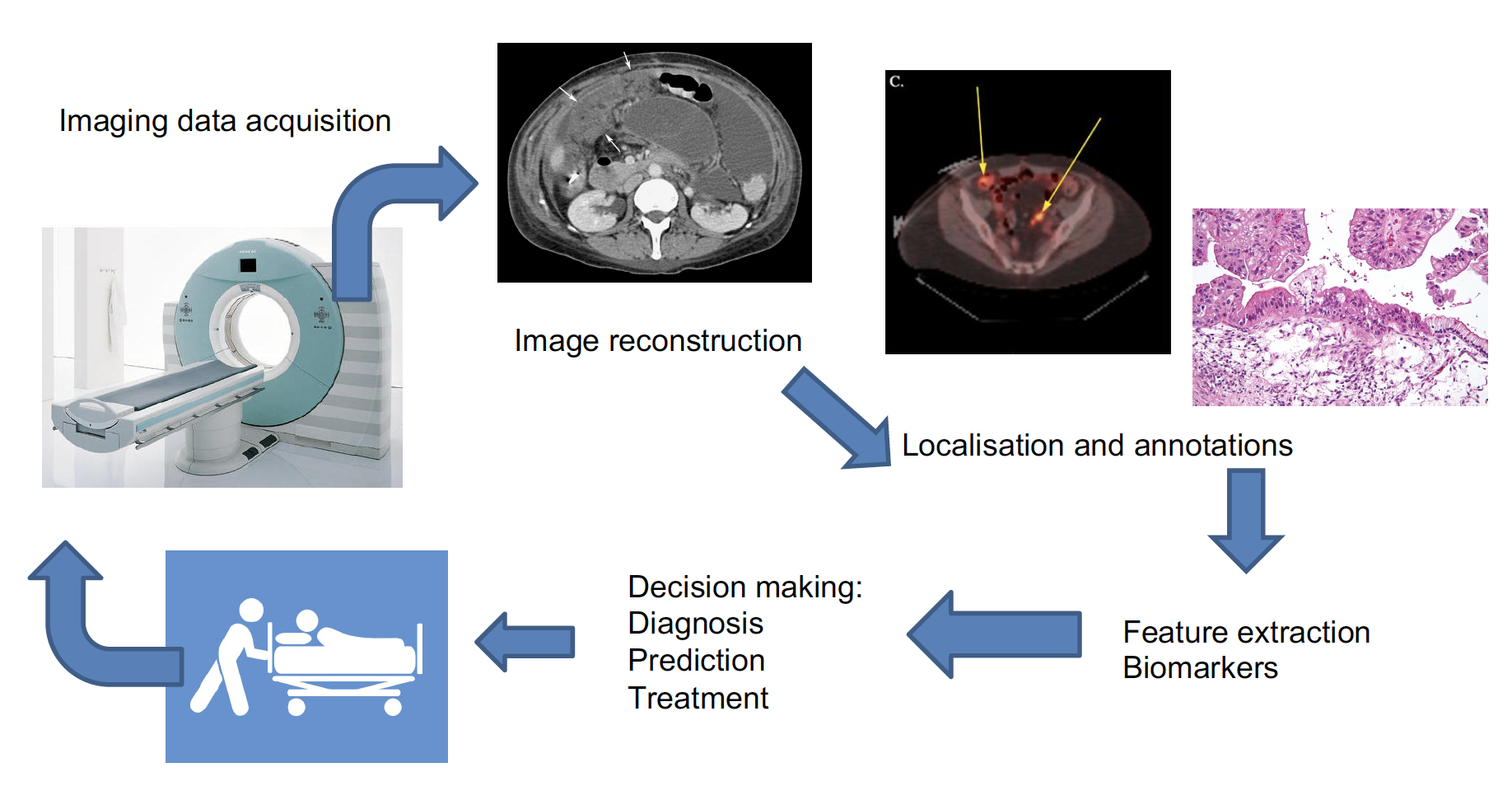}
    \caption{\textbf{Biomedical imaging pathway}: 
The path from imaging data acquisition to prediction; diagnosis; treatment planning, features several processing and analysis steps which usually are performed {\bf sequentially}. CT data and segmentation are courtesy of Evis Sala and Ramona Woitek.}
    \label{fig:pathway}
\end{figure}

The field of inverse problems is increasingly benefiting from the integration of deep learning methodologies. This chapter explores the evolving perspective that moves beyond addressing inverse problems as isolated reconstruction tasks, instead considering their embedding within broader, interconnected workflows and multi-tasking scenarios. Traditionally, inverse problems that can be modeled mathematically comprise only a part of the overall problem. Oftentimes, the full problem involves a sequential pipeline of individual tasks, such as data acquisition, reconstruction, segmentation, and classification, see \Cref{fig:pathway} for example. These stages, while often tackled sequentially and independently, are inherently intertwined. The output quality and characteristics of one stage directly influence the performance and feasibility of subsequent ones. Treating them in isolation can therefore lead to suboptimal overall performance.

\noindent A key observation driving current research is that the metrics used to evaluate the quality of a reconstruction should be intrinsically linked to the ultimate objective of the entire workflow. For instance, in clinical medical imaging, the primary goal is rarely the reconstruction of a visually appealing image, but rather to enable accurate diagnosis, guide treatment planning, or monitor therapeutic response. This motivates the concept of \textbf{task-adapted reconstruction} \citep{adler2022task,wu2018end}, wherein the reconstruction process is explicitly tailored to optimize performance on a specific downstream task, such as segmentation, classification, or quantitative parameter estimation.

\noindent Within the classical, purely model-driven (or knowledge-driven) paradigm of \Cref{cap2}, designing such task-adapted reconstruction algorithms can rapidly become intractable due to the complexity of formulating and solving the coupled optimization problems. Deep learning, however, offers a powerful and flexible framework for realizing task-adapted reconstruction. 

\noindent Training neural networks often involves highly non-convex optimization, and joint training for multiple tasks does not fundamentally alter this characteristic. In fact, combining tasks within a unified learning pipeline can create opportunities for synergy, where the optimization process for one task can provide beneficial regularization or feature representations for another. The inherent non-convexity of end-to-end learned systems means that extending them to incorporate downstream tasks does not necessarily introduce greater optimization challenges than those already present in learning the reconstruction alone.

\noindent Consider, for example, the task of detecting ovarian cancer from medical images (radiomics). This process typically involves reconstructing an image from sensor measurements, followed by segmentation of potential tumorous regions, and then extraction of quantitative imaging features for statistical analysis and classification. By jointly optimizing the reconstruction and segmentation processes within a single deep learning model, it is possible to guide the reconstruction to produce images that are not only faithful to the measured data but are also more amenable to accurate segmentation by the learned segmentation module. 

\begin{example}[Joint Reconstruction Segmentation \citet{adler2022task}]
For example, in tomographic reconstruction, we can jointly optimize the reconstruction and segmentation processes using a combined loss function. Consider a CNN-based MRI reconstruction ($X$) and CNN-based MRI segmentation ($D$).
\begin{figure}[h!]
    \centering
    \begin{subfigure}{\textwidth} %
        \centering
        \includegraphics[width=0.75\textwidth]{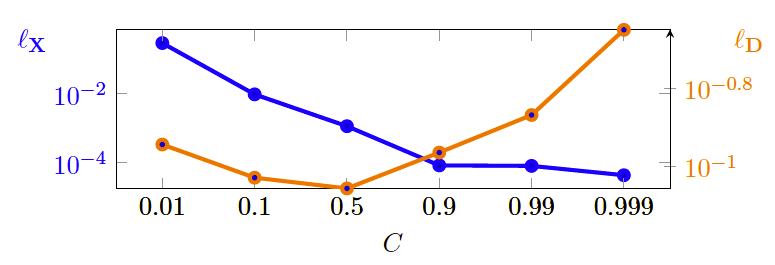} %
        \caption{Minimal loss values for various $C$ values, showcasing that jointly training for reconstruction and segmentation is better.}
        \label{sfig:losses}
    \end{subfigure}
    
    \vspace{1em} %
    
    \begin{subfigure}{\textwidth} %
        \centering
        \includegraphics[width=\textwidth]{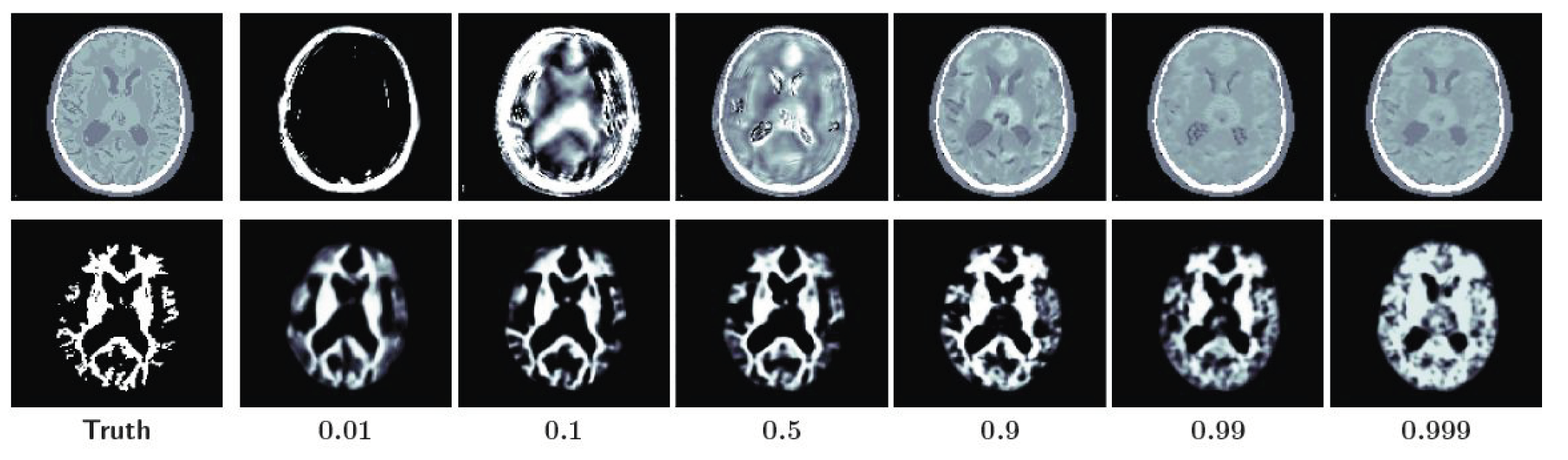} %
        \caption{CNN-based reconstructions (top row) and segmentations (bottom row).}
        \label{sfig:cnn_segmentations}
    \end{subfigure}
    
    \caption{Task-adapted reconstruction, with CNN-based MRI reconstruction (task $X$) and CNN-based MRI segmentation (task $D$). Both are trained jointly with combined loss {$C\ell_X+(1-C)\ell_D$ for varying $C\in[0,1]$}. All figures from \citep{adler2022task}.}
    \label{fig:task_adapted_reconstruction} %
\end{figure}
$$
\begin{gathered}
\left(\theta^*, \vartheta^*\right) \underset{(\theta, \vartheta) \in \Theta \times \Xi}{\arg \min }\left\{\frac{1}{m} \sum_{i=1}^m \ell_{\text {joint }}\left(\left(x_i, \tau\left(z_i\right)\right),\left(\mathcal{A}_\theta^{\dagger}\left(\mathrm{y}_i\right), \mathcal{T}_{\vartheta} \circ \mathcal{A}_\theta^{\dagger}\left(\mathrm{y}_i\right)\right)\right)\right\} \\
\ell_{\text {joint }}\left((x, d),\left(x^{\prime}, d^{\prime}\right)\right):=(1-C) \ell_X\left(x, x^{\prime}\right)+C \ell_D\left(d, d^{\prime}\right) \quad \text { for fixed } C \in[0,1] .
\end{gathered}
$$

\noindent This loss function balances the reconstruction error ($\ell_X$) and the segmentation error ($\ell_D$), allowing for a trade-off between the two tasks. \Cref{sfig:losses} illustrates that whenever segmentation performance is the ultimate task to be solved, training primarily for reconstruction results in poor segmentations, while training primarily for segmentation results in poor reconstructions. What is of most significance is that training for both jointly with equal weighting actually results in better segmentations than if one were to train for segmentation only! 
\end{example}
\noindent This co-adaptation can lead to improved accuracy and robustness in both the reconstruction and the segmentation, ultimately enhancing the reliability of the radiomic analysis and the diagnostic outcome. The degrees of freedom inherent in solving an ill-posed inverse problem can be strategically utilized to favor solutions that, while consistent with the data, also possess features conducive to the success of the subsequent task.

\noindent This task-adapted approach not only enhances the efficiency and accuracy of medical image analysis but also has broader implications for healthcare accessibility. By automating or streamlining certain tasks, such as segmentation, it reduces the reliance on specialized expertise, potentially making advanced imaging techniques more widely available in settings with limited resources.   

\section{The Data Driven - Knowledge Informed Paradigm}

The recent trajectory of research in inverse problems, as explored throughout these notes, signifies a notable paradigm shift. The ascent of deep learning has significantly transformed the field, primarily through its capacity to learn efficient data-driven priors, often surpassing traditional handcrafted models. While empowering, this shift introduces significant challenges, most critically by a potential decline in interpretability especially as model complexity increases, e.g. for multimodal applications. Open questions remain regarding the generalization of these models across diverse datasets and the crucial balance between empirical performance and robust theoretical guarantees.

\noindent It is important to acknowledge the limitations of classical, purely model-based approaches and recognize the advancements that deep learning has brought to the field. However, we must also emphasize the necessity of grounding these powerful data-driven techniques within rigorous mathematical frameworks. Such integration is not merely an academic exercise but a necessary step to ensure interpretability, and provably provide assurances beyond empirical validation, to ultimately foster trust in safety-critical applications.

\paragraph{The Imperative for Guarantees via Structured Learning.}
A central challenge in this new paradigm lies in reconciling the expressive power of deep learning with the need for verifiable guarantees. In many inverse problems, particularly within medical imaging, the concept of an absolute ``ground truth'' does not exist. This ambiguity elevates the importance of model interpretability and reliability. Consequently, while the allure of purely data-driven solutions is strong, a wholesale abandonment of mathematical formalism is untenable. The pursuit of guarantees inherently compels us to impose specific structural or functional properties onto neural network architectures. The pertinent research questions then become: What are the \textit{most effective} properties to instill for achieving meaningful guarantees (e.g., stability, robustness, fairness)? And, how can these properties be integrated into network design and training in a manner that is both principled and computationally tractable, without unduly sacrificing performance? Exploring deeper connections with established mathematical fields, such as the theory of PDEs or optimal transport, continues to be a promising avenue for discovering and formalizing such beneficial structural biases.

\paragraph{Towards more useful theoretical tools:}
The paradigm shift driven by deep learning also necessitates a corresponding evolution in our theoretical approaches. Much of the traditional analysis in inverse problems has focused on model properties, optimization landscapes, and convergence proofs, often treating the model in relative isolation from the data that was used to create it. However, deep learning models are fundamentally data-centric; their behavior, efficacy, and potential failure modes are inextricably linked to the characteristics of the training dataset.
Therefore, future analytical efforts must pivot to more explicitly account for this data dependency. It is no longer sufficient to analyze convex regularizers in abstraction. Rigorous analysis must now encompass the \textbf{training dataset} itself: its size, diversity, representativeness, potential inherent biases, and the precise manner in which these factors influence the learned model's generalization capabilities, and its robustness to distributional shifts.

\noindent As deep learning systems become more complex and their decision-making processes more opaque, \textbf{explainability} emerges as a critical concern. If a model produces a reconstruction specifically optimized for a downstream task, understanding \textit{why} the reconstruction appears as it does, and how specific features (or apparent artifacts) contribute to the downstream decision, is crucial for validation, debugging, and building trust, especially in safety-critical applications like medicine. Future research must focus on developing methods that can provide insights into these complex, end-to-end trained systems.

\paragraph{Beyond Theory:}
And alongside theory, we need to continue working on \textbf{convincing use-cases}. Ultimately, the successful integration of deep learning and mathematics holds the potential to transform medical imaging into a more accessible, efficient and widespread clinical screening tool, benefiting both clinicians and patients. The transformative potential of this research is immense, and we envision a future where its impact is widely recognized, with headlines declaring: 

\begin{displayquote}\centering \huge
    ``Deep Learning \& Maths turn CT/MRI into a clinical screening tool!''
\end{displayquote}